\documentclass[11pt,a4paper,twoside,leqno]{article}

\usepackage{times}
\usepackage{amsmath,amssymb,amsthm,amsfonts}
\usepackage{a4wide}
\usepackage{exscale}
\usepackage{enumerate}
\usepackage{latexsym}

\usepackage[all]{xy}

\newcommand{\restr}{\rvert}   
\newcommand{\abs}[1]{\left\lvert#1\right\rvert}   
\newcommand{\norm}[1]{\left\lVert#1\right\rVert}   

\newcommand{\demph}[1]{{\it #1}}

\DeclareMathOperator{\id}{Id}
\DeclareMathOperator{\ev}{ev}

\newcommand{\Cont}{{\mathcal C}} 
\DeclareMathOperator{\supp}{supp}
\DeclareMathOperator{\top2}{top}

\DeclareMathOperator{\Lin}{L} 
\DeclareMathOperator{\Leb}{L}

\DeclareMathOperator{\Ind}{Ind} 
\DeclareMathOperator{\Komp}{K} 

\DeclareMathOperator{\KTh}{K}
\DeclareMathOperator{\KK}{KK}

\DeclareMathOperator{\Cstar}{C^*}

\DeclareMathOperator{\Mult}{M}

\newcommand{\C}{\ensuremath{{\mathbb C}}}
\newcommand{\E}{\ensuremath{{\mathbb E}}}

\newcommand{\M}{\ensuremath{{\mathbb M}}}
\newcommand{\N}{\ensuremath{{\mathbb N}}}
\newcommand{\R}{\ensuremath{{\mathbb R}}}

\newcommand{\mA}{\ensuremath{{\mathcal A}}}
\newcommand{\mB}{\ensuremath{{\mathcal B}}}
\newcommand{\mF}{\ensuremath{{\mathcal F}}}
\newcommand{\mG}{\ensuremath{{\mathcal G}}}
\newcommand{\mH}{\ensuremath{{\mathcal H}}}
\newcommand{\mK}{\ensuremath{{\mathcal K}}}
\newcommand{\mL}{\ensuremath{{\mathcal L}}}

\theoremstyle{plain}

\newtheorem {theorem} {Theorem}[section]
\newtheorem {lemma}[theorem] {Lemma}
\newtheorem {proposition} [theorem]{Proposition}
\newtheorem {corollary} [theorem]{Corollary}

\newtheorem* {theorem*} {Theorem}
\newtheorem* {proposition*} {Proposition}
\newtheorem* {lemma*}{Lemma}

\theoremstyle{definition}

\newtheorem {definition} [theorem]{Definition}
\newtheorem {defprop} [theorem]{Definition and Proposition}

\newtheorem {remark} [theorem]{Remark}

\newtheorem {Xample}[theorem] {Example}

\DeclareMathOperator{\KKban}{KK^{\ban}}

\DeclareMathOperator{\Eban}{\E^{\ban}}

\DeclareMathOperator{\ban}{ban}

\newcommand{\LazyAnd}{\quad \text{and} \quad}

\newcommand{\rmd}{{\, \mathrm{d}}}
\newcommand{\unital}[1]{\widetilde{#1}} 

\newcommand{\ketbra}[2]{\big|#1\big\rangle\big\langle #2 \big|}

\DeclareMathOperator{\uEgG}{\underline{\rm E}\gG}
\DeclareMathOperator{\uEgH}{\underline{\rm E}\gH}

\newcommand{\flipX}[1]{#1^{-1}} 
\newcommand{\flipx}[1]{#1^{-1}} 
\DeclareMathOperator{\Graph}{Graph}
\DeclareMathOperator{\Morph}{Morph}
\newcommand{\gG}{\mG} 
\newcommand{\gH}{\mH} 
\newcommand{\gK}{\mK} 
\newcommand{\gL}{\mL}

\newcommand{\EbanW}[1]{\E^{\ban}_{#1}} 
\newcommand{\KKbanW}[1]{\KK^{\ban}_{#1}}

\newcommand{\ContSect}{\Gamma}
\newcommand{\loc}{{\rm \bf loc}}
\newcommand{\Linloc}{\Lin^{\!\loc}} 
\newcommand{\Komploc}{\Komp^{\!\loc}}

\newcommand{\simMorita}{\sim_{\text{M}}}

\newcommand{\MbanW}[1]{\M^{\ban}_{#1}} 
\DeclareMathOperator{\MoritaPur}{Mor}

\newcommand{\MoritabanW}[1]{\MoritaPur^{\ban}_{#1}}

\begin{document}

\title{Induction for Banach Algebras, Groupoids and $\KKban$}
\author{Walther Paravicini}
\date{December 17, 2008}
\maketitle

\begin{abstract}
\noindent Given two equivalent locally compact Hausdorff groupoids, the Bost conjecture with Banach algebra coefficients is true for one if and only if it is true for the other. This also holds for the Bost conjecture with C$^*$-coefficients. To show these results, the functoriality of Lafforgue's $\KK$-theory for Banach algebras and groupoids with respect to generalised morphisms of groupoids is established. It is also shown that equivalent groupoids have Morita equivalent $\Leb^1$-algebras (with Banach algebra coefficients).

\medskip

\noindent \textit{Keywords:} Locally compact groupoid, $\KKban$-theory, Banach algebra, Baum-Connes conjecture, Bost conjecture;

\medskip

\noindent \textsc{AMS 2000} \textit{Mathematics subject classification:} Primary 19K35; 46L80; 22A22; 43A20
\end{abstract}

\noindent It is a well-known fact that the C$^*$-algebras $\Cstar(\gG)$ and $\Cstar(\gH)$ are (strongly) Morita equivalent, where $\gG$ and $\gH$ are equivalent locally compact second countable Hausdorff groupoids equipped with Haar systems, see \cite{MuhReWill:87}. It is also true that the $\Leb^1$-algebras $\Leb^1(\gG)$ and $\Leb^1(\gH)$ are Morita equivalent Banach algebras in the sense of V.~Lafforgue; see Paragraph~\ref{Subsection:MoritaEquivalenceOfUnconditionalCompletions} for a definition of this notion. In the present article, we prove this fact in the course of a systematic treatment of groupoid Banach algebras and the descent, and generalise it in two ways: Firstly, we allow for more general unconditional completions $\mA(\gG)$ and $\mA(\gH)$ instead of $\Leb^1(\gG)$ and $\Leb^1(\gH)$ (see Paragraph~\ref{Subsection:MoritaEquivalenceOfUnconditionalCompletions} for details) and we consider, secondly, Banach algebra coefficients, i.e., we consider an $\gH$-Banach algebra $B$, construct the induced $\gG$-Banach algebra $\Ind_{\gH}^{\gG} B$ and show that
\[
\mA(\gH, B) \quad \simMorita \quad \mA(\gG, \Ind_{\gH}^{\gG} B).
\]
From this Morita equivalence it follows that the two Banach algebras have isomorphic $\KTh$-theory. In this article, we fit the isomorphism in $\KTh$-theory into the following commutative diagram:
\[
\xymatrix{
\KTh^{\top2,\ban}_*\left(\gH,\  B\right) \ar[rrr]^{\mu_{\mA}^B} \ar[d]_{\cong} &&& \KTh_*(\mA(\gH,\ B))\ar[d]^{\cong} \\
\KTh^{\top2,\ban}_*\left(\gG,\ \Ind_{\gH}^{\gG} B\right) \ar[rrr]^-{\mu^{\Ind\! B}_{\mA}} &&& \KTh_*(\mA(\gG,\ \Ind_{\gH}^{\gG} B)) \\
}
\]
The vertices on the left-hand side are Banach algebraic versions of C$^*$-algebraic topological $\KTh$-theory which are constructed from $\KKbanW{\gG}$ instead of $\KK_{\gG}$. We show that $\KKbanW{\gG}$ is invariant under equivalences and, more generally, that it is functorial under generalised morphisms of groupoids. We proceed in analogy to \cite{LeGall:94} where the case of $\KK$-theory for C$^*$-algebras and groupoids is treated; the additional technical challenges that arise from the fact that we deal with fields of Banach algebras are met by a systematic development of the theory, what, admittedly, adds not only accuracy but also some length to the article. As a consequence of the invariance under equivalences, we obtain the left vertical arrow in the above diagram and we can deduce that it is an isomorphism.

The horizontal arrows in the diagram are Banach algebraic versions of the Bost assembly map as defined in \cite{Lafforgue:06}. The commutativity of the diagram is shown by an analysis of how equivalences of groupoids and the descent homomorphisms for unconditional completions interact.

The main result of this article can now be read off the above diagram: The Bost conjecture with Banach algebra coefficients for $\gG$ is equivalent to the Bost conjecture with Banach algebra coefficients for $\gH$, because induction is an equivalence of categories between the $\gG$-Banach algebras and the $\gH$-Banach algebras.

In addition to this, there seem to be very few results for the Bost conjecture with Banach algebra coefficients as presented in this article, although this conjecture is a rather obvious variant of the Bost conjecture with C$^*$-coefficients. To my knowledge, there is only the Green-Julg theorem for compact groups, its generalisation to proper groupoids, and, as a consequence, the split surjectivity of the Bost map for proper Banach algebras as coefficients (see \cite{Paravicini:07:GreenJulg}). It might be the case that, so far, the left-hand side of the Bost conjecture with Banach algebra coefficients is not understood well enough.

But for C$^*$-algebra coefficients, we can obtain additional results. First of all, there is a diagram for C$^*$-algebras and $\KTh^{\top2}$ instead of $\KTh^{\top2,\ban}$, analogous to the one above. It can in fact be obtained from the above diagram because the Bost assembly map for C$^*$-algebras factors through $\KTh^{\top2,\ban}$. We can conclude that the Bost conjecture with C$^*$-algebra coefficients for $\gG$ is equivalent to the Bost conjecture with C$^*$-algebra coefficients for $\gH$.

Let us sketch an important consequence of the C$^*$-algebraic result which so far cannot be extended to Banach algebra coefficients: Let $H$ be a closed subgroup of a (second countable) locally compact Hausdorff group $G$. Then the groupoid $G \ltimes G/H$ is equivalent to $H$ (considered as a groupoid with one-point unit space). So the Bost conjecture for $H$ is equivalent to the Bost conjecture for $G\ltimes G/H$. What is left to analyse is the interplay of the Bost conjecture for $G\ltimes G/H$ and the Bost conjecture for $G$, which is much better understood for C$^*$-coefficients. If $H$ is an open subgroup of $G$, then the Bost conjecture with C$^*$-coefficients for $G$ implies the Bost conjecture with C$^*$-coefficients for $G \ltimes G/H$, so it is also true for the open subgroup $H$. This will be the subject of a forthcoming article.

\bigskip

\noindent In the first section, we review the category of locally compact Hausdorff groupoids and generalised morphisms, also introducing the so-called linking groupoid of an equivalence which will prove very useful in the subsequent sections. The second section recalls the definition of upper semi-continuous fields of Banach spaces, Banach algebras etc. Although we mainly follow \cite{Lafforgue:06}, also in Section~3, where we recall the definition of bivariant equivariant $\KTh$-theory for Banach algebras, we put some additional emphasis on compact and locally compact operators between fields of Banach spaces.

In Section~4, we show that $\KKban$ is functorial for generalised morphisms of groupoids; the basic outline of the construction follows \cite{LeGall:99}.

Section~5 first recalls from \cite{Lafforgue:06} the descent in the Banach algebraic setting: If $B$ is a $\gG$-Banach algebra for some locally compact Hausdorff groupoid $\gG$, then we construct a Banach algebra $\mA(\gG,B)$ where $\mA(\gG)$ is a so-called unconditional completion of $\Cont_c(\gG)$, e.g. $\Leb^1(\gG)$. We then show that the Banach algebra $\mA(\gG,B)$ is Morita equivalent to $\mA(\gH, \Ind_{\gG}^{\gH} B)$ for every groupoid $\gH$ equivalent to $\gG$ (see Theorem~\ref{Theorem:EquivalentGroupoidsEquivalentAlgebras} for a precise formulation of the result). We also analyse how the descent on the level of $\KKban$-theory behaves with respect to equivalences of groupoids.

In Section~6, the Bost assembly map for groupoids is presented in a completely Banach algebraic framework (in \cite{Lafforgue:06} it was introduced only for C$^*$-coefficients and not more was needed). We then show that the Bost conjecture is invariant under equivalences of groupoids.

I would like to thank Siegfried Echterhoff for his constructive and encouraging advice; the results of this article are contained in the doctoral thesis \cite{Paravicini:07} which comprises many technical details.

This research has been supported by the Deutsche Forschungsgemeinschaft (SFB 478).

\smallskip

Notation: All Banach spaces and Banach algebras that appear in this article are supposed to be complex. Let $X$, $Y$ and $Z$ be sets and $f\colon X\to Z$ and $g\colon Y\to Z$ be maps. Then we write $X \times_{f,g} Y$ for the fibre product $\{(x,y)\in X\times Y:\ f(x)=g(y)\}$. If $f$ and $g$ are understood, then we also write $X \times_Z Y$ or even $X*Y$. If $X$, $Y$ and $Z$ are topological spaces and $f$ and $g$ are continuous, then $X\times_Z Y$ is a closed subspace of $X\times Y$.

\section{Locally compact groupoids}\label{Section:Groupoids}

In this section, we recall the definitions of locally compact Hausdorff groupoids and generalised morphisms between them. Historically, the theory of locally compact groupoids was inspired by the study of foliations \cite{Connes:80, Haefliger:84, HilSkan:87} but the main application we have in mind concerns the (induction from) closed subgroups of locally compact Hausdorff groups and transformation groupoids, so we concentrate on Hausdorff groupoids with open range and source maps. The notion of a linking groupoid of an equivalence is introduced and the pullback construction for groupoids is recalled from \cite{LeGall:94, LeGall:99}. We then study how Haar systems of groupoids behave under these constructions. A general reference for Sections~\ref{Subsection:Groupoids:BasicDefinitions}, \ref{Subsection:Groupoids:Morphisms}, \ref{Subsection:EquivalencesOfGroupoids} and \ref{Subsection:Groupoids:Pullback} is Section~2.1 of \cite{LeGall:99}. Detailed proofs are contained in \cite{Paravicini:07}

\subsection{Groupoids and groupoid actions}\label{Subsection:Groupoids:BasicDefinitions}

A \demph{groupoid} is a small category such that every morphism is invertible. If $\gG$ is a groupoid, then we will denote the set of composable pairs of morphisms by $\gG^{(2)} \subseteq \gG\times \gG$ or by $\gG * \gG$, and the set of identity morphisms by $\gG^{(0)}\subseteq \gG$. The set $\gG^{(0)}$, called the unit space, can also be regarded as the set of objects of $\gG$. The range and source maps $\gG \to \gG^{(0)}$ will be denoted by $r_{\gG}$ and $s_{\gG}$ (or $r$ and $s$ if $\gG$ is understood).

Often, we will think of $\gG^{(0)}$ as being a set that is not a subset of $\gG$ but a distinct set on which the groupoid ``acts''. If $X$ is a set and $\gG$ is a groupoid such that $\gG^{(0)}=X$, then we say that $\gG$ is a groupoid \demph{over $X$.} The map that sends some $x\in X$ to the associated identity morphism in $\gG$ will usually be called $\epsilon$. In calculations, however, we will always omit this map.

Let $\gG$ be a groupoid. If $K$ and $L$ are subsets of $\gG^{(0)}$, then $\gG^L:=\{\gamma\in \gG:\ r(\gamma)\in L\}$,
$\gG_K:=\{\gamma\in \gG:\ s(\gamma)\in K\}$ and $\gG^L_K:= \gG^L \cap \gG_K$. If $g\in \gG^{(0)}$, then
$\gG_g:=\gG_{\{g\}}$ and $\gG^g:=\gG^{\{g\}}$.

A \demph{topological groupoid} $\gG$ is a groupoid which is at the same time a topological space such that the composition, inversion and the range and source maps are continuous. If $\gG$ is a groupoid over a set $X$, then we also have to assume that $X$ is a topological space and the map $\epsilon\colon X\to \gG$ is continuous.

There are several canonical examples of topological groupoids:
\begin{enumerate}
\item If $X$ is a topological space, then we define the structure of a topological groupoid on $X$ by setting $r:=s:= \id_X$ (so there are only units).
\item If $G$ is a topological group, then $G$ can be regarded as a topological groupoid if we let $r$ and $s$ be the projection on the identity element of $G$.
\item \label{Example:TopologicalGroupoidFibreProduct}

Let $X$ and $Z$ be topological spaces and let $p\colon X\to Z$ be a continuous map. We define the structure of a topological groupoid on the fibre product $X\times_Z X= X \times_p X$ by setting
\[
(X\times_Z X)^{(0)} := X \LazyAnd \epsilon\colon X \to X\times_Z X,\ x\mapsto (x,x),
\]
\[
r\colon X\times_Z X\to X,\ (y,x) \mapsto y \LazyAnd s\colon X\times_Z X\to X, \ (y,x) \mapsto x,
\]
\[
\forall x,y,z \in X, p(x)=p(y)=p(z):\ (z,y) \circ (y,x) := (z,x)\LazyAnd (y,x)^{-1} = (x,y).
\]
If $p$ is open, then $r$ and $s$ are open, too.
\end{enumerate}

\noindent Let $\gG$ and $\gH$ be topological groupoids. Then a \demph{strict morphism} $f$ from $\gG$ to $\gH$ is a continuous map from $\gG$ to $\gH$ which is also a homomorphism of groupoids (i.e., a functor). The topological groupoids, together with the strict morphisms, form a category.

A \demph{locally compact Hausdorff groupoid} is a topological groupoid such that the underlying topological space is locally compact and Hausdorff (which implies that also the unit space is locally compact and Hausdorff). To avoid technicalities, we will only consider locally compact Hausdorff groupoids in the following.\footnote{Although many constructions and results are also available for non-Hausdorff or non-locally compact groupoids; e.g.,  see \cite{Tu:04} for the non-Hausdorff case. In \cite{Lafforgue:06} and \cite{Paravicini:07}, parts of the exposition cover general topological groupoids.} Moreover, we will always assume that the range and source maps of the groupoids we consider are \emph{open}. The openness of these maps is not a dramatic restriction because our main interest is to treat the case that the groupoids carry Haar systems, and in this case, the range and source maps are automatically open.

\textsc{the rest of Section~\ref{Section:Groupoids} let $\gG$, $\gH$ and $\gK$ be locally compact Hausdorff groupoids with open range and source maps.}

We refer to \cite{LeGall:99}, Section~2.1 for the definitions of the following concepts:
\begin{itemize}
\item free, proper and principal $\gG$-spaces (which we always assume to be locally compact Hausdorff);
\item the crossed product groupoid or transformation groupoid $\Omega\rtimes \gG$ where $\Omega$ is a right $\gG$-space;
\item the quotient space $\Omega /  \gG$ where $\Omega$ is a right $\gG$-space; note that the canonical quotient map from $\Omega$ to $\Omega/\gG$ is open (\cite{Tu:04}, Lemma~2.30) and that the quotient space is locally compact Hausdorff if $\Omega$ is proper (\cite{Tu:99}, Proposition~6.3);
\end{itemize}

\noindent Let $\Omega$ be a right $\gG$-space and $\Omega'$ a left $\gG$-space. Then we define $\Omega \times_{\gG} \Omega'$ to be the quotient of $\Omega \times_{\gG^{(0)}} \Omega'$ by the diagonal action of $\gG$. If the action of $\gG$ on $\Omega$ \demph{or} $\Omega'$ is proper, then $\Omega\times_{\gG} \Omega'$ is locally compact Hausdorff.

A \emph{$\gG$-$\gH$-bimodule} or $\gG$-$\gH$-space is a locally compact Hausdorff space $\Omega$ which is at the same time a left $\gG$-space and a right $\gH$-space such that the actions commute. The anchor maps will usually called $\rho\colon \Omega\to \gG^{(0)}$ and $\sigma\colon \Omega\to \gH^{(0)}$, respectively.

Let $\Omega$ be a proper right $\gH$-space and $\Omega'$ an $\gH$-$\gK$-bimodule. Note that the quotient space $\Omega \times_{\gH} \Omega'$ of $\Omega \times_{\gH^{(0)}} \Omega'$ is locally compact and Hausdorff. On this product, there is a canonical structure of a right $\gK$-space, and if the right $\gK$-action on $\Omega'$ is proper, then so is the right action on the product $\Omega \times_{\gH} \Omega'$.

If $\Omega$ is not only a proper right $\gH$ space but also a $\gG$-$\gH$-bimodule, then we can define a canonical left $\gG$-action on the product $\Omega \times_{\gH} \Omega'$ making it a $\gG$-$\gK$-bimodule.

\subsection{Principal fibrations, graphs and morphisms}\label{Subsection:Groupoids:Morphisms}

Because we only consider groupoids with open range and source maps, we can go back to the definitions of principal fibrations and generalised morphisms of \cite{LeGall:94} instead of the more elaborate concepts of \cite{LeGall:99}. Compare also the definitions of graphs and principal fibrations in \cite{Haefliger:84} and \cite{HilSkan:87}.

\begin{definition}[Principal fibration] Let $\gH$ act on the locally compact Hausdorff space $\Omega$ on the right. A map $p$ from $\Omega$ to another topological space $X$ is called \demph{principal fibration} with structure groupoid $\gH$ if
\begin{enumerate}
\item $\Omega$ is a principal $\gH$-space, i.e., the $\gH$-action is free and proper;
\item $p$ is continuous, open and surjective;
\item $p$ is invariant under the action of $\gH$, i.e., $p(\omega) = p(\omega\eta)$ for all $(\omega,\eta) \in \Omega*\gH$.

\item $\gH$ acts transitively on each fibre of $p$, i.e., for all $\omega,\omega'\in \Omega$ such that $p(\omega) = p(\omega')$ there is an $\eta\in \gH$ such that $\omega\eta = \omega'$; note that $\eta$ is unique as $\Omega$ is free (one can even show that the map which associates $\eta$ to the pair $(\omega,\omega') \in \Omega \times_p \Omega$ is continuous).
\end{enumerate}
\end{definition}

\noindent Because $p$ is invariant under the action of $\gH$, it induces a continuous open map $\tilde{p} \colon \Omega/\gH \to X$. Because $\gH$ acts transitively on each fibre, $\tilde{p}$ is injective and hence a homeomorphism.

A generalised morphism of locally compact Hausdorff groupoids is an isomorphism class of graphs: A \demph{graph} $\Omega$ from $\gG$ to $\gH$ is a $\gG$-$\gH$-bimodule (with anchor maps $\rho$ and $\sigma$, say), such that $\rho\colon \Omega \to \gG^{(0)}$ is a principal fibration with structure groupoid $\gH$. Two graphs $\Omega$ and $\Omega'$ from $\gG$ to $\gH$ are called \demph{equivalent} if there is a homeomorphism from $\Omega$ to $\Omega'$ which intertwines the anchor maps and the actions of $\gG$ and $\gH$, i.e., an isomorphism of $\gG$-$\gH$-bimodules.

\begin{definition}[(Generalised) morphism]  A \demph{(generalised) morphism} from $\gG$ to $\gH$ is simply an equivalence class of graphs. If $\Omega$ is a graph, then we denote the corresponding morphism by $[\Omega]$.
\end{definition}

Let $f\colon \gG \to \gH$ be a strict morphism of groupoids. Then we define $\Graph(f)$ to be the following graph from $\gG$ to $\gH$:
\[
\Graph(f):=\Omega:=\gG^{(0)} \times_{\gH^{(0)}} \gH,
\]
where the fibre product is taken over the maps $f\restr_{\gG^{(0)}}$ and $r\colon \gH\to \gH^{(0)}$. The anchor maps
are given by
\[
\rho\colon \Omega \to \gG^{(0)},\ (g,\eta)\mapsto g\LazyAnd \sigma\colon \Omega \to \gH^{(0)},\ (g,\eta)\mapsto
s(\eta).
\]
The action of $\gG$ on $\Omega$ is given by
\[
\gamma (g, \eta) := (r(\gamma),\ f(\gamma)\eta)
\]
for all $\gamma\in \gG$, $g\in \gG^{(0)}$ and $\eta\in \gH$ such that $s(\gamma) = g$ and $f(g) = r(\eta)$. The action of $\gH$ on $\Omega$ is given by multiplication from the right in the second component. The morphism $[\Graph(f)]$ given by $\Graph(f)$ is denoted by $\Morph(f)$.

The \demph{identity morphism of $\gG$} is defined as $\Morph(\id_{\gG})$, where $\id_{\gG}\colon \gG\to \gG$ is the strict identity morphism on $\gG$. It is the equivalence class of the graph $\gG$, where we consider $\gG$ to be a bimodule over itself, because $\gG^{(0)} \times_{\gG^{(0)}} \gG$ is equivalent to $\gG$. We will denote the morphism $\Morph(\id_{\gG})$ also by $\id_{\gG}$.

Let $\Omega$ be a graph from $\gG$ to $\gH$ and $\Omega'$ a graph from $\gH$ to $\gK$. We define on $\Omega'':=\Omega \times_{\gH} \Omega'$ the structure of a $\gG$-$\gK$-bimodule as at the end of Section~\ref{Subsection:Groupoids:BasicDefinitions}. Then this bimodule is a graph from $\gG$ to $\gK$, called the \demph{composition of $\Omega$ and $\Omega'$}.

The definition of the composition of graphs lifts to equivalence classes. Hence we have also defined the \demph{composition of morphisms}. The locally compact Hausdorff groupoids, together with their morphisms, form a category. The assignment $f\mapsto \Morph(f)$ is a functor from the category of locally compact Hausdorff groupoids with the strict morphisms as morphisms to the category of locally compact Hausdorff groupoids with all (generalised) morphisms.

\subsection{Equivalences}\label{Subsection:EquivalencesOfGroupoids}

\begin{definition}[$\gG$-$\gH$-equivalence]\label{DefinitionEquivalenceOfGroupoids} A $\gG$-$\gH$-bimodule $\Omega$ is called a $\gG$-$\gH$-equivalence bimodule if
\begin{enumerate}
    \item it is free and proper both as a $\gG$- and an $\gH$-space;
    \item the anchor map $\rho\colon \Omega\to \gG^{(0)}$ induces a homeomorphisms from $\Omega/\gH$ to $\gG^{(0)}$; and
    \item the anchor map $\sigma\colon \Omega\to \gH^{(0)}$ induces a homeomorphism from $\gG\backslash \Omega$ to $\gH^{(0)}$.
\end{enumerate}
We call $\gG$ and $\gH$ \emph{ equivalent} if such an equivalence exists.,
\end{definition}

\noindent Equivalence of groupoids is an equivalence relation. 
Actually, being equivalent is the same as being isomorphic in the generalised sense (see \cite{Paravicini:07}, Proposition~6.1.30).

If $\gG$ and $\gH$ are equivalent, then the categories of $\gG$-spaces and of $\gH$-spaces are the same: Let $\Omega$ be an equivalence between $\gG$ and $\gH$ and let $X$ be a (locally compact Hausdorff) $\gH$-space. Then $\Ind_{\gH}^{\gG} X:= \Omega^* X:= \Omega \times_{\gH} X$ is a (locally compact Hausdorff) $\gG$-space. Note that $\Ind_{\gG}^{\gH} \Ind_{\gH}^{\gG} X$ is isomorphic to $X$ as an $\gH$-space. If $X$ is proper, then also $\Ind_{\gH}^{\gG} X$ is proper. If $X$ is $\gH$-compact, then $\Ind_{\gH}^{\gG}X$ is $\gG$-compact. If $\uEgH$ is a universal proper locally compact $\gH$-space, then $\Ind_{\gH}^{\gG} \uEgH$ is a universal proper locally compact $\gG$-space.

\subsection{The linking groupoid}\label{SubsectionLinkingGroupoid}\label{SUBSECTIONLINKINGGROUPOID}

In this section we introduce some apparently new notions for locally compact groupoids which we borrow from the theory of Morita equivalences of C$^*$-algebras. It turns out that general equivalences of groupoids can be decomposed into an inclusion and the inverse of an inclusion into a larger locally compact Hausdorff groupoid called the linking groupoid of the given equivalence.

So let $\Omega$ be an equivalence between the locally compact groupoids $\gG$ and $\gH$.

\subsubsection{Definition}

\begin{definition}[The linking groupoid] Let $\gL$ be the locally compact Hausdorff space $\gL:=\gG\sqcup \Omega \sqcup \flipX{\Omega} \sqcup \gH$ and $\gL^{(0)}:= \gG^{(0)} \sqcup \gH^{(0)}$. Define the range and source maps of $\gL$ as
\[
r_{\gL}\colon {\gL} \to {\gL}, \ \left\{\begin{array}{ccrclcl} \gG&\ni&\gamma &\mapsto& r_{\gG}(\gamma) &\in& \gG^{(0)}\\
\Omega&\ni&\omega&\mapsto &\rho(\omega) &\in& \gG^{(0)}\\
\flipX{\Omega}&\ni&\flipx{\omega} &\mapsto& \rho(\flipx{\omega}) = \sigma(\omega) &\in & \gH^{(0)}\\ \gH&\ni& \eta &\mapsto
& r_{\gH}(\eta) & \in & \gH^{(0)} \end{array}\right\},
\]
and
\[
s_{\gL}\colon \gL \to \gL, \ \left\{\begin{array}{rcrclcl} \gG&\ni&\gamma &\mapsto& s_{\gG}(\gamma) &\in& \gG^{(0)}\\ \Omega&\ni&\omega&\mapsto &\sigma(\omega) &\in& \gH^{(0)}\\
\flipX{\Omega}&\ni&\flipx{\omega} &\mapsto& \sigma(\flipx{\omega}) =\rho(\omega) &\in & \gG^{(0)}\\ \gH&\ni& \eta
&\mapsto & s_{\gH}(\eta) & \in & \gH^{(0)} \end{array}\right\}.
\]
With these definitions,
\[
\gL*\gL = \gG\!*\!\gG\  \sqcup \ \gG\!*\! \Omega\  \sqcup\ \Omega\!*\!\flipX{\Omega}\ \sqcup\ \Omega\!*\!\gH\ \sqcup\
\flipX{\Omega} \!*\!\gG \ \sqcup\ \flipX{\Omega} \!*\! \Omega\ \sqcup\  \gH \!*\!\flipX{\Omega}\ \sqcup\ \gH\!*\!\gH.
\]
Define a composition map from $\gL*\gL$ to $\gL$ as the obvious map on the components $\gG\!*\!\gG$, $\gG\!*\!\Omega$, $\Omega\!*\!\gH$, $\flipX{\Omega} \!*\!\gG$, $\gH \!*\!\flipX{\Omega}$, and $\gH\!*\!\gH$; on $\flipX{\Omega} \!*\! \Omega$ and $\Omega\!*\!\flipX{\Omega}$ we take the factor map onto $\flipX{\Omega} \times_{\gG} \Omega$ and $\Omega\times _{\gH} \flipX{\Omega}$, which we identify with $\gH$ and $\gG$, respectively. In other words, a pair $(\omega^{-1}, \omega') \in \flipX{\Omega} \!*\! \Omega$ is mapped to its inner product $\langle \omega,\omega'\rangle_{\gH} \in \gH$, which is the unique element $\eta$ of $\gH$ such that $\omega' = \omega \eta$ (and similarly for $\Omega\!*\!\flipX{\Omega}$).
\end{definition}

\begin{proposition}
$\gL$ is a locally compact Hausdorff groupoid with open range and source maps. The inversion on $\gL$ is the map
\[
\gL \to \gL, \ \left\{\begin{array}{ccrclcl} \gG&\ni&\gamma &\mapsto& \gamma^{-1} &\in& \gG\\ \Omega&\ni&\omega&\mapsto &\flipx{\omega} &\in& \flipX{\Omega}\\
\flipX{\Omega}&\ni&\flipx{\omega} &\mapsto& \omega &\in & \Omega\\ \gH&\ni& \eta &\mapsto & \eta^{-1} & \in & \gH
\end{array}\right\}.
\]
\end{proposition}

\subsubsection{Full subsets}\label{Subsubsection:FullSubsets}

A subset $U\subseteq \gG^{(0)}$ of the unit space of $\gG$ is called \demph{full} if $\gG_U \circ \gG^U = \gG$, i.e., if every element $\gamma$ of $\gG$ can be written as a product $\gamma_1 \gamma_2$ with $\gamma_1$ starting in $U$ (and $\gamma_2$ ending in $U$).

\begin{proposition} Let $U \subseteq \gG^{(0)}$ be a full open subset. Then $\gG_U^U$ is a locally compact Hausdorff groupoid with open range and source maps and $\gG^U$ is a $\gG_U^U$-$\gG$-equivalence.
\end{proposition}
\begin{proof}
We have $\gG^U= r^{-1}(U)$ and $\gG^U_U=r^{-1}(U) \cap s^{-1}(U)$; hence $\gG^U$ and $\gG^{U}_U$ are open subsets of a locally compact Hausdorff space and hence themselves locally compact Hausdorff. The restriction of the range and source maps and of the multiplication on $\gG$ to $\gG^{U}_U$ turn $\gG^{U}_U$ into a topological groupoid; its range and source maps are open because they are restrictions of open maps to an open subset. Moreover, $\gG^U$ is a $\gG^U_U$-$\gG$-bimodule in a canonical way. In the same way as one proves that the action of $\gG$ on itself is proper and free one can show that the left and the right action on $\gG^U$ are proper and free. Because the range map of $\gG$ is open, the map from $\gG^{U} / \gG$ to $U$ induced from the range map is open. It is clearly surjective and continuous. Also injectivity follows by a standard argument: If $\gamma_1,\gamma_2\in \gG$ such that $r(\gamma_1)  =r(\gamma_2) \in U$, then $\gamma_2 = \gamma_1 (\gamma_1^{-1} \gamma_2)$, so $\gamma_1$ and $\gamma_2$ are in the same $\gG$-orbit in $\gG^U$. So far, we have not used that $U$ is a full subset. The map from $\gG^U_U \backslash \gG^U$ to $\gG^{(0)}$ induced by the source map $s$ is open, continuous and injective (by the same argument we have used as above). To see that it is surjective let $x\in \gG^{(0)}$. We identify $x$ with the corresponding unit element in $\gG$. By the fullness of $U$ we can find $\gamma_1,\gamma_2 \in \gG$ such that $x=\gamma_1\gamma_2$ and $s(\gamma_1) = r(\gamma_2)\in U$ (note that $\gamma_1 = \gamma_2^{-1}$). Hence $\gamma_2$ is an element of $\gG^U$ such $s(\gamma_2) =x$.
\end{proof}

\begin{corollary} Let $\Omega$ be a $\gG$-$\gH$-equivalence. Form the linking groupoid $\gL$ as above. Then $U:=\gG^{(0)}$ is a full open and closed subset of $\gL^{(0)}$ and $\gL_U^U$ can be identified with $\gG$. So $\gG$ is equivalent to $\gL$. In a similar fashion, $\gH$ is equivalent to $\gL$.
\end{corollary}

\subsubsection{Induction as restriction}

The categories of $\gG$-spaces, $\gH$-spaces and $\gL$-spaces are mutually equivalent. Given an $\gH$-space $X$, the corresponding $\gG$-space $\Ind_{\gH}^{\gG}X$ is defined as $\Omega \times_{\gH} X$. The $\gL$-space that corresponds to $X$ can be constructed in two alternative ways. On the one hand, one can use the equivalence $\Omega \sqcup \gH$ between $\gL$ and $\gH$ to form $(\Omega \sqcup \gH) \times_{\gH} X$.  On the other hand, $(\Ind_{\gH}^{\gG} X) \sqcup X$ can be equipped with a canonical action of $\gL$; e.g. if $\omega\in \Omega\subseteq \gL$ and $x\in X$, then $\omega \cdot x$ is defined if $\sigma_{\Omega}(\omega) = \rho_X(x)$, in case of which it is defined as the equivalence class of $(\omega,x)$ in $\Ind_{\gH}^{\gG}X$. Note that both ways to define $\Ind_{\gH}^{\gL} X$ yield the same result.

Conversely, if $X$ is an $\gL$-space, then $\Ind_{\gL}^{\gH} X$ can not only be realised as $(\Omega \sqcup \gH)^{-1} \times_{\gL} X$, but also as $X\restr_{\gH^{(0)}}$, i.e., as the subset of $X$ of those elements which are mapped to elements of $\gH^{(0)}$ under the anchor map. This description of induction as restriction is going to be useful in Section~\ref{Subsection:InductionAndAssembly}.

\subsection{The pullback of groupoids}\label{Subsection:Groupoids:Pullback}

Let $X:=\gG^{(0)}$ and $Y$ be locally compact Hausdorff spaces and let $p\colon Y\to X$ be a continuous map. Then we define the \demph{pullback groupoid}\footnote{What we call $p^*(\gG)$ is called $\gG_Y$ in \cite{LeGall:99} and $\gG[Y]$ in \cite{Tu:04}.} $p^*(\gG)$ to be the fibre product of $Y\times Y$ and $\gG$ over $X\times X= \gG^{(0)}\times \gG^{(0)}$, i.e., $p^*(\gG)$ is defined as the pullback in the following diagram:
\[
\xymatrix{p^*(\gG) \ar[d] \ar[r] &  Y\times Y \ar[d]^{p\times p} \\
\gG \ar[r]^-{(r,s)} & X \times X }
\]
It can be realised as follows:
\[
p^*(\gG)\cong  \left\{(z,\gamma,y)\in Y\times \gG\times Y:\ s(\gamma) = p(y),\ r(\gamma) = p(z)\right\}
\]
and the unit space of $p^*(\gG)$ can be identified with $Y$ via $y\mapsto (y,\epsilon(p(y)),y)$. The source and range functions are given by $r_{p^*(\gG)} (z,\gamma,y)=z$ and $s_{p^*(\gG)} (z,\gamma,y) = y$; moreover, the composition is given by
\[
(z,\gamma,y) \circ (z',\gamma',y') = (z, \gamma\circ \gamma', y')
\]
and is defined if and only if $y=z'$. The inverse is given by $(z,\gamma,y)^{-1} = (y,\gamma^{-1},z)$. Note that $p^*(\gG)$ is a locally compact Hausdorff groupoid.

There is a canonical strict morphism from $p^*(\gG)$ to $\gG$, appearing in the above diagram, which we call $p$ (abusing the notation); it is given explicitly by $(z,\gamma,y) \mapsto \gamma$.

\noindent Under certain conditions, the graph of $p\colon p^*(\gG) \to \gG$ is an equivalence:

\begin{proposition}\label{Proposition:GraphPEquivalence} Let $\gG$ have open range and source maps. The strict morphism $p\colon p^*(\gG) \to \gG$ has graph
\[
\Graph(p)=p^*(\gG)^{(0)} \times_{\gG^{(0)}} \gG = Y \times_{\gG^{(0)}} \gG = \{(y,\gamma) \in Y\times \gG:\ p(y) =
r(\gamma)\}.
\]
If $p\colon Y\to X$ is open and surjective, then $\Graph(p)$ is an equivalence.
\end{proposition}

\begin{defprop}[The strict morphism $f_{\Omega}$]\label{DefPropfOmega} Let $\Omega$ be a graph from $\gG$ to $\gH$. Write $\langle \cdot,\cdot \rangle_{\gH}$ for the $\gH$-valued inner product from $\Omega \times_{\rho} \Omega$ to $\gH$, i.e., $\langle \omega, \omega'\rangle_{\gH}$ is defined to be the unique element $\eta$ of $\gH$ such that $\omega \eta = \omega'$.

Define $f_{\Omega}(\omega',\gamma,\omega) := \langle\omega',\gamma \omega\rangle_{\gH}$ for all $(\omega',\gamma,\omega)\in \rho^*(\gG)$. Then $f_{\Omega}\colon \rho^*(\gG) \to \gH$ is a strict morphism extending $\sigma\colon \Omega=\rho^*(\gG)^{(0)} \to \gH^{(0)}$.
\end{defprop}

\noindent Note that every generalised morphism can be written as the composition of an equivalence and a strict morphism:

\begin{proposition}\label{Proposition:ZerlegeMorphismen} Let $\Omega$ be a graph from $\gG$ to $\gH$. Then $\Morph(f_{\Omega})$ makes the following diagram commutative
\begin{equation}\label{EquationGroupoidMorphTriangle}
\xymatrix{ \rho^*(\gG) \ar[dd]_{\Morph(\rho)}^{\cong} \ar[ddrr]^{\Morph\left(f_{\Omega}\right)}& &\\ &&\\ \gG\ar[rr]^{[\Omega]}
&& \gH
}
\end{equation}
\end{proposition}

\subsection{Locally compact groupoids with Haar systems}\label{SubsectionHaarSystems}

We will now analyse how Haar systems behave under the constructions we have introduced above: Can one lift Haar systems to equivalent groupoids, to the pullback of a groupoid or to linking groupoids? To be able to discuss these questions systematically, we introduce Haar systems not only on groupoids but also on spaces on which groupoids act. A basic concept which makes this discussion much simpler is the notion of a continuous field of measures:

Let $X$ and $Y$ be locally compact Hausdorff spaces and let $p\colon Y\to X$ be a continuous map. A \demph{continuous field of measures} on $Y$ over $X$ (with coefficient map $p$) is a family $(\mu_x)_{x\in X}$ of measures on $Y$ such that $\supp \mu_x \subseteq Y_x:=p^{-1}(\{x\})$, for all $x\in X$, and such that, for all $\varphi \in \Cont_c(Y)$,
\begin{equation}\label{Equation:ConditionContinuousFieldOfMeasures}
\mu(\varphi)\colon X\to \C,\ x\mapsto \int_{y\in Y_x} \varphi(y)\rmd \mu_x(y),
\end{equation}
is an element of $\Cont_c(X)$. It is called \demph{faithful} if $\supp \mu_x = Y_x$ for all $x\in X$. Basic facts concerning continuous fields of measures can be found in Appendix B of \cite{Paravicini:07}; note, for example, that $p$ is automatically open if there exists a faithful continuous field of measures with coefficient map $p$.

\begin{definition}[Haar system]A left Haar system on a left $\gG$-space $\Omega$ with (open and) surjective anchor map $\rho$ is a faithful continuous field $(\lambda_{\Omega}^g)_{g\in \gG^{(0)}}$ of measures on $\Omega$ over $\gG^{(0)}$ with coefficient map $\rho$ having the following invariance property
\begin{equation}\label{Equation:HaarSystemInvariance}
\forall \gamma \in \gG\ \forall \varphi \in \Cont_c(\Omega):\  \int_{\omega\in \Omega} \varphi(\omega) \rmd \lambda_{\Omega}^{r(\gamma)}(\omega) = \int_{\omega\in \Omega} \varphi(\gamma \omega) \rmd \lambda_{\Omega}^{s(\gamma)}(\omega).
\end{equation}
\noindent Similarly, right Haar systems are defined, and using that $\gG$ acts on itself on the left, we define a left Haar system on the groupoid $\gG$ to be a left Haar system for this action. Note that such a Haar system need not exist; the range and source maps of a locally compact Hausdorff groupoid admitting a Haar system are automatically open.
\end{definition}

\noindent In \cite{Paravicini:07}, the following constructions are described in detail:

\begin{enumerate}
\item Let $X$ and $Y$ be locally compact Hausdorff spaces and let $p\colon Y\to X$ be an open and surjective continuous map. Then Haar systems on the groupoid $Y \times_X Y$ correspond exactly to the faithful continuous fields of measures on $Y$ over $X$ with coefficient map $p$.\label{Page:HaarsystemsOnYXY}

\item If $\Omega$ is a graph from $\gG$ to $\gH$ (with anchor maps $\rho$ and $\sigma$), then one can construct from a left Haar system $(\lambda_{\gH}^h)_{h\in \gH^{(0)}}$ on $\gH$ a left Haar system $(\lambda_{\Omega}^g)_{g\in \gG^{(0)}}$ on $\Omega$ for the action of $\gG$; if $\varphi \in \Cont_c(\Omega)$ and $g\in \gG^{(0)}$, then
\[
\lambda_{\Omega}^g:= \int_{\eta \in \gH^{\sigma(\omega)}} \varphi(\omega\eta) \rmd \lambda_{\gH}^{\sigma(\omega)}(\eta),
\]
where $\omega$ is some arbitrary element of $\Omega$ such that $\rho(\omega)=g$.

\item If $\Omega$ is an equivalence between $\gG$ and $\gH$, and $\gG$ and $\gH$ carry left Haar systems, then the induced left Haar systems on $\Omega$ (for the $\gG$-action) and $\Omega^{-1}$ (for the $\gH$-action) combine to a left Haar system on the linking groupoid.

\item If $\Omega$ is a left Haar $\gG$-space with anchor map $\rho$ and Haar system $\lambda_{\Omega}$ and also $\gG$ itself carries a left Haar system $\lambda$, then the pullback groupoid $\rho^*(\gG)$ carries a left Haar system as well: If $\omega\in \Omega$ and $\varphi\in \Cont_c(\rho^*(\Omega))$, then define
\[
\lambda^{\omega}_{\rho^*(\Omega)}(\varphi):= \int_{\omega'\in \rho^{-1}(\rho(\omega))} \int_{\gamma\in \gG^{\rho(\omega)}} \varphi(\omega',\gamma,\gamma^{-1}\omega) \rmd \lambda_{\gG}^{\rho(\omega)} (\gamma) \rmd \lambda^{\rho(\omega)}_{\Omega}(\omega').
\]

 This applies in particular to the situation that $\Omega$ is a graph between groupoids with Haar systems because then also $\Omega$ carries a left Haar system (as we have said above in 2.).

\item If $\Omega$ is a left $\gG$-space and $\gG$ carries a left Haar system, then also the groupoid $\gG \ltimes \Omega$ carries a left Haar system.
\end{enumerate}

\section{Fields of Banach spaces, Banach algebras and Banach pairs}\label{Section:Fields}

In preparation of the next section, we introduce upper semi-continuous fields of Banach spaces etc. These concepts are crucial for the definition of $\gG$-Banach algebras where $\gG$ is a locally compact groupoid. Note that the notion of a $\gG$-C$^*$-algebra as defined in \cite{LeGall:94} uses the concept of a $\Cont_0(X)$-C$^*$-algebra introduced in \cite{Kasparov:88} instead, where $X$ denotes the unit space of $\gG$. Every $\Cont_0(X)$-C$^*$-algebra is however an upper semi-continuous field of Banach algebras over $X$ in a canonical way, but for general $\gG$-Banach algebras the viewpoint of fields seems more adequate.

A general reference for most of this section is \cite{Lafforgue:06}; apart from some technical refinements there is also a small but notable novelty: the notion of a compact operator in the context of fields of Banach pairs which we compare to the notion of a ``locally compact operator''.

Throughout Section~\ref{Section:Fields}, let $X$ be a locally compact Hausdorff space.

\subsection{Fields of Banach spaces and linear maps}

\begin{definition}[Upper semi-continuous field of Banach spaces] An \demph{upper semi-continuous field of Banach spaces} over $X$ is a pair $E=\left((E_x)_{x\in X},\ \ContSect\right)$, where $(E_x)_{x\in X}$ is a family of Banach spaces and $\ContSect\subseteq \prod_{x\in X} E_x$ is a set, which has the following properties:
\begin{enumerate}
\item [(C1)] $\ContSect$ is a complex linear subspace of $\prod_{x\in X}E_x$;
\item [(C2)] for all $x\in X$, the evaluation map $\ev_x\colon \ContSect \to E_x,\ \xi\mapsto \xi(x)$, has dense image;
\item [(C3)] for all $\xi\in \ContSect$, the map $\abs{\xi}\colon X\to \R_{\geq 0},\ x\mapsto \norm{\xi(x)}_{E_x}$, is upper semi-continuous;
\item [(C4)]\label{ConditionContFieldC4} if $\xi \in \prod_{x\in X}E_x$ and if, for all $x_0\in X$ and all $\varepsilon>0$, there is an element $\gamma\in \ContSect$ and a neighbourhood $U$ of $x_0$ in $X$ such that for all $x\in U$ we have $\norm{\xi(x)-\gamma(x)}_{E_x} \leq \varepsilon$, then also $\xi$ belongs to $\ContSect$.
\end{enumerate}
\end{definition}

\noindent Condition (C4) just says that an element of $\prod_{x\in X} E_x$ which can be approximated locally by elements of $\ContSect$ is itself in $\ContSect$. Note that all elements of $\ContSect$ are locally bounded by (C3). Instead of ``upper semi-continuous field'' we will usually say ``u.s.c.~field'' of Banach spaces. A \demph{continuous} field of Banach spaces satisfies the extra condition that $\abs{\xi}$ is not only upper semi-continuous but  continuous for all $\xi \in \ContSect$.

There are trivial examples of continuous fields of Banach spaces, namely the constant fields: If $E_0$ is a Banach space, then define $E_x:=E_0$ for every $x\in X$ and let $\ContSect$ be the space $\Cont(X,E_0)$ of all continuous maps from $X$ to $E_0$. Then this gives a continuous field $(E_0)_X$ of Banach spaces, called the \demph{constant} field over $X$ with fibre $E_0$.

For a u.s.c~field $E=\left((E_x)_{x\in X}, \ \ContSect\right)$ of Banach spaces we make the following definitions: The elements of $\ContSect$ are called the \demph{sections of $E$}. We will also write $\ContSect(X,E)$ for $\ContSect$. The Banach space of bounded sections is denoted by $\ContSect_b(X,E)$. The Banach space of all sections of $E$ vanishing at infinity is denoted by $\ContSect_0(X,E)$. The linear space of all sections of $E$ with compact support is denoted by $\ContSect_c(X,E)$. Note that
\[
\ContSect_c(X,E) \subseteq  \ContSect_0(X,E) \subseteq \ContSect_b(X,E) \subseteq \ContSect(X,E).
\]

\noindent If $E$ is a u.s.c.~field of Banach spaces and if $\xi\in \ContSect(X,E)$ is a section of $E$ and $\chi \in \Cont(X)$, then $\chi \xi \in \ContSect(X,E)$. Moreover, $\ContSect_0(X,E)$ is a non-degenerate $\Cont_0(X)$-module in this way. For every $x\in X$, the evaluation map $\ContSect_0(X,E) \to E_x$ has not only dense image but is a metric surjection.

\begin{definition}[Continuous field of linear maps] Let $E$ and $F$ be u.s.c.~fields of Banach spaces. Then a \demph{continuous field of linear maps} from $E$ to $F$ is a family $(T_x)_{x\in X}$ such that
\begin{enumerate}
    \item $T_x \in \Lin\left(E_x,F_x\right)$ for all $x\in X$;
    \item $\forall \xi \in \ContSect(X,E):\quad  T\circ\xi\colon x\mapsto T_x(\xi(x)) \ \in \ \ContSect(X,F)$;
    \item the function $x\mapsto \norm{T_x}$ is locally bounded\footnote{In \cite{Lafforgue:06}, continuous fields of linear maps are defined leaving out our third condition (D\'{e}finition 1.1.7); however, Proposition 1.1.9 of the same article states that Condition 3.~is automatic if $X$ is metrisable. A more general result along these lines is proved in Appendix E.2 of \cite{Paravicini:07}.} on $X$, i.e., for all $x\in X$ there is a neighbourhood $U$ of $x$ in $X$ such that $ \sup_{u\in U} \norm{T_u} < \infty$.
\end{enumerate}
The set of all continuous fields of linear maps from $E$ to $F$ will be denoted by $\Linloc(E,F)$. The subset of (globally) bounded continuous fields of linear maps from $E$ to $F$ is denoted by $\Lin(E,F)$.
\end{definition}

\subsection{Fields of bilinear maps and the tensor product}

In this subsection, let $E$, $F$, $G$ be u.s.c.~fields of Banach spaces over $X$.

The \demph{internal product} $E\times_{X} F$ of $E$ and $F$ is the u.s.c.\ of Banach spaces over $X$ given by the following data: The underlying family of Banach spaces is the family $E\times_X F = (E_x\times F_x)_{x\in X}$,  and the space of sections is
\[
\ContSect:= \left\{x\mapsto (\xi(x),\eta(x)):\ \xi\in \ContSect(X,E), \eta\in \ContSect(X,F)\right\}.
\]
The set $\ContSect$ satisfies condition (C1) - (C4), thus it defines the structure of a u.s.c.~field of Banach
spaces on $E\times_X F$. Similarly, we define the \demph{internal sum} $E\oplus_X F$ of $E$ and $F$ over $X$, the difference being that we take the sum-norm instead of the sup-norm on the fibres.

A \demph{continuous field of bilinear maps from $E\times_X F$ to $G$} is a family $(\mu_x)_{x\in X}$ of continuous bilinear maps $\mu_{x}\in \Mult(E_{x}, F_{x}; G_x)$ for all $x\in X$ such that
\begin{enumerate}
\item $\forall \xi\in \ContSect(X,E)\ \forall \eta\in \ContSect(X,F):\quad x\mapsto \mu_x \left(\xi(x),\ \eta(x)\right) \ \in \ \ContSect(X,G).$

\item $\mu$ is locally bounded.
\end{enumerate}

\noindent We define the (internal) \demph{tensor product $E\otimes_X F$ of $E$ and $F$} to be the following u.s.c.~field of Banach spaces over $X$: for all $x\in X$, the fibre of $E\otimes_X F$ over $x$ is $E_x\otimes^\pi F_x$; to define the sections of $E\otimes_X F$, let $\Lambda$ be the $\C$-linear span of all selections of the family $E\otimes_X F$ given by $x\mapsto \xi(x) \otimes \eta(x)$, where $\xi$ runs through $\ContSect(X,E)$ and $\eta$ runs through $\ContSect(X,F)$. Then $\Lambda$ satisfies conditions (C1), (C2) and (C3), so there is a smallest subset $\Gamma=:\Gamma(X, E\otimes_X F)$ of $\prod_{x\in X} E_x \otimes F_x$ which contains $\Lambda$ and satisfies (C1)-(C4) (use Proposition 1.1.4 of \cite{Lafforgue:06}). There is a canonical contractive continuous field of bilinear maps $\pi=(\pi_x)_{x\in X}$ from $E\times_X F$ to $E\otimes_X F$. Note that there property (C3) is not completely trivial to show, compare page 6 of \cite{Lafforgue:06}.

\subsection{Fields of Banach algebras}

A \demph{u.s.c.~field of Banach algebras over $X$} is an upper semi-continuous field $A$ of Banach spaces over $X$ together with a continuous field of bilinear maps $\mu\colon A\times_X A \to A$ such that $(A_x, \mu_x)$ is a Banach algebra for all $x\in X$. In particular, this means that $\mu$ is bounded by one. A field of Banach algebras $A$ over $X$ (with multiplication $\mu$) is called \demph{non-degenerate} if $\mu_x$ is non-degenerate for all $x\in X$, i.e., the span of $A_xA_x$ is dense in $A_x$.

Let $A$ and $B$ be u.s.c.~fields of Banach algebras over $X$. Then a \demph{homomorphism (of fields of Banach algebras) from $A$ to $B$} is a continuous field of homomorphisms of Banach algebras from $A$ to $B$, i.e., a continuous field $(\varphi_x) _{x\in X}$ of linear maps from $A$ to $B$ such that $\varphi_x$ is a (contractive) homomorphism of Banach algebras from $A_{x}$ to $B_x$. In particular, such a $\varphi$ is bounded by one.

Let $B$ be a u.s.c.~field of Banach algebras over $X$. Then we define the \demph{fibrewise unitalisation}
\[
\unital{B} = B\oplus_X \C_X= \left(\unital{B_x}\right)_{x\in X}
\]
to be the following u.s.c.~field of Banach algebras: For all $x\in X$, the fibre of $\unital{B}$ is the unitalisation $\unital{B_x}$ of the fibre $B_x$ of $B$. The sections of $\unital{B}$ are $\ContSect(X,B) \oplus \Cont(X)$.

\subsection{Fields of Banach modules}\label{Subsection:FieldsOfBanachModules}

Let $A$, $B$ and $C$ be u.s.c.~fields of Banach algebras over $X$. 

A \demph{right Banach $B$-module} is an upper semi-continuous field $E$ of Banach spaces over $X$ together with a continuous field of bilinear maps $\mu^E\colon E\times_X B \to E$ such that, for all $x\in X$, $E_x$ is a right Banach $B_x$-module with the $B_x$-action $\mu^E_x$. In particular, this means that $\mu^E$ is bounded by one. The module $E$ is called \demph{non-degenerate} if $\mu^E_x$ is non-degenerate for all $x\in X$, i.e., the span of $E_xB_x$ is dense in $E_x$. Left Banach $A$-modules and Banach $A$-$B$-bimodules are defined similarly.

Let $E$ and $F$ be right Banach $B$-modules. Then a \demph{$B$-linear field of operators from $E$ to $F$} (or just a \demph{$B$-linear operator from $E$ to $F$}) is a continuous field $T$ of linear maps from $E$ to $F$ such that $T_x$ is $B_x$-linear (on the right) for all $x\in X$. We denote the space of all such $T$ by $\Linloc_B(E,F)$.

As usual, the field $T$ is called bounded if $\norm{T}:= \sup_{x\in X}\norm{T_x} <\infty$. We denote the bounded $B$-linear operators from $E$ to $F$ by $\Lin_B(E,F)$.

Let $B'$ be another field of Banach algebras over $X$ and let $\psi\colon B\to B'$ be a continuous field of homomorphisms. Let $E$ be a right Banach $B$-module and $E'$ be a right Banach $B'$-module. Then a \demph{homomorphism $\Phi_\psi$ (of u.s.c.~fields of Banach modules) from $E_B$ to $E'_{B'}$ with coefficient map $\psi$} is a \emph{contractive} continuous field $\Phi$ of linear maps from $E$ to $E'$ such that $\Phi_x$ is a homomorphism with coefficient map $\psi_x$ from $(E_x)_{B_x}$ to $(E'_x)_{B'_x}$ for all $x\in X$, compare \cite{Paravicini:07:Morita:erschienen}, Definition~1.1. An analogous definition can be formulated for left Banach modules and Banach bimodules.

Let $E_1$ be a right Banach $B$-module and $E_2$ a left Banach $B$-module. Let $F$ be a u.s.c.~field of Banach spaces over $X$. A continuous field $\mu$ of bilinear maps from $E_1 \times_X E_2$ to $F$ is called \demph{$B$-balanced} if $\mu_x\colon (E_1)_x \times (E_2)_x \to F_x$ is $B_x$-balanced for all $x\in X$. The \demph{balanced tensor product $E_1 \otimes_B E_2$ of $E_1$ and $E_2$} is defined in complete analogy to the internal tensor product of fields of Banach spaces introduced above; its fibre over $x\in X$ is $(E_1)_x \otimes_{B_x} (E_2)_x$. Compare Section~4 of \cite{Paravicini:07:Morita:erschienen}.

Let $E$ and $E'$ be right Banach $B$-modules and $F$ a Banach $B$-$C$-bimodule. For all $T\in \Linloc_B(E,E')$ define $T\otimes 1\in \Linloc_C\left(E\otimes_B F,\ E'\otimes_B F\right)$ as the family $\left(T_x \otimes_{B_x} \id_{F_x}\right)_{x\in X}$. Note that the assignment $T\mapsto T\otimes 1$ is linear and functorial. If $T$ is bounded, then $\norm{T\otimes 1} \leq \norm{T}$.

Let $B'$ be a u.s.c.~field of Banach algebras and $\psi\colon B\to B'$ a continuous field of homomorphisms. Let $E$ be a right Banach $B$-module. Then $\psi_*(E):=E\otimes_{\unital{B}} \unital{B'}$ is a right Banach $B'$-module, called the \demph{pushout of $E$ along $\psi$}. The fibre of $\psi_*(E)$ at $x$ is $(\psi_x)_*(E_x)$.

\subsection{Fields of Banach pairs}

Let $A$ and $B$ be u.s.c.~fields of Banach algebras over $X$. 

\begin{definition}[Field of Banach pairs] A \demph{Banach $B$-pair} is a pair $E=(E^<,E^>)$ such that $E^<$ is a left Banach $B$-module and $E^>$ is a right Banach $B$-module, together with a contractive continuous field of bilinear maps $\langle, \rangle \colon E^< \times_X E^> \to B$, $B$-linear on the left and on the right. $E$ is called non-degenerate if $E^<$ and $E^>$ are non-degenerate Banach $B$-modules.

Define $E_x:=(E_x^<,E_x^>)$, which is a $B_x$-pair when equipped with the bracket $\langle,\rangle_x$.
\end{definition}

\begin{definition}[Linear operator between fields of Banach pairs] Let $E$ and $F$ be Banach $B$-pairs. Then a \demph{continuous field of $B$-linear operators from $E$ to $F$} (or just a \demph{$B$-linear operator from $E$ to $F$}) is a pair $(T^<,T^>)$ where $T^>$ is a locally bounded continuous field of $B$-linear operators from $E^>$ to $F^>$ and $T^<$ is a locally bounded continuous field of $B$-linear operators from $F^<$ to $E^<$ such that $T_x:=(T_x^<,T_x^>)$ is in $\Lin_{B_x}(E_x,F_x)$ for all $x\in X$. We denote the linear space of all such $T$ by $\Linloc_B(E,F)$.

A $B$-linear operator from $E$ to $F$ is called \demph{bounded} if $T^<$ and $T^>$ are bounded. The space of all bounded $B$-linear operators from $E$ to $F$ will be denoted by $\Lin_B(E,F)$. It is a Banach space when equipped with the obvious operations and the norm $\norm{T}:= \max\{\norm{T^<},\norm{T^>}\} = \sup_{x\in X}\norm{T_x}$.
\end{definition}

\noindent The notion of a homomorphism of Banach pairs was introduced in \cite{Paravicini:07}, see also \cite{Paravicini:07:Morita:erschienen}.

\begin{definition}[Homomorphism between fields of Banach pairs] Let $B$ and $B'$ be u.s.c.~fields of Banach algebras over $X$ and let $\psi\colon B\to B'$ be a continuous field of homomorphisms of Banach algebras. Let $E_B$ and $E'_{B'}$ be Banach pairs. Then a \demph{continuous field $\Phi$ of homomorphisms from $E$ to $E'$ with coefficient map $\psi$} is a pair $(\Phi^<,\Phi^>)$ where $\Phi^>$ is a continuous field of homomorphisms from $E^>$ to $E'^>$ and $\Phi^<$ is a continuous field of homomorphisms from $E^<$ to $E'^<$, both with coefficient map $\psi$, such that $\Phi_x:=(\Phi_x^<,\Phi_x^>)$ is a homomorphism with coefficient map $\psi_x$ from the pair $(E_x)_{B_x}$ to the pair $(E'_x)_{B'_x}$.
\end{definition}

\noindent Note that the composition of linear operators is again a linear operator and the composition of homomorphisms is again a homomorphism.

\begin{definition}[Banach $A$-$B$-pair] A Banach $A$-$B$-pair $E=(E^<,E^>)$ is a Banach $B$-pair $E$ such that $E^<$ is a Banach $B$-$A$-bimodule and $E^>$ is a Banach $A$-$B$-bimodule and the bracket $\langle,\rangle \colon E^< \times_X E^> \to B$ is $A$-balanced (which means that for all $x\in X$ the map $\langle,\rangle_x\colon E^<_x\times E^>_x \to B_x$ is $A_x$-balanced).
\end{definition}

\noindent There is an obvious notion of a homomorphism with coefficient maps between Banach $A$-$B$-pairs.
\smallskip

Using the definition of the balanced tensor product of fields of Banach modules we can define the balanced tensor product of fields of Banach pairs. Similarly, we can define the pushout of fields of Banach pairs along continuous fields of homomorphisms between u.s.c.~fields of Banach algebras. It has the usual functorial properties.

\subsection{Locally compact operators}

Let $B$ be a u.s.c.~field of Banach algebras over the locally compact Hausdorff space $X$. We now recall a concept from \cite{Lafforgue:06} which we call ``locally compact operator'' between $B$-pairs. It is complemented by the notion of a ``compact operator'', the difference being that compact operators are fields that vanish at infinity and can hence be approximated globally by finite rank operators, whereas locally compact operators are only locally bounded fields of operators which can nevertheless be approximated locally by finite rank operators (or by compact operators). Both notions can (and will) be used in the definition of $\KKban$ for groupoids.

\begin{definition}[Rank one operator] Let $E$ and $F$ be Banach $B$-pairs. Then we define for all $\xi^<\in \ContSect(X,E^<)$ and all $\eta^>\in \ContSect(X,F^>)$ the continuous field of operators $\ketbra{\eta^>}{\xi^<} := \left(\ketbra{\eta^>}{\xi^<}_x\right)_{x\in X}\in \Linloc_B(E,F)$ by
\[
\ketbra{\eta^>}{\xi^<}_x := \ketbra{\eta^>(x)}{\xi^<(x)}\in \Komp_{B_x}\left(E_x,F_x\right)
\]
for all $x\in X$.
\end{definition}

\noindent If $\xi^<$ and $\eta^>$ are bounded, then $\ketbra{\eta^>}{\xi^<}$ is bounded by $\norm{\xi^<} \norm{\eta^>}$. If $\xi^<$ and $\eta^>$ vanish at infinity, then so does $\ketbra{\eta^>}{\xi^<}$.

\begin{definition}[Locally compact operator]\label{DefProp:LocallyCompactFieldsOfOperators} \footnote{V.~Lafforgue calls such operators ``compact partout'' in \cite{Lafforgue:06}.} Let $E$ and $F$ be Banach $B$-pairs. A continuous field $T$ of $B$-linear operators is called \demph{locally compact} if, for all $x\in X$ and all $\varepsilon>0$, there is an open neighbourhood $U$ of $x$, an $n\in \N$ and $\xi_1^<,\ldots,\xi_n^< \in \ContSect(X,E^<)$ and $\eta_1^>,\ldots,\eta_n^>\in \ContSect(X,F^>)$ such that $\norm{T_u - \sum_{i=1}^n \ketbra{\eta_i^>(u)}{\xi_i^<(u)}} \leq  \varepsilon$ for all $u\in U$. The space of all locally compact operators from $E$ to $F$ is denoted by $\Komploc_B(E,F)$.
\end{definition}

\noindent In other words: If $\mF$ denotes the linear span of the operators of the form $\ketbra{\eta^>}{\xi^<}$, with $\xi^<\in \ContSect(X,E^<)$ and $\eta^>\in \ContSect(X,F^>)$, in the space $\Linloc_B(E,F)$, then $\Komploc_B(E,F)$ is the space of all operators that are locally approximable by elements of $\mF$.

\begin{lemma}\label{Lemma:FieldsCompositionOfLocallyCompactOperators} Let $E$, $F$ and $G$ be Banach $B$-pairs. Then $\Linloc_B(F,G) \circ \Komploc_B(E,F) \subseteq \Komploc_B(E,G)$ and $\Komploc_B(F,G) \circ \Linloc_B(E,F) \subseteq \Komploc_B(E,G)$.
\end{lemma}

\begin{Xample}\label{XampleStandardPairLocallyCompactAction} Let $B$ be a non-degenerate u.s.c.~field of Banach algebras over $X$. Then $\ContSect(X,B)$ acts by locally compact operators on the Banach $B$-pair $(B,B)$.
\end{Xample}

\subsection{Compact operators}\label{Subsection:CompactOperatorsOnFieldsOfBanachPairs}

Let $B$ be a u.s.c.~field of Banach algebras over the locally compact Hausdorff space $X$ and let $E$ and $F$ be Banach $B$-pairs.

\begin{definition}[Compact operators] A continuous field $T$ of $B$-linear operators from $E$ to $F$ is called \demph{compact} if, for all $\varepsilon>0$, there is an $n\in \N$ and $\xi_1^<,\ldots,\xi_n^< \in \ContSect_0(X,E^<)$ and $\eta_1^>,\ldots,\eta_n^>\in \ContSect_0(X,F^>)$ such that
\[
\norm{T-\sum_{i=1}^n \ketbra{\eta_i^>}{\xi_i^<}} = \sup_{x\in X} \norm{T_x-\sum_{i=1}^n
\ketbra{\eta_i^>(x)}{\xi_i^<(x)}} \leq \varepsilon.
\]
The compact operators from $E$ to $F$ are denoted by $\Komp_B(E,F)$.
\end{definition}

\noindent Note that the sections are taken to be vanishing at infinity. This means that, if $T$ is compact, then $(\norm{T_x})_{x\in X}$ is also vanishing at infinity. It follows that $\Komp_B(E,F) \subseteq \Lin_B(E,F)$ and $\Komp_B(E,F)$ is the closed linear span in $\Lin_B(E,F)$ of all operators of the form $\ketbra{\eta^>}{\xi^<}$. In particular, $\Komp_B(E,F)$ is a Banach space.

We will now justify the choice of the name ``locally compact operator'':

\begin{proposition}[Characterisation of locally compact operators]\label{Proposition:Characterisation:LocallyCompactOperator} Let $T$ be a continuous field of $B$-linear operators from $E$ to $F$. Then the following are equivalent:
\begin{enumerate}
\item $T$ is locally compact.

\item For all compact subsets $K$ of $X$ and all $\varepsilon>0$, there is an $n\in \N$ and $\xi_1^<, \ldots, \xi_n^< \in \ContSect(X,E^<)$ and $\eta_1^>,\ldots,\eta_n^>\in \ContSect(X,F^>)$ such that $\norm{T_k - \sum_{i=1}^n \ketbra{\eta_i^>(k)} {\xi_i^<(k)}} \leq  \varepsilon$ for all $k\in K$.

\item For all $x\in X$ and all $\varepsilon>0$, there is an open neighbourhood $U$ of $x$ and a compact operator $S\in  \Komp_B(E,F)$ such that $\norm{T_u -S_u} \leq \varepsilon$ for all $u\in U$.

\item For all compact subsets $K\subseteq X$ and all $\varepsilon>0$, there is an operator $S\in  \Komp_B(E,F)$ such that $\norm{T_k -S_k} \leq \varepsilon$ for all $k\in K$.

\item For all $\varphi \in \Cont_c(X)$, the field $\varphi T$ is compact.

\end{enumerate}
\end{proposition}
\begin{proof}
Assume that 1.~holds. Let $K\subseteq X$ be a compact subset. Let $\varepsilon>0$. For all $x\in X$, find $U_x$,
$n_x$, $\xi_{x,1}^<,\ldots,\xi_{x,n_x}^< \in \ContSect(X,E^<)$ and $\eta_{x,1}^>,\ldots,\eta_{x,n_x}^>\in \ContSect(X,F^>)$ as in the definition of local compactness for $T$. Then $\{U_x:\ x\in K\}$ is an open cover of $K$ so we can find a finite subset $A$ of $K$ such that $K\subseteq \bigcup_{a\in A} U_a$. Find a partition of unity $(\chi_a)_{a\in A}$ on $K$ subordinate to the cover $(U_a)_{a\in A}$. Then for all $k\in K$:
\[
\norm{T_k- \sum_{a\in A} \chi_a(k)\sum_{i=1}^{n_a}\left|\eta_{a,i}^>(k)\right\rangle\left\langle \xi_{a,i}^<(k)\right|} \leq
\varepsilon.
\]
This shows 1.~$\Rightarrow$ 2..

The same argument shows 3.~$\Rightarrow$ 4.. Since $X$ is locally compact, the implications 2.~$\Rightarrow$ 1.~and 4.~$\Rightarrow$ 3.~are trivial. Moreover, it is clear that 4.~implies 2.~ and 3.~implies 1.. Cutting down the sections used in the approximation in 2.~easily shows 2.~$\Rightarrow$ 4.. So the first four conditions are mutually equivalent. It is straightforward to show 4.~$\Leftrightarrow$ 5.~(note that if $S$ is compact, then $\varphi S$ is also compact for all $\varphi \in \Cont_c(X)$).
\end{proof}

\begin{proposition}\label{Proposition:CompactOperatorIfVanishesAtInfinity} Let $E$ and $F$ be Banach $B$-pairs and let $T\colon E\to F$ be a continuous field of operators. Then $T$ is compact if and only if $T$ is locally compact and $x\mapsto \norm{T_x}$ vanishes at infinity.
\end{proposition}
\begin{proof}
Let $T$ be compact. It is clear from the definitions that $T$ is locally compact. Moreover, we have already noted that $x\mapsto \norm{T_x}$ vanishes at infinity.

Conversely, let $T$ be locally compact and let $x\mapsto \norm{T_x}$ vanish at infinity. Let $\varepsilon>0$. Find a compact set $K\subseteq X$ such that $\norm{T_x}\leq \varepsilon$ for all $x\notin K$. Find a function $\chi\in \Cont_c(X)$ such that $0\leq \chi\leq 1$ and $\chi\equiv 1$ on $K$. Find a compact operator $S\in \Komp_B(E,F)$ such that $\norm{T_l-S_l}\leq \varepsilon$ for all $l\in \supp \chi$ (using the above characterisation of local compactness). Then also $\norm{T_l -(\varphi S)_l}\leq \varepsilon$ for all $l\in \supp \varphi$ and $T_x = (\varphi S)_x =0$ for all $x\notin \supp \varphi$. Hence $\norm{T-\varphi S}\leq \varepsilon$. So $T$ can be approximated by compact operators and is therefore compact.
\end{proof}

\begin{lemma} Let $E_1$, $E_2$ and $E_3$ be Banach $B$-pairs. Then we have $\Lin_B(E_2,E_3) \circ \Komp_B(E_1,E_2) \subseteq \Komp_B(E_1,E_3)$ and $\Komp_B(E_2,E_3) \circ \Lin_B(E_1,E_2) \subseteq \Komp_B(E_1,E_3)$.
\end{lemma}
\begin{proof}
The composition of a compact operator and a bounded linear operator is locally compact and vanishes at infinity. Hence it is compact. One can also easily prove this by direct calculation.
\end{proof}

\subsection{Operators of the form $T\otimes 1$}

\begin{definition} Let $E$ and $E'$ be Banach $B$-pairs and $F$ a Banach $B$-$C$-pair. For all operators $T\in \Linloc_B(E,E')$, define $T\otimes 1\in \Linloc_C\left(E\otimes_B F,\ E'\otimes_B F\right)$ as $\left(1\otimes T^<,\ T^>\otimes 1\right)$.
\end{definition}

\noindent The assignment $T\mapsto T\otimes 1$ is linear and functorial, and if $T$ is bounded, then $\norm{T\otimes 1} \leq \norm{T}$. The following proposition generalises Proposition 4.2 of \cite{Paravicini:07:Morita:erschienen} and is proved in \cite{Paravicini:07}, Proposition 3.1.59.

\begin{proposition}\label{Proposition:FieldsTensorProductsOfLocallyCompactOperators} Let $E$ and $E'$ be Banach $B$-pairs and $F$ a Banach $B$-$C$-pair. Assume that $\ContSect(X,B)$ acts on $F$ by locally compact operators, and call the action $\pi\colon \ContSect(X,B) \to \Komploc_C(F)$. Assume, moreover, that $E$ or $E'$ is non-degenerate. Then
\[
T\in \Komploc_B(E,E')\ \Rightarrow \ T\otimes 1 \in \Komploc_C\left(E\otimes_B F,\ E' \otimes_B F\right)
\]
and
\[
T\in \Komp_B(E,E')\ \Rightarrow \ T\otimes 1 \in \Komp_C\left(E\otimes_B F,\ E' \otimes_B F\right).
\]
\end{proposition}

There is a variant of this proposition for the pushforward:

\begin{proposition}\label{Proposition:FieldsDirectImageOfLocallyCompactOperators} Let $B'$ be another u.s.c.~field of Banach algebras and $\psi\colon B\to B'$ a continuous field of homomorphisms. Let $E$ and $F$ be Banach $B$-pairs. For all operators $T\in \Komploc_B(E,F)$, the operator $\psi_*(T)=T\otimes 1$ is contained in $\Komploc_C\left(\psi_*(E),\ \psi_*(F)\right)$ (and the same is true for compact operators instead of locally compact operators).
\end{proposition}

\subsection{The pullback}

The pullback construction that we recall and detail here is an important technical tool which is used for the definition of groupoids actions on Banach spaces and Banach algebras.

In this section, let $X$ and $Y$ be locally compact Hausdorff spaces and let $p\colon Y\to X$ be continuous.

\subsubsection{The pullback of fields of Banach spaces and Banach algebras}

\begin{definition}[The pullback]\footnote{See \cite{Lafforgue:06}, page 3.} Let $E$ be a u.s.c.~field of Banach spaces over $X$. Then we define a u.s.c.~field $p^*(E)=p^*E$ of Banach spaces over $Y$ as follows: The underlying family of Banach spaces is $(E_{p(y)})_{y\in Y}$ and the set of sections is the subset of $\prod_{y\in Y} E_{p(y)}$ of elements which are locally approximable by elements of $\{\xi \circ p:\ \xi \in \ContSect(X,E)\}$.
\end{definition}

\noindent If $E$ and $F$ are u.s.c.~fields of Banach spaces over $X$ and  $T$ is a continuous field of linear maps from $E$ to $F$, then we define
\[
p^*(T)_y:= T_{p(y)} \in \Lin\left(E_{p(y)},\ F_{p(y)}\right)
\]
for all $y\in Y$. Then $p^*(T)=p^*T$ is a continuous field of linear maps from $p^*(E)$ to $p^*(F)$. If $T$ is bounded, then so is $p^*T$ with $\norm{p^*T}\leq \norm{T}$. The assignment $T\mapsto p^*T$ is a functor from the category of fields of Banach spaces over $X$ to the category of fields of Banach spaces over $Y$.

If $E$ and $F$ are u.s.c.~fields of Banach spaces over $X$, then the internal product $p^*(E) \times_Y p^*(F)$ and $p^*\left(E\times_X F\right)$ are identical.

If $E$, $F$, $G$ are u.s.c.~fields of Banach spaces over $X$ and if $\mu$ is a continuous field of bilinear maps from $E\times_X F$ to $G$, then the family  $p^*(\mu):=(\mu_{p(y)})_{y\in Y}$ is a continuous field of bilinear maps from $p^*(E) \times_Y p^*(F)=p^*(E\times_X F)$ to $p^*(G)$. If $\mu$ is bounded, then so is $p^*\mu$  with $\norm{p^*\mu} \leq \norm{\mu}$, and if $\mu$ is non-degenerate (i.e., the image of $\mu_y$ spans a dense subset of $G_y$ for all $y\in Y$), then $p^*\mu$ is non-degenerate as well.

Note that the pullback of fields is compatible with the tensor product. Moreover, the pullback respects all kind of associativity of bilinear maps. Therefore we can pull back algebras and modules and obtain algebras and modules again:

\begin{definition}[The pullback of a field of Banach algebras] Let $A$ be a field of Banach algebras over $X$ with multiplication $\mu$. Then we equip $p^*A$ with the multiplication $p^*\mu$ to give a field of Banach algebras over $Y$. If $A$ is non-degenerate, then $p^*A$ is non-degenerate as well.
\end{definition}

\noindent Let $A$ and $B$ be fields of Banach algebras over $X$ and $\varphi\colon A\to B$ a homomorphism. Then $p^*\varphi$ is a homomorphism of fields of Banach algebras from $p^*A$ to $p^*B$, and this defines a functor from the category of fields of Banach algebras over $X$ to those over $Y$. 

Similar constructions can be done for Banach modules: if $E$ is a left Banach $A$-module then $p^*E$ is a left Banach $p^*A$-module in an obvious way. The pullback of non-degenerate modules will then be non-degenerate, the pullback of balanced bilinear maps will be balanced, and the pullback of a balanced tensor product will be a balanced tensor product over the pulled-back Banach algebra. Moreover, the pullback commutes with the direct image of Banach modules.

\subsubsection{The pullback of fields of Banach pairs}

Let $B$ be a field of Banach algebras over $X$.

\begin{definition}[The pullback of a field of Banach pairs] Let $E=(E^<,E^>)$ be a Banach $B$-pair. Then $p^*E:=\left(p^*(E^<),\ p^*(E^>)\right)$ is a Banach $p^*B$-pair when equipped with the obvious bracket.
\end{definition}

\noindent This defines a functor from the category of Banach $B$-pairs to the category of Banach $p^*B$-pairs, linear and contractive on the spaces of linear operators. As for Banach modules, the pullback of a homomorphism is a homomorphism and the pullback commutes with the tensor product and the direct image.

We now study how the pullback and locally compact operators are related (considering that the norm of compact operators vanishes at infinity, it is quite obvious that the pullback of a compact operator need not be compact).

\begin{lemma} Let $E$ and $F$ be Banach $B$-pairs. If $\eta^>\in \ContSect(X,F^>)$ and $\xi^<\in \ContSect(X,E^<)$, then
\[
\left(p^*\ketbra{\eta^>}{\xi^<}\right)_y = \ketbra{\eta^>(p(y))}{\xi^<(p(y))}
\]
for all $y\in Y$.
\end{lemma}

\begin{proposition}\label{Proposition:ContFieldsOfPairsCompactOperatorsPullBack}  Let $E$ and $F$ be Banach $B$-pairs and let $T\in \Linloc_B(E,F)$. If $T$ is locally compact, then so is $p^*T\colon p^*E \to p^*F$. Conversely, every operator $\tilde{T} \in \Komploc_{p^*B}\left(p^*E,\ p^*F\right)$ can be locally approximated by operators of the form $p^*T$ with $T\in \Komploc_B(E,F)$.
\end{proposition}
\begin{proof}
Let $T$ be locally compact. Let $y_0\in Y$. Let $\varepsilon>0$. Find a neighbourhood $U$ of $x_0:= p(y_0)$ in $X$ and $n\in \N$ and $\xi_1^<,\ldots,\xi_n^<\in \ContSect(X,E^<)$, $\eta_1^>,\ldots,\eta_n^>\in \ContSect(X,F^>)$ such that
\[
\norm{T_u -\sum_{i=1}^n \ketbra{\eta_i^>(u)}{\xi_i^<(u)}} \leq \varepsilon
\]
for all $u\in U$. Let $V:= p^{-1}(U)$. Then $V$ is a neighbourhood of $y_0$ in $Y$. For all $i\in \{1,\ldots,n\}$, the sections $\xi_i^< \circ p$ and $\eta_i^>\circ p$ belong to $\ContSect(Y,p^*E^<)$ and $\ContSect(Y,p^*F^>)$, respectively. Let $v\in V$ and define $u:=p(v) \in U$. Then
\[
\norm{p^*(T)_v - \sum_{i=1}^n \ketbra{\eta_i^>(p(v))}{\xi_i^<(p(v))}} = \norm{T_u - \sum_{i=1}^n \ketbra{\eta_i^>(u)}{\xi_i^<(u)}} \leq \varepsilon.
\]
Hence $p^*(T)$ is locally compact.

Now let $\tilde{T} \in \Komploc_{p^*B}\left(p^*E,\ p^*F\right)$. Without loss of generality we can assume that $\tilde{T}$ is of the form $\ketbra{\tilde{\eta}^>}{\tilde{\xi}^<}$ with $\tilde{\eta}^>\in \ContSect(Y, p^*F^>)$ and $\tilde{\xi}^<\in \ContSect(Y, p^*E^<)$. Let $y_0\in Y$ and $\varepsilon>0$. Find a neighbourhood $V_{\eta}$ of $y_0$ in $Y$ such that $\tilde{\eta}^>$ is bounded on $V_{\eta}$ by some constant $C_{\eta}>0$. Find an analogous neighbourhood $V_{\xi}$ for $\tilde{\xi}^<$ and a constant $C_{\xi}>0$. Find a neighbourhood $V$ contained in $V_{\eta}\cap V_{\xi}$ and $\eta^>\in \ContSect(X,F^>)$, $\xi^<\in \ContSect(X,E^<)$ such that $\|\tilde{\eta}^>(v) - \eta^>(p(v))\| \leq \varepsilon / (3C_{\eta})$ and $\|\tilde{\xi}^<(v) - \xi^<(p(v))\| \leq \varepsilon / (3C_{\xi})$ and $\|\tilde{\eta}^>(v) - \eta^>(p(v))\| \ \|\tilde{\xi}^<(v) - \xi^<(p(v))\| \leq \varepsilon/3$ for all $v\in V$. Then
\begin{eqnarray*}
\norm{\ \ketbra{\tilde{\eta}^>(v)}{\tilde{\xi}^<(v)} - \ketbra{\eta^>(p(v))}{\xi^<(p(v))}\ } &\leq & \varepsilon
\end{eqnarray*}
for all $v\in V$.
\end{proof}

\section{$\KK$-theory for Banach algebras and groupoids}

In this section, we recall the basics of $\gG$-equivariant $\KKban$-theory as introduced in \cite{Lafforgue:06}, where $\gG$ is a locally compact Hausdorff groupoid with unit space $X$. We proceed in a systematic fashion using the language that we have introduced above. In addition, we briefly discuss a sufficient condition for homotopy of $\KKbanW{\gG}$-cycles and a consequence of it, the Morita invariance of $\KKbanW{\gG}$-theory in the second component; these results are in complete analogy to the corresponding results for groups instead of groupoids which can be found in \cite{Paravicini:07:Morita:erschienen}.

\subsection{$\gG$-Banach spaces}

\begin{definition}[$\gG$-Banach space] A $\gG$-Banach space $E$ is a u.s.c.~field $E$ of Banach spaces over $\gG^{(0)}$ together with an isometric isomorphism $\alpha\colon s^*E \to r^*E$ such that
\begin{enumerate}
    \item $\forall g\in \gG^{(0)}:\ \alpha_{g} = \id_{E_g}$;
    \item $\forall (\gamma,\gamma') \in \gG*\gG:\ \alpha_{\gamma\circ \gamma'} = \alpha_\gamma \circ \alpha_{\gamma'}$;
    \item $\forall \gamma\in \gG:\ \alpha_{\gamma^{-1}} = \alpha_\gamma^{-1}$.
\end{enumerate}
\end{definition}

\noindent The Axioms 1.~and 3.~follow from Axiom 2. They are just stated to give a clearer impression of what a $\gG$-Banach space is. For notational convenience, we also write $\gamma e$ for $\alpha_{\gamma} (e)$ if $\gamma\in \gG$ and $e\in E_{s(\gamma)}$.

\begin{Xample} If we regard the locally compact space $X$ as a groupoid with unit space $X$, then every u.s.c.~field of Banach spaces over $X$ is, canonically, an $X$-Banach space (and every $X$-Banach space is, trivially, a u.s.c.~field over $X$).
\end{Xample}

\begin{definition}[$\gG$-equivariant fields of linear maps] Let $E$ and $F$ be $\gG$-Banach spaces with actions $\alpha$ and $\beta$, respectively. A \demph{$\gG$-equivariant continuous field of linear maps from $E$ to $F$} is a continuous field $(T_x)_{x\in X}$ of linear maps from $E$ to $F$ such that the following diagram commutes
\[
\xymatrix{s^*E \ar[rr]^-{s^*T} \ar[d]^{\alpha} && s^*F \ar[d]^{\beta} \\
r^*E \ar[rr]^-{r^*T} && r^*F
}
\]
This means that $T_{r(\gamma)} \circ \alpha_{\gamma} = \beta_{\gamma} \circ T_{s(\gamma)}$ for all $\gamma \in \gG$.
\end{definition}

\begin{definition}[The product and the sum of $\gG$-Banach spaces] Let $E$ and $F$ be $\gG$-Banach spaces with actions $\alpha$ and $\beta$, respectively. Then $r^*(E\times_X F) = r^*E \times_{\gG} r^*F$ and $s^*(E\times_X F) = s^*E \times_{\gG} s^*F$. We hence get a continuous field of isomorphisms $\alpha \times_{\gG} \beta\colon s^*(E \times_X F) \to r^*(E\times_X F)$. It is an action on $E\times_X F$ which we call the \demph{product action} of $\alpha$ and $\beta$. Similarly, we define an action $\alpha \oplus_{\gG} \beta$ on $E\oplus_X F$.
\end{definition}

\begin{definition}[Equivariant bilinear maps between $\gG$-Banach spaces] Let $E_1$, $E_2$ and $F$ be $\gG$-Banach spaces with $\gG$-actions $\alpha_1$, $\alpha_2$ and $\beta$, respectively. Let $\mu \colon E_1\times_X E_2 \to F$ be a continuous field of bilinear maps. Then $\mu$ is called \demph{$\gG$-equivariant} if the following diagram commutes
\[
\xymatrix{
s^*(E_1 \times_X E_2) \ar[rr]^-{s^*(\mu)} \ar[d]^{\alpha_1 \times_{\gG} \alpha_2} && s^*(F) \ar[d]^{\beta} \\
r^*(E_1\times_X E_2) \ar[rr]^-{r^*(\mu)} && r^*(F)
}
\]
This means that $\gamma \mu_{s(\gamma)}\left(e_1,e_2\right) = \mu_{r(\gamma)}\left(\gamma e_1,\ \gamma e_2\right)$ for all $\gamma \in \gG$ and $e_1\in (E_1)_{s(\gamma)}$ and $e_2\in (E_2)_{s(\gamma)}$.
\end{definition}

\begin{definition}[The tensor product of $\gG$-Banach spaces] Let $E$ and $F$ be $\gG$-Banach spaces with actions $\alpha$ and $\beta$, respectively. Then we can form the tensor product $E\otimes_X F$ of the continuous fields of Banach spaces $E$ and $F$. Now
\[
s^*\left( E\otimes_X F\right) = s^*(E) \otimes_{\gG} s^*(F) \LazyAnd r^*\left( E\otimes_X F\right) = r^*(E)
\otimes_{\gG} r^*(F).
\]
Now $\alpha \otimes \beta$ is a continuous field of isometric isomorphisms from $s^*(E) \otimes_{\gG} s^*(F)$ to
$r^*(E) \otimes_{\gG} r^*(F)$. This induces on $E\otimes_X F$ the structure of a $\gG$-Banach spaces.
\end{definition}

\noindent Note that $E\otimes_X F$ has the universal property for $\gG$-equivariant continuous fields of bilinear maps.

\begin{definition}[The trivial $\gG$-Banach space] Let $\C_X$ denote the constant field of Banach spaces over $X$ with fibre $\C$. Note that $s^*(\C_X) = \C_{\gG} = r^*(\C_X)$. So $\C_X$ is a $\gG$-Banach space if we take
$\left(\id_{\C}\right)_{\gamma\in \gG}$ as the action of $\gG$.
\end{definition}

\subsection{$\gG$-Banach algebras, $\gG$-Banach modules and $\gG$-Banach pairs}

\begin{definition}[$\gG$-Banach algebra] A $\gG$-Banach algebra $A$ is a u.s.c.~field $A$ of Banach algebras over $\gG^{(0)}$ together with a continuous field of isometric Banach algebra isomorphisms between the continuous fields of Banach algebras $s^*A$ and $r^*A$ which makes $A$ a $\gG$-Banach space.
\end{definition}

\noindent If $A$ and $B$ are $\gG$-Banach algebras, then a $\gG$-equivariant homomorphism from $A$ to $B$ is, by definition, a homomorphism of fields of Banach algebras over $\gG^{(0)}$ which is at the same time a $\gG$-equivariant continuous field of linear maps.

For the rest of this paragraph, let $B$ be a $\gG$-Banach algebra with $\gG$-action $\alpha\colon s^*B \to r^*B$.

\begin{definition}[$\gG$-Banach module] A right $\gG$-Banach $B$-module $E$ is a right Banach module $E$ over the u.s.c.~field $B$ of Banach algebras over $\gG^{(0)}$ together with a continuous field of isometric isomorphisms $\alpha^E \colon s^*E \to r^*E$ with coefficient map $\alpha$ between the Banach $s^*B$-module $s^*E$ and the Banach $r^*B$-module $r^*E$ which makes $E$ a $\gG$-Banach space.
\end{definition}

\noindent Analogously, one defines left $\gG$-Banach modules and $\gG$-Banach bimodules. There is an obvious definition of $\gG$-equivariant homomorphisms with coefficient maps between $\gG$-Banach modules and $\gG$-Banach bimodules. The balanced tensor product of $\gG$-Banach modules is defined analogously to the tensor product of $\gG$-Banach spaces, using that the balanced tensor product commutes with the pullback along $r$ and $s$. Similarly, the pushout along a continuous equivariant field of homomorphisms of Banach algebras is defined.

\begin{definition}[$\gG$-Banach $B$-pair] A $\gG$-Banach $B$-pair $E$ is a Banach $B$-pair $E=(E^<,E^>)$ together with an isometric isomorphisms $\alpha^E \colon s^*E \to r^*E$ with coefficient map $\alpha$ between the Banach $s^*B$-pair $s^*E$ and the Banach $r^*B$-pair $r^*E$ which makes $E^<$ and $E^>$ into $\gG$-Banach spaces.
\end{definition}

\noindent There is a more familiar but less systematic way to express what a $\gG$-Banach $B$-pair is: Assume that $E$ is a Banach $B$-pair. Then one can first turn the $s^*B$-pair $s^*E$ into an $r^*B$-pair by the use of the isomorphism $s^*B \cong r^*B$ which defines the $\gG$-action on $B$. Then the above isometric isomorphism $\alpha^E$ defining a $\gG$-action on $E$ can be interpreted as an element of $\Lin_{r^*B}(s^*E, r^*E)$; to turn the concurrent homomorphism $\alpha^E$ into an operator you have to invert the left-hand part, i.e., you have to substitute $\alpha_E=(\alpha_E^<,\alpha_E^>)$ with $V_E:=((\alpha_E^<)^{-1},\ \alpha_E^>)$. Then this operator $V_E$ is invertible and isometric. In the case of Hilbert modules considered by Le~Gall in \cite{LeGall:99} this operator $V_E$ is the unitary defining the groupoid action.

Although the operator notation is less systematic, we will use it for notational convenience and assume that $\gG$-actions on Banach $B$-pairs $E$ are given by ``unitaries'' $V_E$.

As above, one defines $\gG$-equivariant homomorphisms with coefficient maps. The definitions of the balanced equivariant tensor product of $\gG$-Banach pairs and the definition and properties of the pushout are straightforward.

If $E$ is a $\gG$-Banach $A$-$B$-pair, then the action of $A$ on $E$ regarded as a homomorphism from $A$ to $\Lin_B(E)$, is $\gG$-equivariant in the following sense:

\begin{lemma}\label{LemmaLeftActionEquivariantBP} Let $E$ be a $\gG$-Banach $A$-$B$-pair with $A$ and $B$ being $\gG$-Banach algebras. Let $\tilde{a} \in \ContSect(\gG, s^*A)$. Then
\[
V_E \circ \pi_{s^*A}(\tilde{a}) \circ V_E^{-1} = \pi_{r^*A}\left(\alpha^A \circ \tilde{a}\right)
\]
where $\pi_{s^*A}$ and $\pi_{r^*A}$ are the actions of $s^*A$ on $s^*E$ and $r^*A$ on $r^*E$ (regarded as homomorphisms into the linear operators) and $\alpha^A$ is the action of $\gG$ on $A$.
\end{lemma}

\subsection{$\KK$-theory for $\gG$-Banach algebras}

\subsubsection{$\KKbanW{\gG}$-cycles}

Let $E$ be a $\gG$-Banach space. Then a \demph{grading automorphism} $\sigma_E$ of $E$ is a $\gG$-equivariant contractive continuous field of linear maps from $E$ to $E$ such that $\sigma_E^2 = \id_E$. A $\gG$-Banach space endowed with a grading automorphism is called a \demph{graded $\gG$-Banach space} (compare the graded formalism of \cite{Kasparov:80}).

Just as for gradings of ordinary Banach spaces or Banach spaces with group actions we can define the notions of graded ( =even) and odd $\gG$-equivariant continuous fields of linear maps between graded $\gG$-Banach spaces, graded $\gG$-Banach algebras, graded $\gG$-Banach modules and graded $\gG$-Banach pairs. All the above constructions are compatible with this additional structure, e.g., the tensor product.

We will usually assume that our Banach algebras are trivially graded, but gradings are nevertheless important for the definition of $\KKban$-cycles.

For the rest of this paragraph, let $A$ and $B$ be (trivially graded) $\gG$-Banach algebras.

Recall from Section~\ref{Subsection:FieldsOfBanachModules} that a Banach $B$-module $E$ is called non-degenerate if $E_xB_x$ is dense in $E_x$ for all $x\in \gG^{(0)}$.

\begin{definition}[$\KKbanW{\gG}$-cycle]\label{DefinitionFieldGroupoidKasparovCycle}
A \demph{$\KKban$-cycle} from $A$ to $B$ is a pair $(E,T)$ such that $E$ is a non-degenerate graded $\gG$-$A$-$B$-bimodule and $T$ is an odd element of $\Lin_B(E)$ such that
\[
\left[\pi_A(a), T\right],\ \pi_A(a)\ (\id-T^2) \in \Komploc_B(E)
\]
for all $a\in \ContSect(X,A)$ and
\[
\pi(\tilde{a}) \left(V_E \circ s^*T \circ V_E^{-1}  - r^*T\right) \in \Komploc_{r^*B} \left(r^*E\right)
\]
for all $\tilde{a}\in \ContSect\left(\gG, r^*A\right)$, where $V_E \in \Lin_{r^*B}(s^*E,r^*E)$ denotes the ``unitary'' operator defining the action of $\gG$ on $E$. We write $\EbanW{\gG}(A,B)$ for the class of all $\KKbanW{\gG}$-cycles from $A$ to $B$.
\end{definition}

\noindent In this definition, we have used locally compact operators. Because the underlying space $X$ is locally compact Hausdorff, we can actually use compact operators instead. More precisely, we have the following characterisation of $\KKban$-cycles:

\begin{proposition}\label{Proposition:FieldsOfKasparovCyclesAlsoWithCompactOperators} A pair $(E,T)$ such that $E$ is a non-degenerate graded $\gG$-$A$-$B$-bimodule and $T$ is an odd element of $\Lin_B(E)$ is an element of $\EbanW{\gG}(A,B)$ if and only if
\[
\left[\pi_A(a), T\right],\ \pi_A(a)\ (\id-T^2) \in \Komp_B(E)
\]
for all $a\in \ContSect_0(X,A)$ and
\[
\pi(\tilde{a}) \left(V_E \circ s^*T \circ V_E^{-1}  - r^*T\right) \in \Komp_{r^*B} \left(r^*E\right)
\]
for all $\tilde{a}\in \ContSect_0\left(\gG, r^*A\right)$.
\end{proposition}
\begin{proof}
If $(E,T)$ is a $\KKban$-cycle, then we know that $\left[\pi_A(a), T\right]$ is locally compact for all $a\in \ContSect(X,A)$. In particular, this is true if $a\in \ContSect_0(X,A)$.  Since $T$ is bounded and $x\mapsto \norm{\pi_A(a)_x}$ vanishes at infinity, also $x\mapsto \left[\pi_A(a), T\right]_x$ vanishes at infinity. So $\left[\pi_A(a), T\right]$ is compact by Proposition~\ref{Proposition:CompactOperatorIfVanishesAtInfinity}. The same argument works for the other operators which have to be shown to be compact.

For the other direction use 5.\ of Proposition~\ref{Proposition:Characterisation:LocallyCompactOperator}.
\end{proof}

\noindent If $(E_1,T_1)$ and $(E_2,T_2)$ are elements of $\EbanW{\gG}(A,B)$, then we define the direct sum of cycles $(E_1, T_1) \oplus (E_2,T_2) := \left(E_1 \oplus E_2, T_1 \oplus T_2\right) \in \EbanW{\gG}(A,B)$. And if $(E,T)$ is in $\EbanW{\gG}(A,B)$, then we define $-(E,T)$ to be $(E,T)$, but equipped with the opposite grading. This is also an element of $\EbanW{\gG}(A,B)$.

Using the facts that the pushout of locally compact operators is again locally compact and that the pullback commutes with the pushout, we can define the pushout for cycles: Let $B'$ be another $\gG$-Banach algebra and $\psi\colon B\to B'$ a $\gG$-equivariant homomorphism from $B$ to $B'$. Let $(E,T)$ be an element of $\EbanW{\gG}(A,B)$. Then the pushout $\psi_*(E,T)$ of $(E,T)$ along $\psi$ is defined as $(\psi_*(E), T\otimes 1)$.  It is contained in $\EbanW{\gG}(A,B')$.

\subsubsection{Morphisms between $\KKbanW{\gG}$-cycles}

Let $A$, $A'$ and $B$, $B'$ be $\gG$-Banach algebras. Let $\varphi\colon A\to A'$ and $\psi\colon B\to B'$ be $\gG$-equivariant homomorphisms.

\begin{definition}[Morphism between $\KKbanW{\gG}$-cycles] Let $(E,T)$ and $(E',T')$ be elements of $\EbanW{\gG}(A,B)$ and $\EbanW{\gG}(A',B')$, respectively. Then a \demph{morphism} from $(E,T)$ to $(E',T')$ with coefficient maps $\varphi$ and $\psi$ is a pair $\Phi=(\Phi^<,\Phi^>)$ such that
\begin{itemize}
\item $(\Phi^<,\Phi^>)$ is an equiv.~homomorphism of graded Banach pairs with coefficient maps $\varphi$ and $\psi$;
\item we have
\[
T'^< \circ \Phi^< = \Phi^< \circ T^< \LazyAnd T'^> \circ \Phi^> = \Phi^> \circ T^>.
\]
\end{itemize}
\end{definition}

\noindent The class $\EbanW{\gG}(A,B)$, together with the morphisms of cycles (with $\id_A$ and $\id_B$ as coefficient maps), forms a category. This gives us an obvious notion of \emph{isomorphic $\KKban$-cycles} in $\EbanW{\gG}(A,B)$. Just as for ordinary $\KKban$-cycles, the sum of cycles is associative and the pushout is functorial up to isomorphism.

\subsubsection{The $\gG$-Banach algebra $B[0,1]$}\label{Subsubsection:ThegGBABNullEins}

 \begin{definition}[{The $\gG$-Banach space $E[0,1]$}]\label{DefinitionTrivialBundleOfFieldsOfBanachSpaces} Let $E$ be a $\gG$-Banach space with $\gG$-action $\alpha\colon s^*E \to r^*E$. Then we define the $\gG$-Banach space $E[0,1]$ by the following data:
\begin{enumerate}
    \item the underlying family of Banach spaces is $\left(E_x[0,1]\right)_{x\in X}$;
    \item a section $\xi$ of $E[0,1]$ is continuous if and only if $(x,t)\mapsto \xi(x)(t)$ is a continuous
    section in $p_1^*(E)$, where $p_1\colon X\times [0,1] \to X$ denotes the projection onto the first component;
    \item the action $\alpha[0,1] \colon s^*(E[0,1]) \to r^*(E[0,1])$ is defined by
    \[
     E[0,1]_{s(\gamma)}= E_{s(\gamma)}[0,1]\ \ni \ \xi_\gamma  \quad \mapsto \quad \left(t\mapsto
    \alpha_{\gamma}(\xi_{\gamma}(t))\right) \ \in \ E_{r(\gamma)}[0,1].
    \]
\end{enumerate}
For all $t\in [0,1]$, define the continuous family of linear contractions $\ev_t\colon E[0,1] \to E$ given by $(\ev_t)_x \colon E_x[0,1] \to E_x, \ \xi_x \mapsto \xi_x(t)$ for all $x\in X$.
\end{definition}

\begin{proposition} If $B$ is a $\gG$-Banach algebra, then $B[0,1]$ is a $\gG$-Banach algebra as well (when equipped with the obvious multiplication). For every $t\in [0,1]$, the field $\ev_t\colon B[0,1] \to B$ is a continuous field of homomorphisms in this case. Similar statements hold for Banach modules and pairs.
\end{proposition}

\noindent Note that $\left(\ev_{t,*} E\right)_x = (\ev_t)_{x,*}E_x$  for all $t\in [0,1]$ and $x\in X$ and all $\gG$-Banach $B[0,1]$-pairs $E$.

\subsubsection{Homotopies between $\KKban$-cycles}

Let $A$, $B$ be $\gG$-Banach algebras.

\begin{definition}[Homotopies] A \demph{homotopy} between cycles $(E_0,T_0)$ and $(E_1,T_1)$ in $\EbanW{\gG}(A,B)$ is a cycle $(E,T)$ in $\EbanW{\gG}(A,B[0,1])$ such that $\ev_{0,*}(E,T)$ is isomorphic to $(E_0,T_0)$ and $\ev_{1,*}(E,T)$ is isomorphic to $(E_1,T_1)$. If such a homotopy exists, then $(E_0,T_0)$ and $(E_1,T_1)$ are called \demph{homotopic}. We will denote by $\sim$ the equivalence relation on $\EbanW{\gG}(A,B[0,1])$ generated by homotopy. The equivalence classes for $\sim$ are called homotopy classes.
\end{definition}

\begin{defprop}[$\KKbanW{\gG}(A,B)$]\label{DefProp:KKbanGruppoid} The class of all homotopy classes in $\EbanW{\gG}(A,B)$ is denoted by $\KKbanW{\gG}(A,B)$. The addition of cycles induces a law of composition on $\KKbanW{\gG}(A,B)$ making it an abelian group (at least if we restrict the cardinality of dense subsets of the involved Banach modules by some cardinality to obtain a set $\KKbanW{\gG}(A,B)$ rather than just a class). The assignment $(A,B)\mapsto \KKbanW{\gG}(A,B)$ is functorial in both variables with respect to $\gG$-equivariant continuous fields of homomorphisms of Banach algebras.
\end{defprop}

\noindent The fact that $\KKbanW{\gG}(A,B)$ has inverses should be proved by adjusting Lemme 1.2.5 in \cite{Lafforgue:02} to the situation of $\gG$-Banach algebras. The above definition is part of D\'{e}finition-Proposition 1.2.6 in \cite{Lafforgue:06}.

The following Lemma suggests itself as a generalisation of Lemme 1.2.3 in \cite{Lafforgue:02}; a proof can be found in \cite{Paravicini:07}, Lemma~3.5.11.

\begin{lemma}\label{LemmaKasparovCyclesDifferingByCompacts}
Let $(E,T) \in \EbanW{\gG}(A,B)$ and assume that $T' \in \Lin(E)$ is odd bounded operator such that $a(T-T'),\ (T-T')a \in \Komploc_B(E)$ for all $a\in \ContSect(X,A)$. Then $(E,T') \in \EbanW{\gG}(A,B)$ and there is a homotopy from $(E,T)$ to $(E,T')$.
\end{lemma}

\subsubsection{A sufficient condition for homotopy}

As mentioned in \cite{Paravicini:07:Morita:erschienen} and proved in \cite{Paravicini:07}, Section~3.7, there is a sufficient condition for the homotopy of $\KKbanW{\gG}$-cycles, the basic idea being the following: If there is a homomorphism between two cycles, then under certain conditions the mapping cylinder of this homomorphism is a homotopy between the cycles. We formulate these conditions here and refer the reader to \cite{Paravicini:07} for the proofs.

Let $\psi\colon B\to B'$ be a continuous field of homomorphisms between u.s.c.~fields of Banach algebras over $X$. Let $\Phi_{\psi} \colon E_B \to E'_{B'}$ and $\Psi_{\psi}\colon F_B \to F'_{B'}$ be contractive continuous fields of concurrent homomorphisms with coefficient map $\psi$ between u.s.c.~fields of Banach pairs over $X$.

\begin{itemize}
\item The Banach space $\Lin_{\psi}(\Phi,\Psi)$ of ``bounded linear operators'' from $\Phi_{\psi}$ to $\Psi_{\psi}$ is defined to be the set of pairs $(T,T')$ such that $T\in \Lin_B(E,F)$, $T'\in \Lin_{B'}(E',F')$  satisfying
\[
\Psi^> \circ T^> = T'^> \circ \Phi^> \LazyAnd T'^< \circ \Psi^< = \Phi^< \circ T^<.
\]
\item The Banach space $\Komp_{\psi}(\Phi,\Psi)$ is the space of pairs $(T,T') \in \Lin_B(E,F) \times \Lin_{B'}(E',F')$ such that for all $\varepsilon>0$ there is an $n\in \N$, $\xi_1^<,\ldots,\xi_n^<\in \ContSect_0(X, E^<)$ and $\eta_1^>,\ldots,\eta_n^> \in \ContSect_0(X, F^>)$ such that
\[
\norm{T_x-\sum_{i=1}^n \ketbra{\eta_i^>(x)}{\xi_i^<(x)}} \leq \varepsilon \LazyAnd \norm{T'_x-\sum_{i=1}^n
\ketbra{\Psi^>_x(\eta_i^>(x))}{\Phi^<_x(\xi_i^<(x))}}
\]
for all $x\in X$. Note that $\Komp_{\psi}(\Phi,\Psi) \subseteq \Lin_{\psi}(\Phi,\Psi)$.
\end{itemize}

\noindent Now we are ready to formulate the sufficient condition:

\begin{theorem} Let $A$ and $B$ be $\gG$-Banach algebras. Let $(E,T), (E',T')$ be elements of $\EbanW{\gG}(A,B)$. If there is a morphism $\Phi$ from $(E,T)$ to $(E',T')$ (with coefficient maps $\id_A$ and $\id_B$) such that
\begin{enumerate}
\item $\forall a\in \ContSect_0\left(X, A\right):\ \left[a,(T,T')\right] = \left([a,T],\ [a,T']\right) \in \Komp_{\id_B}(\Phi,\Phi),$

\item $\forall a\in \ContSect_0\left(X, A\right):\ a((T,T')^2-1) = \left(a(T^2-1),\ a({T'}^2-1) \right) \in \Komp_{\id_B}(\Phi,\Phi),$

\item $\forall a\in \ContSect_0\left(\gG, r^*A\right):\ a\left(\left(\alpha^{\Lin(E,F)} s^*T,\ \alpha^{\Lin(E',F')} s^*T'\right)-(r^*T,r^*T')\right) \in \Komp_{\id_{r^*B}}(r^*\Phi,r^*\Phi),$
\end{enumerate}
then $(E,T) \sim (E',T')$.
\end{theorem}

\subsubsection{Morita cycles}

Let $A$, $B$ and $C$ be non-degenerate $\gG$-Banach algebras.

A \demph{$\gG$-equivariant Morita equivalence} between $A$ and $B$ is a pair $\left({_BE^<_A},{_AE^>_B}\right)$ of $\gG$-Banach bimodules endowed with an equivariant continuous field of  bilinear maps $\langle\cdot,\cdot \rangle_B\colon E^< \times E^> \to B$ and an equivariant continuous field of bilinear maps ${_A}\langle\cdot,\cdot\rangle\colon E^>\times E^< \to A$ such that for all $x\in X$ the pair $\left(E_x^<,\ E^>_x\right)$ with the brackets $\langle\cdot,\cdot\rangle_{B,x}$ and ${_A}\langle\cdot,\cdot\rangle_x$ is a Morita equivalence between $A_x$ and $B_x$, compare \cite{Paravicini:07:Morita:erschienen}, Definition~5.1. This notion of Morita equivalence is an equivalence relation on the class of non-degenerate $\gG$-Banach algebras.

Morita equivalences are special cases of what is called a Morita cycle: A \demph{$\gG$-equivariant Morita cycle} $F$ from $A$ to $B$ is a non-degenerate $\gG$-Banach $A$-$B$-pair $F$ such that $\ContSect_0(X,A)$ acts on $F$ by compact operators, i.e., if $\pi_A\colon \ContSect_0(X,A) \to \Lin_B(F)$ is the action of $\ContSect_0(X,A)$ on $F$, then $\pi_A(\ContSect(X,A))\subseteq \Komp_B(F)$. The class of all Morita cycles from $A$ to $B$ is denoted by $\MbanW{\gG}(A,B)$. Note that the elements of $\MbanW{\gG}(A,B)$ are elements of $\EbanW{\gG}(A,B)$ with trivial operator and trivial gradings, and we also have a canonical homotopy relation on $\MbanW{\gG}(A,B)$. The quotient after this relation is called $\MoritabanW{\gG}(A,B)$.
One can show that there is an action
\[
\otimes_B\colon \KKbanW{\gG}(A,B) \times \MoritabanW{\gG}(B,C) \to \KKbanW{\gG}(A,C)
\]
which is given on the level of cycles by the simple formula $(E,T) \otimes_B F := (E \otimes_B F, T\otimes 1)$ if $(E,T) \in \EbanW{\gG}(A,B)$ and $F\in \MbanW{\gG}(B,C)$; note that we can formulate this definition without having a Banach algebras substitute of Kasparov's Technical Theorem.

As a consequence, we have:

\begin{theorem}\label{TheoremFieldsMoritaEquivalenceRespectsKKban}
Let $E$ be a $\gG$-equivariant Morita equivalence between $B$ and $C$. Then $\cdot \otimes_B [E]$ is an isomorphism from $\KKbanW{\gG}(A,B)$ to $\KKbanW{\gG}(A,C)$.
\end{theorem}

\noindent See \cite{Paravicini:07} for the details.

\section{Induction}

In this section, we show that $\KKban$ is functorial with respect to generalised morphisms. To do this, we use that every generalised morphism can be written as a product of a special kind of equivalence and of a strict morphism (see Proposition~\ref{Proposition:ZerlegeMorphismen}). We hence first show the functoriality under strict morphisms, then under the this special kind of equivalences, and in Paragraph~4.3 we combine these two constructions. This idea appeared in the C$^*$-algebraic setting in \cite{LeGall:99}, and the basic constructions evolved from \cite{Rieffel:74, Rieffel:76, Renault:80, MuhReWill:87} to name but a few.

\subsection{The pullback along strict morphisms}

Let $\gG$ and $\gH$ be locally compact Hausdorff groupoids and let $f\colon \gH\to \gG$ be a strict morphism.

\subsubsection{The pullback of $\gG$-Banach spaces along strict morphisms}

Let $E$ be an $\gG$-Banach space with action $\alpha$. Write $f_0$ for $f\restr_{\gH^{(0)}} \colon \gH^{(0)} \to \gG^{(0)}$. Then $f_0^*(E)$ is a u.s.c.~field of Banach spaces over $\gH^{(0)}$. Now $s_{\gG} \circ f = f_0 \circ s_{\gH}$ and $r_{\gG} \circ f = f_0 \circ r_{\gH}$, so
\[
s_{\gH}^*(f^*_0(E)) = (f_0 \circ s_{\gH})^*(E) = (s_{\gG} \circ f)^*(E) = f^*(s_{\gG}^*(E))
\]
and similarly for the range maps. So $f^*(\alpha)$ is a continuous field of isometric isomorphisms from
$s_{\gH}^*(f^*_0(E))$ to $r_{\gH}^*(f^*_0(E))$. It is an action of $\gH$. The $\gH$-Banach space $f_0^*(E)$ with the action $f^*(\alpha)$ is called the \demph{pullback of $E$ along $f$} and is denoted by $f^*(E)$. The pullback commutes with the tensor product: Let $F$ be another $\gG$-Banach space. Then $f^*\left(E \otimes_{\gG^{(0)}} F\right) = f^*(E) \otimes_{\gH^{(0)}} f^*(F)$ as $\gH$-Banach spaces.

The pullback also preserves equivariance, more precisely: If $T\in \Linloc (E,F)$ is $\gG$-equivariant, then $f^*T \in \Linloc(f^*E, f^*F)$ is $\gH$-equivariant. An analogous statement is true for equivariant bilinear maps.

The pullback along $f$ is a functor from the category of $\gG$-Banach spaces to the category of $\gH$-Banach spaces, linear and contractive on the sets of bounded continuous fields of linear maps, and sending equivariant continuous fields of linear maps to equivariant continuous fields.

The assignment $f\mapsto f^*$ has the expected functorial properties: Let $\gK$ be another topological groupoid and let $g\colon \gK \to \gH$ be a strict morphism. Then $(f\circ g)^*=g^* \circ f^*$ as functors from the category of $\gG$-Banach spaces to the category of $\gK$-Banach spaces. Secondly, $\id_{\gG}^*$ is the identity functor of the category of $\gG$-Banach spaces.

\subsubsection{The pullback of $\gG$-Banach algebras etc.\ along strict morphisms}

Let $B$ be a $\gG$-Banach algebra. Then $f^*B$ is an $\gH$-Banach algebra. Also, the pullback along $f$ of a $\gG$-equivariant homomorphism of Banach algebras is a $\gH$-equivariant homomorphism.

If $E$ is a $\gG$-Banach $B$-module, then $f^*E$ is an $\gH$-Banach $f^*B$-module in an obvious way. The situation is similar for $\gG$-Banach bimodules. The pullback along $f$ of a $\gG$-equivariant linear operator or of a $\gG$-equivariant homomorphism with coefficient maps is an $\gH$-equivariant linear operator or an $\gH$-equivariant homomorphism with coefficient maps, respectively.

The pullback along $f$ respects balanced equivariant bilinear maps and balanced tensor products of equivariant Banach modules.

The functor $f^*$ on Banach modules induces a functor $f^*$ from the category of $\gG$-Banach $B$-pairs to the category of $\gH$-Banach $f^*(B)$-pairs. It sends a $\gG$-Banach $B$-pair $E=(E^<,\ E^>)$ to the $\gH$-Banach $f^*B$-pair $f^*(E)=\left(f^*(E^<),\ f^*(E^>)\right)$. The ``unitary'' operator $V_{f^*E}$ defining the action of $\gH$ on $f^*E$ is given by $f^*V_E$. A ($\gG$-equivariant) $B$-linear operator $T=(T^<,T^>)$ is sent to the ($\gH$-equivariant) $f^*(B)$-linear operator $f^*(T)=(f^*(T^<),f^*(T^>))$.

One proceeds similarly for $\gG$-Banach $A$-$B$-pairs and homomorphisms with coefficient maps. The functor respects the tensor product of Banach pairs. Also the pushout of Banach pairs is preserved:

\begin{proposition}\label{PropositionDirectImageAndPullBackofgGBanachPairs}
Let $B$ be a $\gG$-Banach algebra and $E$ a $\gG$-Banach $B$-pair. Let $B'$ be another $\gG$-Banach algebra and let $\psi\colon B\to B'$ be a $\gG$-equivariant homomorphism. Then
\[
f^*\left(\psi_*(E)\right) = \left(f^*(\psi)\right)_* \left(f^*(E)\right)
\]
as $\gH$-Banach $f^*(B')$-pairs.
\end{proposition}

\subsubsection{The pullback of $\KKban$-cycles along strict morphisms}

Let $\gG$ and $\gH$ be topological groupoids over $X$ and $Y$, respectively, and let $f\colon \gH\to \gG$ be a strict morphism of topological groupoids. Let $A$ and $B$ be $\gG$-Banach algebras.

\begin{proposition} Let $(E,T) \in \EbanW{\gG}(A,B)$. Then $f^*(E,T):= \left(f^*E, f^*T\right)$ is in $\EbanW{\gH}\left(f^*A,f^*B\right)$.
\end{proposition}
\begin{proof}
We already know that $f^*E$ is a non-degenerate $\gH$-Banach $f^*A$-$f^*B$-pair. If $\sigma_E$ is the grading automorphism of $E$, then $f^*\sigma_E=(f^*\sigma_E^<, \ f^*\sigma_E^>)$ is a grading automorphism for $f^*E$. The operator $f^*T$ is odd for this grading. Let $a\in \ContSect(X,A)$. Then $a\circ f \in \ContSect\left(Y,\ f^* A\right)$. Now Proposition~\ref{Proposition:ContFieldsOfPairsCompactOperatorsPullBack} says that the pullback of locally compact operators is again locally compact, so
\[
\left[\pi(a\circ f),\ f^*T\right] = \left[f^*(\pi(a)), \ f^*T\right] = f^*\left[\pi(a), T\right] \in \Komploc_{f^*B}\left(f^*E\right).
\]
Now let $b\in \ContSect(Y,\ f^*A)$. Let $\varepsilon>0$ and $y_0\in Y$. Then we can find an $a\in \ContSect(X,A)$ and a neighbourhood $V$ of $y_0$ in $Y$ such that $\norm{T} \norm{b(v) - a(f(v))} \leq \varepsilon$ for all $v\in V$. For all $v\in V$, we have
\begin{eqnarray*}
\norm{\left[\pi(b),\ f^*T\right]_v - \left[\pi(a\circ f),\ f^*T\right]_v} &=&\norm{\left[\pi(b-a\circ f),\ f^*T\right]_v}\\ &=& \norm{\left[\pi_{A_{f(v)}}(b(v)-a(f(v))),\ T_{f(v)}\right]}\\ &\leq& \norm{T} \norm{b(v) - a(f(v))} \leq \varepsilon.
\end{eqnarray*}
So $\left[\pi(b),\ f^*T\right]$ is locally approximable by locally compact operators, so it is itself locally compact. Analogously one shows that $\pi(b)\left(\id - f^*T^2\right)$ is locally compact.

Now let $\tilde{a}\in \ContSect\left(\gG,\ r_{\gG}^*A\right)$. Then $\tilde{a} \circ f \in \ContSect\left(\gH,\ f^*r_{\gG}^*A\right) = \ContSect\left(\gH,\ r_{\gH}^*f^*A\right)$; note that $f^* r_{\gG}^* A=r_{\gH}^* f^*A$. Now
\begin{eqnarray*}
&& \pi(\tilde{a}\circ f)\left( V_{f^*E}\circ \left(s_{\gH}^*f^*T\right)\circ V_{f^*E}^{-1}- r_{\gH}^*f^*T\right)\\
&=& f^*\pi(\tilde{a}) \left((f^*V_E) \circ \left( f^*s_{\gG}^*T \right)\circ (f^*V_E^{-1}) - f^* r_{\gG}^* T\right)\\ &=& f^*\left(\pi(\tilde{a})\left(V_E \circ ( s_{\gG}^*T)\circ V_E^{-1} - r_{\gG}^* T\right) \right)\\
&\in& \Komploc_{f^*r_{\gG}^*B} \left(f^*r_{\gG}^*E\right) = \Komploc_{r_{\gH}^* f^*B} \left(r_{\gH}^*f^*E\right).
\end{eqnarray*}
As above, one can extend this to all $\tilde{b} \in \ContSect\left(\gH,\ r_{\gH}^*f^*A\right)$ (instead of $\tilde{a} \circ f$).

\noindent So $f^*(E,T) \in \EbanW{\gH}\left(f^*A,f^*B\right)$.
\end{proof}

\noindent The pullback along $f$ respects the direct sum of cycles, the pushout, and we have $f^*(B[0,1]) = (f^*B)[0,1]$. It follows that the pullback also respects homotopies. Hence we get the following theorem:

\begin{theorem}\label{TheoremPullBackOfKKAlongStrictMorphisms}
The pullback along the strict morphism $f\colon \gH\to \gG$ induces a homomorphism
\[
f^* \colon \KKbanW{\gG}\left(A,B\right) \to \KKbanW{\gH}\left(f^*A,\ f^*B\right).
\]
\end{theorem}

\noindent It is natural with respect to $\gG$-equivariant homomorphisms in both variables.

\subsection{The functor $p_!$}\label{SubsectionTheFunctorPShriek}\label{SUBSECTIONTHEFUNCTORPSHRIEK}

We now analyse the second part of the functoriality with respect to generalised morphisms. Every generalised morphism can be expressed as the product of a strict morphism and an equivalence of the following type (see Proposition~\ref{Proposition:ZerlegeMorphismen}):

Let $Y$ and $X$ be locally compact Hausdorff spaces and let $p\colon Y\to X$ be continuous, open, and surjective. Let $\gG$ be a locally compact Hausdorff groupoid over $X$. We denote the canonical strict morphism from $p^*(\gG)$ onto $\gG$ also by $p$. According to Proposition~\ref{Proposition:GraphPEquivalence}, the graph of $p$ is an equivalence between $p^*(\gG)$ and $\gG$. Also the category of $p^*(\gG)$-Banach spaces and the category of $\gG$-Banach spaces are equivalent, but here, we have to be more precise:

The space $X$ itself can be regarded as a groupoid over $X$, and we have
\[
p^*(X) \cong Y \times_X Y.
\]
The isomorphism from $p^*(X)$ to $Y \times_X Y$ sends $(y',x,y)$ to $(y',y)$, where $y,y'\in Y$, $x\in X$ and $p(y')=x=p(y)$. Note that $p^*(X) =Y\times_X Y$ is contained as a closed subgroupoid in $p^*(\gG)$.

If $E$ is a u.s.c.~field of Banach spaces over $X$, then $p^*(E)$ is not only a u.s.c.~field of Banach spaces over $Y$, but also a $Y\times_X Y$-Banach space. As a consequence, a condition on the linear operators between $p^*(\gG)$-Banach spaces which is natural in our context is $Y\times_X Y$-equivariance. Every continuous field of linear maps between $p^*(\gG)$-Banach spaces which is $p^*(\gG)$-equivariant is also $Y\times_X Y$-equivariant. In this section we show that the pullback functor $p^*$ implements the following one-to-one correspondences:
\begin{enumerate}
\item $\gG$-Banach spaces correspond to $p^*(\gG)$-Banach spaces;
\item continuous fields of linear maps between $\gG$-Banach spaces correspond to $Y\times_X Y$-equivariant continuous fields of linear maps between $p^*(\gG)$-Banach spaces;
\item $\gG$-equivariant continuous fields of linear maps correspond to $p^*(\gG)$-equivariant fields of linear maps.
\end{enumerate}
We reach this goal by defining a functor $p_!$ which inverts $p^*$; it points in the opposite direction, from the $p^*(\gG)$-Banach spaces to the $\gG$-Banach spaces. The functor $p_!$ is obtained by ``factoring out'' the $Y\times_X Y$-action on the given $p^*(\gG)$-Banach space.

For technical reasons, \textsc{we assume that there exists a faithful continuous field of measures on $Y$ over $X$ with coefficient map $p$}. As discussed on page~\pageref{Page:HaarsystemsOnYXY}, this condition is equivalent to the condition that the locally compact Hausdorff groupoid $Y\times_X Y$ admits a left Haar system. Note that such a faithful continuous field of measures on $Y$ (and hence a Haar system on $Y\times_X Y$) exists if $Y$ is second countable.\footnote{This can bee deduced from Proposition 3.9 in \cite{Blanchard:96}.} In the situation we are interested in, the space $Y$ is actually a graph $\Omega$ from $\gG$ into some other locally compact Hausdorff groupoid $\gH$. We have learned above that such an $\Omega$ carries a canonical left Haar system if $\gH$ carries a left Haar system, so the existence of a faithful continuous field of measures on $Y=\Omega$ will be automatic in this case.

The proofs of the results of this section can be found in \cite{Paravicini:07}, most of them in Section 6.5.

\subsubsection{The functor $p_!$ for $p^*(\gG)$-Banach spaces}

Let $E$ be a $p^*(\gG)$-Banach space with action $\alpha$; as $Y\times_X Y$ is contained in $p^*(\gG)$, it also acts on $E$. We define a u.s.c.~field of Banach spaces $p_!E$ over $X$ as follows: For every $x\in X$, define
\[
(p_! E)_x:= \bigg\{\left(e_y\right)_{y\in Y_x}\Big| \ \forall y,y'\in Y_x:\ e_y\in E_y \wedge
\alpha_{(y',y)}\left(e_y\right) = e_{y'}\bigg\} \subseteq \prod_{y\in Y_x} E_y,
\]
where we take the $\sup$-norm on $\prod_{y\in Y_x} E_y$. Note that $(p_! E)_x$ is a closed linear subspace of the product. Since $\alpha$ is a field of isometries, it follows that the norm of a family $(e_y)_{y\in Y_x} \in (p_! E)_x$ equals the norm of each $e_y$, $y\in Y_x$; hence $(p_! E)_x$ is isometrically isomorphic to $E_y$ for each $y\in Y_x$ (note that $Y\times_X Y$ acts freely on $Y$).

To define the structure of a u.s.c.~field of Banach spaces over $X$ on $(p_!E_x)_{x\in X}$, we set
\[
\Delta:= \Delta_E := \bigg\{\delta \in \ContSect\left(Y, E\right)\Big| \ \forall (y,y')\in Y \times_X Y: \ \alpha_{(y',y)} \left(\delta(y)\right) = \delta(y')\bigg\}.
\]
In other words: $\Delta$ consists of those sections of $E$ which are invariant under the action of $Y\times_X Y$. If $\delta\in \Delta$ and $x\in X$, then define
\[
(p_!\delta)(x):= \left(\delta(y)\right)_{y\in Y_x} \in (p_!E)_x
\]
Now
\[
\Gamma := \left\{p_!\delta: \delta \in \Delta\right\}
\]
satisfies conditions (C1)-(C4), so $(p_!E,\ \Gamma)$ is a u.s.c.~field of Banach spaces over $X$.

There is an action $p_! \alpha$ of $\gG$ on $p_!E$. It has the (defining) property that for all $\gamma\in \gG$, $e=(e_y)_{y\in Y_{s(\gamma)}}\in (p_! E)_{s(\gamma)}$, and $y\in Y_{s(\gamma)}$:
\begin{equation}\label{EquationDefiningRelationForTheActionOfG}
(p_! \alpha)_{\gamma} (e) = (\alpha_{(z,\gamma,y)} e_y)_{z\in Y_{r(\gamma)}}.
\end{equation}

Now let $E$ and $F$ be $p^*(\gG)$-Banach spaces. Let $T$ be an $Y\times_X Y$-equivariant continuous field of linear maps from $E$ to $F$. Define for all $x\in X$ and $e=(e_y)_{y\in Y_x} \in (p_!E)_x$:
\[
\left( p_! T\right)_x \left(e\right) := \left(T_y e_y\right)_{y\in Y_x} \in \left(p_! F\right)_x.
\]
Then $p_! T$ is a continuous field of linear maps from $p_! E$ to $p_! F$. If $T$ is bounded, then $\norm{p_! T} =
\norm{T}$. If $T$ is $p^*(\gG)$-equivariant, then $p_!T$ is $\gG$-equivariant.

The maps $E\mapsto p_{!}E$ and $T\mapsto p_!T$ define a functor from the category of $p^*(\gG)$-Banach spaces with the bounded $Y\times_X Y$-equivariant continuous fields of linear maps to the category of $\gG$-Banach spaces with the bounded continuous fields of linear operators, isometric and linear on the morphism sets and respecting the tensor product.

\begin{proposition}
This functor $p_!$ from the category of $p^*(\gG)$-Banach spaces to the category of $\gG$-Banach spaces is an equivalence which inverts $p^*$; more precisely:
\begin{enumerate}
\item Define for all $p^*(\gG)$-Banach spaces $E$ and all $y\in Y$ the linear map
\[
I^E_y\colon \left(p^*p_! E\right)_y = \left(p_! E\right)_{p(y)} \to E_y,\ \left(e_z \right)_{z\in Y_{p(y)}} \to e_y.
\]
Then
\[
I^E\colon p^*p_! E \cong E
\]
is a natural $p^*(\gG)$-equivariant isometric isomorphism, compatible with the tensor product (=``multiplicative'').

\item For all $\gG$-Banach spaces $E$ there is a natural multiplicative $\gG$-equivariant isometric isomorphism
\[
J^E \colon p_!p^* E \cong E.
\]
To define $J^E$, let us analyse the action of $Y \times_X Y$ on $p^*E$ and the fibres of $p_! p^*E$: The action of $Y\times_X Y$ is the pullback of the trivial action of $X$ on $E$, so for all $(z,y)\in Y\times_X Y$ we have $p^*(E)_z = E_{p(z)} = E_{p(y)} = p^*(E)_y$ and $\alpha_{(z,y)} = \id_{E_{p(y)}}$. So if $x\in X$, then the elements of $(p_! p^* E)_x$ are of the form $(e)_{y\in Y_x}$ with $e\in E_x$; so it makes sense to define
\[
J^E_x \colon \left(p_! p^*E\right)_x \to E_x,\ \left(e\right)_{y\in Y_x} \mapsto e.
\]
\end{enumerate}
\end{proposition}

\subsubsection{The functor $p_!$ for Banach algebras, etc.}

The functor $p_!$ is multiplicative and contractive on the morphism sets. The multiplicativity gives us a way to define the functor also for equivariant fields of bilinear maps. We can therefore also define $\gG$-Banach algebras $p_!A$ for $p^*(\gG)$-Banach algebras $A$ and $\gG$-equivariant homomorphisms $p_!\varphi$ for $p^*(\gG)$-equivariant homomorphisms of Banach algebras. Similarly, we can define $p_!$ for $p^*(\gG)$-Banach modules and $p^*(\gG)$-equivariant homomorphisms of Banach modules. Moreover, if $T$ is a $Y\times_X Y$-equivariant continuous field of linear operators between $p^*(\gG)$-Banach modules $E_B$ and $F_B$, then $p_!T$ is a continuous field of linear operators between $p_!E_{p_!B}$ and $p_!F_{p_!B}$ (where $B$ is some $p^*(\gG)$-Banach algebra). All this culminates in the following definition:

\begin{definition} Let $B$ be a $p^*(\gG)$-Banach algebra and let $E=(E^<,E^>)$ be a $p^*(\gG)$-Banach $B$-pair. Then $p_!E=(p_!E^<,\ p_!E^>)$ is a $\gG$-Banach $p_!B$-pair. If $F$ is another $p^*(\gG)$-Banach $B$-pair and $T\in \Linloc_B(E,F)$ is $Y\times_X Y$-equivariant, then $p_!T=(p_!T^<,\ p_!T^>)$ is in $\Linloc_{p_!B}(p_!E,p_!F)$.
\end{definition}

\noindent This defines a functor form the category of $p^*(\gG)$-Banach $B$-pairs to the category of $\gG$-Banach $p_!B$-pairs. It inverts the functor $p^*$ and respects grading automorphisms.

\medskip
\noindent As a variant of Proposition~\ref{Proposition:ContFieldsOfPairsCompactOperatorsPullBack} one proves:

\begin{proposition}\label{Proposition:PShriekAndLocallyCompactOperators} Let $B$ be a $p^*(\gG)$-Banach algebra and let $E$ and $F$ be $p^*(\gG)$-Banach $B$-pairs. If $T \in \Komploc_B(E,F)$ is $Y\times_X Y$-equivariant, then $p_!T\in \Komploc_{p_!B}(p_!E,\ p_!F)$.
\end{proposition}

It is obvious that the functor $p_!$ is compatible with the direct sum of $p^*(\gG)$-Banach spaces and of $\gG$-Banach spaces and that the same is true for Banach modules and Banach pairs. Because $p_!$ is also compatible with the (balanced) tensor product, we obtain:

\begin{proposition} Let $B$ and $C$ be $p^*(\gG)$-Banach algebras and let $\psi\colon B\to C$ be a $p^*(\gG)$-equivariant homomorphism. Let $E$ be a right $p^*(\gG)$-Banach $B$-module. Then $p_!\C_Y$ is isomorphic to $\C_X$, $p_!\tilde{C} = p_!(C \oplus_Y \C_Y)$ is isomorphic to $\widetilde{p_!C} = p_!C \oplus_X \C_X$ and, finally, $p_!(\psi_*(E)) = p_!(E\otimes_{\tilde{B}} \tilde{C})$ is isomorphic to $(p_!\psi)_*(p_! E)$.
\end{proposition}

\noindent Moreover, $p_!$ is also compatible with the construction of trivial fields over $[0,1]$; in particular, we have:

\begin{proposition} Let $B$ be $p^*(\gG)$-Banach algebra. Then $p_!(B[0,1])$ is isomorphic to $(p_!B)[0,1]$. The isomorphism in the fibre over $x\in X$ sends $(\beta_y)_{y\in Y_x} \in p_!(B[0,1])_x$ to $t\mapsto (\beta_y(t))_{y\in Y_x} \in (p_!B)[0,1]_x$.
\end{proposition}

\subsubsection{The functor $p_!$ and $\KKban$-cycles}

This section is a translation of Section 7.2 in \cite{LeGall:99} into the language of Banach algebras; in particular, the method to make the operator of a $\KKban$-cycle equivariant is borrowed from Lemma 7.1 of that article.

Let $A$ and $B$ be $p^*(\gG)$-Banach algebras. Let $\E_{p^*(\gG)}^{\ban, Y\times_X Y}(A,\, B)$ be the class of those cycles $(E,T)$ in $\EbanW{p^*(\gG)}(A,\, B)$ such that $T$ is $Y\times_X Y$-equivariant. In an obvious manner, we define $\KK_{p^*(\gG)}^{\ban, Y\times_X Y}(A,\, B)$.

\begin{proposition} Let $(E,T) \in \E_{p^*(\gG)}^{\ban, Y\times_X Y}(A,\, B)$. Then
\[
p_!(E,T):=(p_!E,\ p_!T) \in \EbanW{\gG}\left(p_! A,\, p_! B\right).
\]
\end{proposition}
\begin{proof}
Let $a\in \ContSect(X,p_!A)$. Then we can find a $\tilde{a}\in \ContSect(Y, A)$ which is invariant under the action of $Y\times_X Y$ such that $p_!\tilde{a} = a$. Now
\[
[a,p_!T] = \left[p_! \tilde{a},\ p_!T\right] = p_! \left[\tilde{a},\ T\right] \in \Komploc_{p_!B}\left(p_!E\right)
\]
where we have used Proposition~\ref{Proposition:PShriekAndLocallyCompactOperators} and the fact that the action of $a$ on $p_!E$ is $p_!$ of the action of $\tilde{a}$ on $E$. Similarly, the other two conditions are checked (here we also use that $V_{p_! E} = p_!V_E$ in an appropriate sense).
\end{proof}

\noindent Up to isomorphism of cycles, $p_!$ inverts $p^*$ as a map from $\EbanW{\gG}(p_! A,\, p_! B)$ to $\E_{p^*(\gG)}^{\ban, Y\times_X Y}(A,\, B)$. And up to isomorphism, $p_!$ commutes with the pushforward and the pullback of cycles. It also commutes with homotopies. We therefore get:

\begin{proposition} The map $p_!$ defines an isomorphism
\[
p_! \colon \KK_{p^*(\gG)}^{\ban, Y\times_X Y}\left(A,\ B\right) \cong \KKbanW{\gG}\left(p_!A,\ p_!B\right),
\]
inverting $p^*$.
\end{proposition}

\textsc{Now let $X$ be $\sigma$-compact.}

\begin{lemma}\label{LemmaMakeOperatorsEquivariantShriek} Let $(E,T) \in \EbanW{p^*(\gG)}(A,\, B)$. Then there is an odd $Y\times_X Y$-equivariant linear operator $\tilde{T}$ on $E$ such that $(E,T)$ is homotopic to $(E,\tilde{T})$ in $\EbanW{p^*(\gG)}(A,\, B)$. The construction is compatible with the pullback and hence with homotopies of cycles.
\end{lemma}
\begin{proof}
Let $\mu$ be a faithful continuous field of measures on $Y$ over $X$, which can be regarded as a Haar system of the locally compact Hausdorff groupoid $Y\times_X Y$. Note that the quotient space $Y/ (Y\times_X Y)$ is homeomorphic to the $\sigma$-compact space $X$, so it is itself $\sigma$-compact, and we can find\footnote{See, for example, Theorem~6.3 in \cite{Tu:04}.} a cut-off function $c\colon Y\to [0,\infty[$ for $Y\times_X Y$, i.e., a continuous function $c$ on $Y$ such that $\int_{y\in Y_x} c(y)\rmd \mu_x(y) =1$ for all $x\in X$ and $p^{-1}(K) \cap \supp c$ is compact for all compact $K\subseteq X$.

\noindent Define
\[
\tilde{T}_y:= \int_{y' \in Y_{p(y)}} c(y') \ \alpha_{(y,y')} T_{y'} \alpha_{(y',y)} \rmd \mu_{p(y)} (y')
\]
for all $y\in Y$, where $\alpha$ denotes the action of $p^*(\gG)$ (and hence also of $Y\times_X Y$) on $E$ (actually, the formula makes sense for the right-hand side of the pair $E$ and should be interpreted appropriately for the left-hand side). Now $\tilde{T}$ is an odd $Y\times_X Y$-equivariant bounded continuous field of linear operators on $E$. It is not hard to show that $a(T-\tilde{T})$ and $(T-\tilde{T})a$ are locally compact for all $a\in \ContSect(Y, A)$, compare Lemma 7.2.6 of \cite{Paravicini:07}. It follows from Lemma~\ref{LemmaKasparovCyclesDifferingByCompacts} that $(E,\tilde{T})$ is a cycle and that $(E,T)$ and $(E, \tilde{T})$ are homotopic.
\end{proof}

\noindent The preceding lemma implies:

\begin{proposition} The obvious homomorphism from $\KK_{p^*(\gG)}^{\ban, Y\times_X Y}(A,\, B)$ to $\KK_{p^*(\gG)}^{\ban}(A,\, B)$ is an isomorphism.
\end{proposition}

\begin{corollary}\label{CorollaryPShriekForPlainKKTheory}
$p_!$ is a well-defined isomorphism
\[
p_! \colon \KK_{p^*(\gG)}^{\ban}\left(A,\ B\right) \cong \KKbanW{\gG}\left(p_!A,\ p_!B\right),
\]
inverting $p^*$.
\end{corollary}

\subsection{The pullback along generalised morphisms}\label{Subsection:PullbackAlongGeneralisedMorphisms}

Let $\gG$ and $\gH$ be locally compact Hausdorff groupoids (with open range and source maps) carrying left Haar systems. Recall that the existence of a left Haar system on $\gH$ implies the existence of a left Haar system on each graph from $\gG$ to $\gH$.

\subsubsection{The pullback of Banach spaces along generalised morphisms}

\begin{definition} Let $\Omega$ be a graph from $\gG$ to $\gH$ with anchor maps $\rho$ and $\sigma$. Then $f_{\Omega}$ as defined in \ref{DefPropfOmega} is a strict morphism from $\rho^*(\gG)$ to $\gH$, which extends $\sigma\colon \Omega \to \gH^{(0)}$. For all $\gH$-Banach spaces $E$, define
\[
\Omega^* \left( E\right) := \rho_! f_{\Omega}^* \left(E\right).
\]
This will also be written as $\rho_! \sigma^* E$ or as $\Ind_{\gH}^{\gG} E$ and be called induction functor, especially if $\Omega$ is an equivalence.
\end{definition}

\noindent If $\Omega$ is as above, then $E\mapsto \Omega^* E$ is a functor from the category of $\gH$-Banach spaces with the $\gH$-equivariant (bounded, contractive) continuous fields of linear maps to the category of $\gG$-Banach spaces with the $\gG$-equivariant (bounded, contractive) continuous fields of linear maps. It commutes with the tensor product and has the (characterising) property that $\rho^* \Omega^* E$ is naturally isomorphic to $f_{\Omega}^*(E)$.

\begin{proposition}\label{PropositionCompositionOfGeneralisedPullBackSpaces} Let $\gK$ be another locally compact Hausdorff groupoid carrying a left Haar system. Let $\Omega$ be a graph from $\gG$ to $\gH$ and $\Omega'$ a graph from $\gH$ to $\gK$. Then
\[
\Omega^*\circ (\Omega')^* \cong \left(\Omega \times_{\gH} \Omega'\right)^*
\]
as multiplicative functors from the $\gK$-Banach spaces to the $\gG$-Banach spaces, i.e.,
\[
\Ind_{\gH}^{\gG} \circ \Ind_{\gK}^{\gH} = \Ind_{\gK}^{\gG}.
\]
\end{proposition}
\begin{proof}
Let $\rho$ and $\sigma$ be the anchor maps of $\Omega$ and $\rho'$ and $\sigma'$ those of $\Omega'$. Let $\pi_1$ and
$\pi_2$ denote the projections from $\Omega \times_{\gH^{(0)}} \Omega'$ to the first and second component. As $\rho'$
is open and surjective, so is $\pi_1$. Write $p$ for the (open and surjective) quotient map from $\Omega
\times_{\gH^{(0)}} \Omega'$ onto $\Omega'':= \Omega \times_{\gH} \Omega'$, and denote the anchor maps of $\Omega''$ by
$\rho''$ and $\sigma''$. Consider the diagram
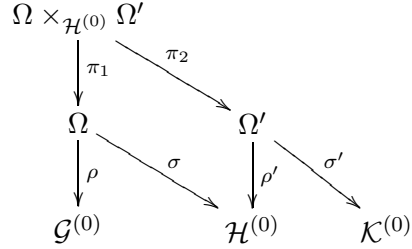
\[
\xymatrix{\Omega \times_{\gH^{(0)}} \Omega' \ar[d]^{\pi_1} \ar[rd]^{\pi_2}&&\\
\Omega \ar[d]^{\rho} \ar[rd]^{\sigma} & \Omega'\ar[d]^{\rho'} \ar[rd]^{\sigma'}&\\
\gG^{(0)} & \gH^{(0)} & \gK^{(0)}
}
\]
This is a diagram just for the unit spaces, but of course there is a corresponding commutative diagram also for the
groupoids themselves:
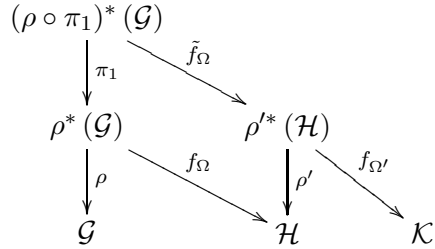
\[
\xymatrix{(\rho \circ \pi_1)^*\left(\gG\right) \ar[d]^{\pi_1} \ar[rd]^{\tilde{f}_{\Omega}}&&\\
\rho^*\left(\gG\right) \ar[d]^{\rho} \ar[rd]^{f_{\Omega}} & \rho'^*\left(\gH\right) \ar[d]^{\rho'} \ar[rd]^{f_{\Omega'}}&\\
\gG & \gH & \gK
}
\]
Here the strict morphism $\tilde{f}_{\Omega}$ is defined as follows: It sends
$((\omega_2,\omega'_2),\, \gamma,\,(\omega_1,\omega_1')) \in (\rho \circ \pi_1)^* (\gG)$ to $(\omega'_2,\,
f_{\Omega}(\omega_2,\gamma,\omega_1),\, \omega_1')$. It follows that
\[
\Omega^*\circ (\Omega')^* =  \rho_!  \circ \left(f_{\Omega}^* \circ \rho'_!\right) \circ f_{\Omega'}^* \cong \rho_!
\circ \left((\pi_1)_! \circ \tilde{f}_{\Omega}^* \right) \circ f_{\Omega'}^* \cong (\rho \circ \pi_1)_! \circ
\left(f_{\Omega'} \circ \tilde{f}_{\Omega}\right)^*.
\]
On the other hand, also the following diagrams commute
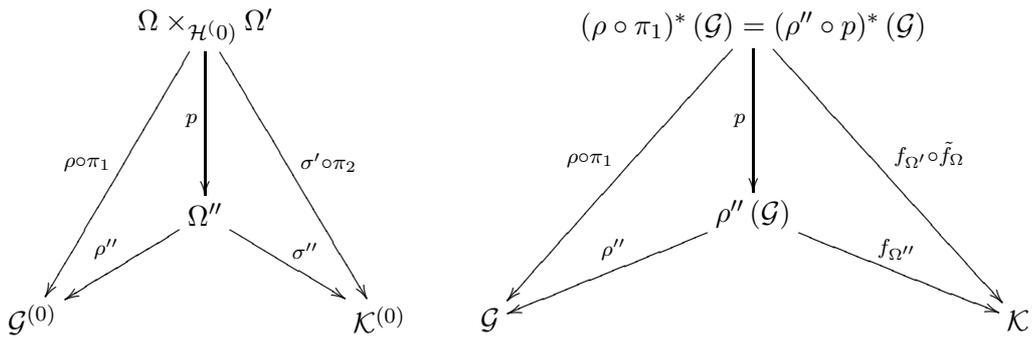
\[
\xymatrix{ &\Omega \times_{\gH^(0)} \Omega' \ar[dd]_p \ar[dddl]_{\rho \circ \pi_1} \ar[dddr]^{\sigma' \circ \pi_2} &\\
&&\\
& \Omega'' \ar[dl]_{\rho''} \ar[dr]^{\sigma''}& \\
\gG^{(0)} && \gK^{(0)}
}
\qquad
\xymatrix{ &(\rho \circ \pi_1)^*\left(\gG\right) = (\rho'' \circ p)^*\left(\gG\right) \ar[dd]_p \ar[dddl]_{\rho \circ \pi_1} \ar[dddr]^{f_{\Omega'} \circ \tilde{f}_{\Omega}} &\\
&&\\
& \rho''\left(\gG\right) \ar[dl]_{\rho''} \ar[dr]^{f_{\Omega''}}& \\
\gG && \gK
}
\]
To check that $f_{\Omega'} \circ \tilde{f}_{\Omega}= f_{\Omega''} \circ p$ let
$((\omega_2,\omega'_2),\gamma,(\omega_1,\omega_1'))$ be an element of $(\rho \circ \pi_1)^* (\gG)$. Then $f_{\Omega''} (p((\omega_2,\omega'_2),\gamma,(\omega_1,\omega_1')))$ is defined to be the unique element $\kappa \in \gK$ such that $[\omega_2,\omega'_2] \kappa = \gamma [\omega_1,\omega'_1]$. Also $f_{\Omega}(\omega_2,\gamma,\omega_1)$ is the unique element $\eta\in \gH$ such that $\omega _2 \eta =\gamma \omega_1$
and $f_{\Omega'}(\omega'_2,\eta,\omega'_1)$ is the unique element $\kappa'\in \gK$ such that $\omega'_2
\kappa' = \eta \omega'_1$. Now
\[
[\omega_2, \omega'_2] \kappa' = [\omega_2,\omega'_2\kappa'] = [\omega_2,\eta \omega'_1] = [\omega_2 \eta, \omega'_1] =
[\gamma \omega_1, \omega'_1] = \gamma [\omega_1,\omega'_1],
\]
so $\kappa = \kappa'$, which is what we wanted to verify.

So it follows that
\[
\left(\Omega''\right)^* = \rho''_! \circ f_{\Omega''}^* \cong \left((\rho \circ \pi_1)_! \circ p^* \right) \circ
\left(p_! \circ \left(f_{\Omega'} \circ \tilde{f}_{\Omega}\right)^* \right) \cong (\rho \circ \pi_1)_! \circ
\left(f_{\Omega'} \circ \tilde{f}_{\Omega}\right)^*.\qedhere
\]
\end{proof}

\begin{proposition}\label{PropositionPullBackOfMorphismUnambiguousForBanachSpaces} Let $f$ be a strict morphism from $\gG$ to $\gH$. Then $\Graph(f)^* \cong f^*$. In particular we have $\gG^* \cong \id_{\gG}^*$.
\end{proposition}
\begin{proof}
Write $\Omega$ for $\Graph(f)=\gG^{(0)}\times_{\gH^{(0)}} \gH$ and denote the anchor maps of $\Omega$ by $\rho$ and $\sigma$. Then $\rho^*(\gG) = \Omega \times_{\gG^{(0)}} \gG \times_{\gG^{(0)}} \Omega$, and $f_{\Omega} \colon \rho^*(\gG) \to \gH$ sends $(g,\eta,\gamma,g',\eta')$ to $\eta^{-1}f(\gamma)\eta'$. If $E$ is an $\gH$-Banach space and $g\in \gG^{(0)}$, then the fibre $(\Omega^* E)_g$ of $\Omega^*E$ at $g$ is, by definition, given by
\[
\bigg\{\left(e_{(g,\eta)}\right)_{(g,\eta)\in \Omega}\Big|\ \forall (g,\eta,g,g,\eta')\in \rho^*(\gG):\ e_{(g,\eta)} \in (\sigma^*E)_{(g,\eta)} \wedge e_{(g,\eta)} = (g,\eta,g,g,\eta') e_{(g,\eta')} \bigg\}.
\]
Analysing the action of $\rho^*(\gG)$ on $\sigma^*(E)$ gives $(g,\eta,\gamma,g',\eta') e = (\eta^{-1}f(\gamma)\eta')e$ for all elements $(g,\eta,\gamma,g',\eta') \in \rho^*(\gG)$ and $e\in (\sigma^*(E))_{(g',\eta')} = E_{s(\eta')}$. We can therefore simplify the above expressions:
\[
(\Omega^* E)_g=\bigg\{\left(e_{(g,\eta)}\right)_{(g,\eta)\in \Omega}\Big|\ \forall \eta,\eta' \in \gH^{f(g)} :\ e_{(g,\eta)} \in E_{s(\eta)} \wedge e_{(g,\eta)} = \eta^{-1}\eta' e_{(g,\eta')} \bigg\}.
\]
For all $g\in \gG^{(0)}$, the fibre of $f^*E$ at $g$ is simply $(f^*E)_g = E_{f(g)}$. If $e\in E_{f(g)}$, then define
\[
\Phi_g(e):= \left(\eta^{-1} e\right)_{(g,\eta)\in \Omega} \in (\Omega^*E)_g.
\]
This defines an isometric bijection between $(f^*E)_g$ and $(\Omega^*E)_g$; the inverse sends $(e_{(g,\eta)})_{(g,\eta)\in \Omega}$ to $e_{(g,f(g))} \in E_{f(g)}$. It can be shown that $\Phi$ is a $\gG$-equivariant continuous field of isometric linear maps and that this construction is compatible with the tensor product.
\end{proof}

\begin{corollary}\label{CorollaryEquivalentGroupoidsEquivalentCategoriesOfSpaces} Let $\Omega$ be an equivalence between $\gG$ and $\gH$. Then $E\mapsto \Omega^* E$ is an equivalence of the categories of $\gH$-Banach spaces and $\gG$-Banach spaces, isometric and linear on the morphism sets of equivariant bounded continuous fields of linear maps and compatible with the tensor product.
\end{corollary}

\subsubsection{The pullback of $\KKban$-cycles along generalised morphisms}

\textsc{For the rest of this chapter, assume that all the unit spaces of the appearing groupoids are $\sigma$-compact.}\medskip

\noindent Because the functor $\Omega^*$ is compatible with the tensor product, we can define a $\gG$-Banach algebra $\Omega^*A=\Ind_{\gH}^{\gG}A$ for every $\gH$-Banach algebra $A$.  This defines an (induction) functor from the category of $\gH$-Banach algebras together with the $\gH$-equivariant homomorphisms to the category of $\gG$-Banach algebras with the $\gG$-equivariant homomorphisms. If $\Omega$ is an equivalence, then $\Omega^*=\Ind_{\gH}^{\gG}$ is an equivalence of these categories.

Similar statements are true for Banach modules and equivariant homomorphisms of Banach modules, and for Banach pairs and equivariant homomorphisms of Banach pairs. Note that $\Omega^*$ is \emph{not} defined for linear operators between Banach modules or between Banach pairs. The problem is that $f_{\Omega}^*$ makes sense for linear operators, but the resulting operator between, say, $\rho^*{\gG}$-Banach modules is not necessarily $\Omega \times_{\rho} \Omega$-invariant. So $\rho_!$ of this operator cannot be defined in general.

However, we still get a map on the level of $\KK$-groups because in the intermediate step, we can \emph{make} the operator of the $\KKban$-cycle $\Omega \times_{\rho} \Omega$-invariant (recall that we have assumed $\gG^{(0)}$ to be $\sigma$-compact). This was done in Lemma~\ref{LemmaMakeOperatorsEquivariantShriek}, which enables us to define $\Omega^*$ on the level of $\KKban$-groups.

\begin{definition} Let $\Omega$ be a graph from $\gG$ to $\gH$. Then Theorem~\ref{TheoremPullBackOfKKAlongStrictMorphisms} gives a homomorphism
\[
f_{\Omega}^* \colon \KKbanW{\gH}\left(A,B\right) \to \KKbanW{\rho^*(\gG)}\left(f^*_{\Omega}A,\ f_{\Omega}^*B\right).
\]
Corollary~\ref{CorollaryPShriekForPlainKKTheory} gives us an isomorphism
\[
\rho_! \colon \KK_{\rho^*(\gG)}^{\ban}\left(f_{\Omega}^*A,\ f_{\Omega}^*B\right) \cong \KKbanW{\gG}\left(\Omega^*A,\ \Omega^*B\right).
\]
Define
\[
\Ind_{\gH}^{\gG}:=\Omega^*:= \rho_! \circ f_{\Omega}^* \colon \KKbanW{\gH}\left(A,B\right) \to \KKbanW{\gG}\left(\Omega^*A,\ \Omega^*B\right).
\]
\end{definition}

\noindent A variant of the proof of Proposition~\ref{PropositionCompositionOfGeneralisedPullBackSpaces}, the corresponding statement for Banach spaces, shows:

\begin{proposition}\label{PropositionCompositionOfGeneralisedPullBackCycles} Let $\gK$ be another locally compact Hausdorff groupoid carrying a left Haar system. Let $\Omega$ be a graph from $\gG$ to $\gH$ and $\Omega'$ a graph from $\gH$ to $\gK$. Then
\[
\Omega^*\circ (\Omega')^* = \left(\Omega \times_{\gH} \Omega'\right)^* \colon \KKbanW{\gK}\left(A,B\right) \to \KKbanW{\gG}\left(\Omega^*\Omega'^*A,\ \Omega^*\Omega'^*B\right)
\]
which could also be written as
\[
\Ind_{\gH}^{\gG} \circ \Ind_{\gK}^{\gH} = \Ind_{\gK}^{\gG}\colon \KKbanW{\gK}\left(A,B\right) \to \KKbanW{\gG}\left(\Ind_{\gK}^{\gG}A,\ \Ind_{\gK}^{\gG}A\right)
\]
\end{proposition}

The following proposition is proved in Appendix~D.2 of \cite{Paravicini:07}.

\begin{proposition}\label{PropositionPullbackStrictIsPullbackGraph} Let $f\colon \gG\to \gH$ be a strict morphism. Then
\[
f^* = \Graph(f)^* \colon \KKbanW{\gH}\left(A,B\right) \to \KKbanW{\gG}\left(f^*A,\ f^*B\right)
\]
if we identify $f^*A$ with $\Graph(f)^*A$ and $f^*B$ with $\Graph(f)^*B$ (which is possible according to Proposition~\ref{PropositionPullBackOfMorphismUnambiguousForBanachSpaces}).
\end{proposition}

\begin{corollary}
The homomorphism
\[
\gG^*\colon \KKbanW{\gG}\left(A,\ B\right) \to \KKbanW{\gG}\left(A,\ B\right)
\]
is the identity.
\end{corollary}

\begin{corollary} Let $\Omega$ be an equivalence from $\gG$ to $\gH$. Then
\[
\Omega^*\colon  \KKbanW{\gH}\left(A,B\right) \cong \KKbanW{\gG}\left(\Omega^*A,\ \Omega^*B\right)
\]
is an isomorphism with inverse map $(\Omega^{-1})^*$.
\end{corollary}

\noindent This can also be expressed as
\[
\Ind_{\gH}^{\gG} \colon \KKbanW{\gH}\left(A,B\right) \cong \KKbanW{\gG}\left(\Ind_{\gH}^{\gG}A,\ \Ind_{\gH}^{\gG} B\right)
\]

\subsubsection{$\KKban$-cycles and the linking groupoid}

Let $\Omega$ be an equivalence between $\gG$ and $\gH$. Let $\gL$ denote the linking groupoid as defined in Section~\ref{SubsectionLinkingGroupoid}. The categories of $\gG$-Banach spaces, $\gH$-Banach spaces and $\gL$-Banach spaces are mutually equivalent, and it seems worthwhile to analyse how to turn $\gL$-Banach spaces into $\gH$-Banach spaces:

Let $E$ be a $\gL$-Banach space. The graph of the inclusion $\iota_{\gH}$ of $\gH$ into $\gL$ is isomorphic to $\Omega \sqcup \gH$ as an $\gH$-$\gL$-space. Hence we can identify, using Proposition~\ref{PropositionPullBackOfMorphismUnambiguousForBanachSpaces}, the $\gH$-Banach spaces $\iota_{\gH}^* E= E\restr_{\gH^{(0)}}$ and $\Ind_{\gL}^{\gH} E = (\Omega \sqcup \gH)^*E$; here the restriction $E \restr_{\gH^{(0)}}$ is simply the part of the field $E$ over $\gL^{(0)}$ which is indexed over $\gH^{(0)}$. Hence we can identify induction with restriction in this case.

What we have just said about Banach spaces is also valid for Banach algebras, Banach modules and Banach pairs. Moreover, it also holds for the $\KKban$-groups. More precisely, we have the following result: If $A$ and $B$ are non-degenerate $\gL$-Banach algebras, then $\Ind_{\gL}^{\gH} A$ can be identified with $A\restr_{\gH^{(0)}}$ and $\Ind_{\gL}^{\gH} B$ can be identified with $B\restr_{\gH^{(0)}}$, and the induction homomorphism
\[
\Ind_{\gL}^{\gH}\colon \KKbanW{\gL}\left(A,\  B\right) \cong \KKbanW{\gH} \left(A\restr_{\gH^{(0)}},\ B\restr_{\gH^{(0)}}\right)
\]
can also be identified with the homomorphism given by restriction onto $\gH^{(0)}$ on the level of cycles.

This observation shows that it suffices to understand the restriction onto an open and closed subspace because every induction along a general equivalence can be written as the product of a restriction and the inverse of such a restriction (using the linking groupoid).

\section{The descent}\label{Section:TheDescent}

The descent for unconditional completions and locally compact groupoids was defined in \cite{Lafforgue:06}. We repeat the basic definitions because we decided to slightly modify them on the technical level to make them more systematic (Paragraphs~5.1, 5.3 and 5.4). Moreover, we show that the Banach algebras which we assign to locally compact groupoids using unconditional completions are Morita equivalent for equivalent groupoids (Paragraph~5.3). We also analyse how the descent on the level of $\KKban$-theory behaves under equivalences of groupoids and the pushforward construction
(Paragraphs~5.5 and 5.6).

Let $\gG$ be a locally compact Hausdorff groupoid with unit space $X$. Let $\gG$ carry a Haar system $\lambda$.

\subsection{The convolution product and unconditional completions}

\subsubsection{The convolution product}

Let $E_1$, $E_2$ and $F$ be $\gG$-Banach spaces. Let $\mu\colon E_1 \times_X E_2 \to F$ be a continuous field of bilinear maps (so that $\mu_x\colon (E_1)_x \times (E_2)_x \to F_x$ for all $x\in X=\gG^{(0)}$). We define
\begin{equation}\label{Equation:GeneralConvolution}
\mu(\xi_1,\xi_2) (\gamma'):= \int_{\gG^{r(\gamma')}}\mu_{r(\gamma')} \Big(\,\xi_1(\gamma),\  \gamma \left(\xi_2(\gamma^{-1} \gamma')\right)\,\Big) \rmd \lambda^{r(\gamma')} (\gamma)
\end{equation}
for all $\xi_1 \in \ContSect_c(\gG,r^* E_1)$, $\xi_2 \in \ContSect_c(\gG,r^* E_2)$ and $\gamma'\in \gG$. Then $\mu(\xi_1,\xi_2)$ is in $\ContSect_c(\gG, r^* F)$ and $(\xi_1,\xi_2) \mapsto \mu(\xi_1,\xi_2)$ defines a bilinear map which is separately continuous for the inductive limit topologies and non-degenerate if $\mu$ is non-degenerate.\footnote{The proof of this latter fact is elementary but a bit delicate and can be found in \cite{Paravicini:07}, Subsection~5.1.1.}

If $\mu$ is written as a product, then we simply write $\xi_1 * \xi_2$ for $\mu(\xi_1,\xi_2)$. If $\mu$ is written as a bracket $\langle \cdot,\cdot\rangle$ then we write $\langle \xi_1,\xi_2\rangle$ for $\mu(\xi_1,\xi_2)$.

By direct calculation one can prove that this convolution product respects associativity laws if the involved bilinear maps are $\gG$-equivariant.

\subsubsection{Unconditional completions}

The notion of an unconditional norm for $\Cont_c(\gG)$ was first defined in \cite{Lafforgue:02} for the group case and in \cite{Lafforgue:06} for $\gG$ being a groupoid.

\begin{definition}[Unconditional completion]\label{DefinitionUnconditionalNorm} An \demph{unconditional completion} $\mA(\gG)$ of $\Cont_c(\gG)$ is a Banach algebra containing $\Cont_c(\gG)$ as a dense subalgebra and having the following property
\begin{equation}\label{PropertyUnconditionalNorm}
\forall f_1,f_2\in \Cont_c(\gG):\quad  \big(\forall \gamma\in \gG:\ \abs{f_1(\gamma)} \leq \abs{f_2(\gamma)} \big) \
\Rightarrow \ \norm{f_1}_{\mA(\gG)} \leq \norm{f_2}_{\mA(\gG)}.
\end{equation}
In this case we say that the norm of $\mA(\gG)$ is unconditional. We also write $\norm{\cdot}_{\mA}$ for the norm on $\mA(\gG)$.
\end{definition}

\noindent \emph{For the rest of Section~\ref{Section:TheDescent} we fix an unconditional completion $\mA(\gG)$ of $\Cont_c(\gG)$.}

For technical reasons, we extend the norm $\norm{\cdot}_{\mA}$ as follows: Let $\mF_c^+\left(\gG\right)$ be the set of all non-negative (locally) \demph{bounded} functions $\varphi\colon \gG\to \R$ with compact support. Define
\[
\norm{\varphi}_\mA := \inf\left\{\norm{\psi}_\mA:\ \psi \in \Cont_c(\gG),\ \psi\geq \varphi\right\}
\]
for all $\varphi \in \mF_c^+\left(\gG\right)$. Note that by Property~(\ref{PropertyUnconditionalNorm}) the new semi-norm agrees on $\Cont_c^+(\gG)$ with the norm we started with.

\begin{definition}[The Banach space $\mA(\gG,E)$] Let $E$ be a $\gG$-Banach space. Then we define the following semi-norm on $\ContSect_c(\gG,\, r^*E)$:
\[
\norm{\xi}_{\mA}:= \norm{\gamma \mapsto \norm{\xi(\gamma)}_{E_{r(\gamma)}}}_{\mA}.
\]
The Hausdorff completion of $\ContSect_c(\gG,\, r^*E)$ with respect to this semi-norm will be denoted by $\mA(\gG,E)$.
\end{definition}

\noindent Note that the function $\gamma \mapsto \norm{\xi(\gamma)}$ is not necessarily continuous but has at least compact support and is non-negative upper semi-continuous, so we can apply the extended norm on $\mF_c^+(\gG)$ to it. If $E$ is the trivial bundle over $\gG^{(0)}$ with fibre $E_0$, then $\ContSect_c(\gG,\, r^*E)$ is $\Cont_c(\gG, E_0)$ and $\mA(\gG,E)$ could also be denoted as $\mA(\gG,E_0)$; in particular, if $E_0 =\C$, then $\mA(\gG,E) = \mA(\gG,\C)= \mA(\gG)$.

\begin{defprop}\label{DefpropUnconditionalDescentContinuityofConvolution} Let $E_1$, $E_2$, $F$ be $\gG$-Banach spaces and let $\mu\colon E_1  \times_X E_2 \to F$ be a bounded continuous field of bilinear maps. Then for all $\xi_1\in \ContSect_c(\gG,\, r^* E_1)$ and $\xi_2\in \ContSect_c(\gG,\, r^* E_2)$:
\[
\norm{\mu\left(\xi_1,\xi_2\right)}_{\mA(\gG,F)} \leq \norm{\mu}_{\infty}
\norm{\xi_1}_{\mA(\gG,E_1)}\norm{\xi_2}_{\mA(\gG,E_2)}.
\]
So $\mu$ lifts to a continuous bilinear map $\mA(\gG,\mu)$ from $\mA(\gG,E_1) \times \mA(\gG,E_2)$ to $\mA(\gG,F)$
(with norm less than or equal to $\norm{\mu}_\infty$). If $\mu$ is non-degenerate, then so is $\mA(\gG,\mu)$.
\end{defprop}

\noindent The non-degeneracy result can be deduced from the corresponding result on the level of sections with compact support: The canonical map from $\ContSect_c(\gG,\, r^* F)$ to $\mA(\gG, F)$ is continuous with respect to the inductive limit topology on $\ContSect_c(\gG,\,  r^*F)$ and the norm topology on $\mA(\gG,F)$. In particular, if a subset $\Xi$ of $\ContSect_c(\gG,r^*F)$ is dense for the inductive limit topology, then its canonical image in $\mA(\gG,F)$ is dense for the norm topology.

Note that the convolution product between unconditional completions preserves all kinds of associativity laws for bilinear between the underlying $\gG$-Banach spaces. In particular, we have the following result:

\begin{proposition}
Let $B$ be a $\gG$-Banach algebra.

\begin{enumerate}
\item The space $\ContSect_c(\gG,r^*B)$ is an associative algebra with the convolution product
\[
(\xi_1 *\xi_2) (\gamma'):= \int_{\gG^{r(\gamma')}} \xi_1(\gamma)\ \gamma \xi_2(\gamma^{-1} \gamma') \rmd \lambda^{r(\gamma')} (\gamma)
\]
for all $\gamma'\in \gG$, $\xi_1,\xi_2\in \ContSect_c(\gG,r^*B)$. The completion $\mA(\gG, B)$ is a Banach algebra which is non-degenerate if $B$ is non-degenerate.

\item Let $E$ be a right $\gG$-Banach $B$-module. Then $\mA(\gG,E)$ is a right Banach $\mA(\gG,B)$-module which is non-degenerate if $E$ is non-degenerate. An analogous statement is true for left modules.

\item Let $E$ be a $\gG$-Banach $B$-pair. The bracket of $E$ induces a bilinear map
\[
\langle\cdot ,\cdot \rangle_{\mA(\gG,B)}\colon \mA(\gG, E^<) \times \mA(\gG, E^>) \to \mA(\gG, B)
\]
that turns $(\mA(\gG,E^<), \mA(\gG,E^>))$ into a Banach $\mA(\gG,B)$ which we call $\mA(\gG,E)$. The bracket of $\mA(\gG,E)$ is non-degenerate if the bracket of $E$ is non-degenerate.
\end{enumerate}
\end{proposition}

\subsection{Morita equivalence of unconditional completions}\label{Subsection:MoritaEquivalenceOfUnconditionalCompletions}

It is well-known that the C$^*$-algebras of equivalent groupoids are Morita equivalent, see \cite{MuhReWill:87}. We are now going to transfer this result to unconditional completions and we are also going to allow for coefficients in Banach algebras, using the following concept of Morita equivalence for Banach algebras introduced by Vincent Lafforgue in an unpublished note: A \demph{Morita equivalence} between Banach algebras $A$ and $B$ is a pair $E=(E^<,E^>)$, where $E^<$ and $E^>$ are Banach spaces, such that the direct sum $L:=A\oplus B \oplus E^< \oplus E^>$ forms a Banach algebra with a multiplication which operates as the multiplication of two-by-two matrices if one writes the direct sum as
\[
L=\left(\begin{array}{cc} A& E^>\\ E^< & B\end{array}\right)
\]
i.e., there are binary operations $A \times E^< \to E^<$,\ \ $E^< \times E^> \to B$ etc.\ which satisfy a number of associativity and norm conditions; in addition, the binary operations are assumed to be non-degenerate, i.e., the closed linear span of $AE^<$ is all of $E^<$ etc. A more rigorous definition can be found in \cite{Paravicini:07:Morita:erschienen}.

Note that the Banach algebra $A$ above is not only Morita equivalent to $B$, but also to the ``linking algebra'' $L$. The point is that the inclusion of $A$ into $L$ identifies $A$ with a \emph{full corner} of $L$: There is an idempotent $P$ in $\Mult(L)$ such that $A=PLP$ and the closed linear span of $LPL$ is all of $L$. The Morita equivalence between $A=PLP$ and $L$ is given by $(LP, PL)$.

The same pattern can be observed in the case of groupoids as we have seen in Paragraph~\ref{SubsectionLinkingGroupoid}. Working with the linking groupoid we can hence reduce the treatment of general equivalences to the case of inclusion as a subgroupoid.

Let $U$ be an open subset of the unit space $\gG^{(0)}$ of the locally compact Hausdorff groupoid $\gG$. Then $\gG^{U}_U$ is an open subgroupoid of $\gG$, inheriting a Haar system from $\gG$, and the unconditional norm $\norm{\cdot}_{\mA}$ on $\Cont_c(\gG)$ restricts to an unconditional norm on $\Cont_c(\gG^U_U)$; the resulting unconditional completion will be called $\mA(\gG^U_U)$. If $E$ is a $\gG$-Banach space, then we can consider $E\restr_U = (E_u)_{u\in U}$ which is a $\gG^{U}_U$-Banach space in a canonical way. We can imbed $\Cont_c(\gG_U^U,r^* E\restr_U)$ canonically into $\Cont_c(\gG,r^*E)$, and hence $\mA(\gG^U_U, E\restr_U)$ is a closed subspace of $\mA(\gG, E)$. This applies in particular to the case that $E=B$ is a $\gG$-Banach algebra; then $\mA(\gG^U_U, B\restr_U)$ is a closed subalgebra of $\mA(\gG, B)$. 

Recall from Paragraph~\ref{SubsectionLinkingGroupoid} that a subset $U\subseteq \gG^{(0)}$ is called \demph{full} if $\gG_U \circ \gG^U = \gG$.

\begin{proposition}\label{Proposition:InclusionDescentMoritaEquivalence} Let $B$ be a non-degenerate $\gG$-Banach algebra and let $U$ be an open, closed and \emph{full} subspace of $\gG^{(0)}$. Then there is a projection $P$ in the multiplier algebra $\Mult(\mA(\gG,B))$ such that $\mA(\gG_U^U, B\restr_U)=P \mA(\gG,B) P$ and such that the span of $\mA(\gG,B) P\mA(\gG,B)$ is dense in $\mA(\gG,B)$, i.e., $\mA(\gG_U^U, B\restr_U)$ is a \emph{full corner} in $\mA(\gG,B)$. In particular, $\mA(\gG_U^U, B\restr_U)$ and $\mA(\gG,B)$ are \emph{Morita equivalent}.
\end{proposition}
\begin{proof}
Define continuous linear maps
\[
p^> \colon \ContSect_c(\gG,\, r^*B) \to \ContSect_c(\gG,\, r^*B), \ \xi \mapsto \xi \restr_{\gG_U}
\]
and
\[
p^< \colon \ContSect_c(\gG,\, r^*B) \to \ContSect_c(\gG,\, r^*B), \ \xi \mapsto \xi \restr_{\gG^U},
\]
where the restricted sections should be extended by zero to all of $\gG$. Then $(p^<)^2=p^<$ and $(p^>)^2=p^>$. Moreover, $p^<$ is $\ContSect_c(\gG,\, r^*B)$-linear on the right, $p^>$ is $\ContSect_c(\gG,\, r^*B)$-linear on the left. Finally, for all $\xi_1, \xi_2\in \ContSect_c(\gG,\, r^*B)$:
\[
p^>(\xi_1) * \xi_2 = \xi_1 * p^<(\xi_2).
\]
So $p=(p^<,p^>)$ could be called an (idempotent) multiplier of $\ContSect_c(\gG,\, r^*B)$. We have
\[
p \ContSect_c(\gG,\, r^*B) p = \ContSect_c(\gG_U^U,\, r^*B).
\]
Note that the maps $p^<$ and $p^>$ are contractive on the level of sections with compact support because $\mA(\gG)$ is unconditional. Hence $p^<$ and $p^>$ give contractive operators $P^<$ and $P^>$ on $\mA(\gG,B)$ such that $\iota \circ p^< = P^< \circ \iota$ and $\iota \circ p^> = P^> \circ \iota$ where $\iota$ denotes the inclusion of $\Cont_c(\gG,r^*B)$ into $\mA(\gG,B)$. The pair $P=(P^<,P^>)$ inherits the algebraic properties of the $p$, so $P$ is an idempotent multiplier of $\mA(\gG,B)$ and we have
\[
P \mA(\gG, B) P = \mA(\gG^U_U, B\restr_U).
\]
What is left to show is that $P$ is a full projection, i.e., that the span of $\mA(\gG,B) P\mA(\gG,B)$ is dense in $\mA(\gG,B)$. This is a consequence of the fact that the span of $\ContSect_c(\gG, r^*B) p \ContSect_c(\gG,r^*B)$ is dense in $\ContSect_c(\gG,r^*B)$ for the inductive limit topology. To see this, first observe that $\ContSect_c(\gG, r^*B)p$ is the same as $\ContSect_c(\gG_U, r^*B)$ and $p\ContSect_c(\gG, r^*B)$ is the same as $\ContSect_c(\gG^U, r^*B)$. We thus have to show that the span of $\ContSect_c(\gG_U, r^*B)*\ContSect_c(\gG^U, r^*B)$ is dense in $\ContSect_c(\gG, r^*B)$. This follows because $U$ is a full subset and $B$ is non-degenerate; it can be proved in much that same way as the fact that $\ContSect_c(\gG, r^*B)$ is a non-degenerate algebra. See \cite{Paravicini:07}, Appendix~C.1.2 for the technical details.
\end{proof}

\noindent Now consider a locally compact Hausdorff groupoid $\gH$ equivalent to $\gG$ and let it also carry a Haar system. Let $B$ be a non-degenerate $\gG$-Banach algebra. We want to show that $\mA(\gG, B)$ and $\mA(\gH, \Ind_{\gG}^{\gH}B)$ are Morita equivalent. For this to make sense we first have to say what $\mA(\gH)$ is, and we do this by considering the linking groupoid of $\gG$ and $\gH$. Let $\gL$ be the linking groupoid for the equivalence of $\gG$ and $\gH$ and let $\mA(\gL)$ be an unconditional completion of $\Cont_c(\gL)$. Then $\mA(\gL)$ induces unconditional completions $\mA(\gG)$ of $\Cont_c(\gG)$ and $\mA(\gH)$ of $\Cont_c(\gH)$.

By Proposition~\ref{Proposition:InclusionDescentMoritaEquivalence}, the inclusions $\mA(\gG, B) \hookrightarrow \mA(\gL,\Ind_{\gG}^{\gL} B)$ and $\mA(\gH, \Ind_{\gG}^{\gH} B) \hookrightarrow \mA(\gL,\Ind_{\gG}^{\gL} B)$ induce Morita equivalences
\[
\mA(\gG, B) \simMorita \mA(\gL,\Ind_{\gG}^{\gL} B) \simMorita \mA(\gH, \Ind_{\gG}^{\gH} B).
\]
Note that there is also a canonical equivalence directly between $\mA(\gG, B)$ and $\mA(\gH, \Ind_{\gG}^{\gH} B)$ as these algebras are contained as full corners in $\mA(\gL, \Ind_{\gG}^{\gL} B)$, see Section~5.2 of \cite{Paravicini:07:Morita:erschienen}. 

We state the main result that we have shown in this section for further reference:
\begin{theorem}\label{Theorem:EquivalentGroupoidsEquivalentAlgebras}
Let $\gG$ and $\gH$ be equivalent locally compact Hausdorff groupoids with Haar measures and let $\mA(\gG)$ and $\mA(\gH)$ be as above. Then we have a Morita equivalence
\[
\mA(\gG, B) \simMorita \mA(\gH, \Ind_{\gG}^{\gH} B).
\]
In particular, we have
\[
\Leb^1(\gG, B) \simMorita \Leb^1(\gH, \Ind_{\gG}^{\gH} B).
\]
\end{theorem}

\noindent By a result of V.~Lafforgue, the $\KTh$-theory of Morita equivalent Banach algebras is Morita equivalent; compare \cite{Paravicini:07:Morita:erschienen}, where this is shown for $\KTh_0$; note that the corresponding statement for $\KTh_1$ follows straight, because it is clear that if $E$ is a Morita equivalence between two Banach algebras $A$ and $A'$, then $\Cont_0(\R,E)$ defines a Morita equivalence between $\Cont_0(\R, A)$ and $\Cont_0(\R,A')$. Hence we have:

\begin{corollary}
In the situation of the preceding theorem we have
\[
\KTh_*(\mA(\gG, B)) \cong \KTh_*(\mA(\gH, \Ind_{\gG}^{\gH} B)).
\]
In particular, we have
\[
\KTh_*(\Leb^1(\gG, B)) \cong \KTh_*(\Leb^1(\gH, \Ind_{\gG}^{\gH} B)).
\]
\end{corollary}

\subsection{The descent and linear maps}

Let $E$ and $F$ be $\gG$-Banach spaces and let $T$ be a bounded continuous field of linear maps between them. We are now constructing linear maps between $\mA(\gG,E)$ to $\mA(\gG,F)$; there are two different ways to do this depending on whether the operator is thought of acting on the left or on the right.

\begin{defprop}
\begin{enumerate}
\item The continuous linear operator $\xi \mapsto T * \xi$ from $\mA(\gG,E)$ to $\mA(\gG,F)$ is defined by the assignment
\[
\ContSect_c(\gG,r^*E)\ni \xi \quad \mapsto \quad (\gamma \mapsto T_{r(\gamma)} \xi(\gamma))\in \ContSect_c(\gG,r^*F)
\]
and satisfies $\norm{T*\xi}_{\mA(\gG,F)}\leq \norm{T}_\infty \norm{\xi}_{\mA(\gG,E)}$ for all $\xi\in \ContSect_c(\gG,E)$.

\item The continuous linear operator $\xi \mapsto T * \xi$ from $\mA(\gG,E)$ to $\mA(\gG,F)$ is defined by the assignment
\[
\ContSect_c(\gG,r^*E)\ni \xi \quad \mapsto \gamma \quad \left[T_{s(\gamma)} (\gamma^{-1} \xi(\gamma))\right] \in \ContSect_c(\gG,r^*F)
\]
and satisfies $\norm{\xi*T}_{\mA(\gG,F)}\leq \norm{T}_\infty \norm{\xi}_{\mA(\gG,E)}$ for all $\xi\in \ContSect_c(\gG,E)$.

\item If $T$ is $\gG$-equivariant, then $T*\cdot = \cdot *T$ and the map $T\mapsto T*\cdot = \cdot*T$ makes $E\mapsto \mA(\gG,E)$ a functor from the $\gG$-Banach spaces to the Banach spaces.
\end{enumerate}
\end{defprop}

As might be expected, the two ways of forming the descent of a linear map are compatible with the general convolution product; we do not formulate a general rule for this compatibility but state the consequences which are of interest in what follows:

\begin{proposition}
\begin{enumerate}
\item If $B$ and $B'$ are $\gG$-Banach algebras and $\varphi$ denotes a $\gG$-equivariant field of homomorphisms between them, then $\mA(\gG, \varphi):=\varphi*\cdot = \cdot*\varphi$ is a continuous homomorphism from $\mA(\gG, B)$ to $\mA(\gG,B')$. Similarly, $\gG$-equivariant homomorphisms with coefficient maps between $\gG$-Banach modules or $\gG$-Banach pairs descend to continuous homomorphisms.

\item Let $B$ be a $\gG$-Banach algebra and let $E$ and $F$ be $\gG$-Banach $B$-pairs. Let $T=(T^<,T^>)$ be an element of $\ \Lin_B(E,F)$. Then
\begin{enumerate}
\item $T^> *\cdot$ is a linear operator from $\mA(\gG, E^>)$ to $\mA(\gG, F^>)$ being $\mA(\gG, B)$-linear on the right and of norm $\norm{T^>*\cdot}\leq \norm{T^>}$;
\item $\cdot*T^<$ is a linear operator from $\mA(\gG, F^<)$ to $\mA(\gG, E^<)$ being $\mA(\gG, B)$-linear on the left and of norm $\norm{\cdot*T^>} \leq \norm{T^<}$;
\item The pair $(\cdot * T^<,\ T^>*\cdot)$ is in $\Lin_{\mA(\gG,B)} (\mA(\gG, E),\ \mA(\gG, F))$ and of norm less than or equal to $\norm{T}$. It will be denoted by $\mA(\gG,\ T)$.
\end{enumerate}
The assignment $E\mapsto \mA(\gG,E)$ and $T\mapsto \mA(\gG,T)$ defines a functor from the category of $\gG$-Banach $B$-pairs to the category of Banach $\mA(\gG, B)$-pairs.
\end{enumerate}
\end{proposition}

\subsection{The descent and $\KKban$-cycles}\label{SubsectionTheDescentAndKasparovCycles}

Let $A$ and $B$ be $\gG$-Banach algebras. If $E$ is a $\gG$-Banach $A$-$B$-pair, then there is\footnote{Compare Proposition 1.3.3 of \cite{Lafforgue:06}.} a canonical action of the Banach algebra $\mA(\gG,A)$ on the Banach $\mA(\gG,B)$-pair $\mA(\gG,E)$.

\begin{defprop}\label{DefPropDescentRespectsKasparovCycles}\footnote{Compare D\'{e}finition-Proposition 1.3.4 of \cite{Lafforgue:06}.} Let $(E,T) \in \EbanW{\gG}(A,B)$. Then define
\[
j_{\mA}(E,T):= \left(\mA(\gG,E),\ \mA(\gG,T)\right) \in \Eban\left(\mA(\gG,A),\ \mA(\gG,B)\right).
\]
\end{defprop}

\noindent The idea of the proof, as given in \cite{Lafforgue:06}, is to express operators of the form $[\alpha, \mA(\gG,T)]$ and $\alpha(\mA(\gG,, T)^2-1)$, where $\alpha \in \ContSect_c(\gG, r^*A)$, as the convolution with a suitable field of operators; by this we mean the following: If $S =(S^<,S^>) = (S_{\gamma})_{\gamma\in \gG}$ is in $\Lin_{r^*B} (r^*E)$ and has compact support, then there is a canonical definition of a convolution operator $\hat{S} = (\cdot *S^<, S^>* \cdot)$ in $\Lin_{\mA(\gG,B)}(\mA(\gG,E))$. If $S$ is a \emph{compact} operator with compact support, then one can show that also $\hat{S}$ is compact. For details, see \cite{Lafforgue:06}, Section~1.3, or \cite{Paravicini:07}, Section~5.2.7 and Section~5.2.8.

One can show that the descent map respects homotopies and direct sums of cycles, the pullback and the pushforward along homomorphisms of Banach algebras (at least up to homotopy). The descent therefore induces a homomorphism on the level of $\KKban$-groups:

\begin{theorem}\label{Theorem:TheDescentIsAHomomorphism}
Let $A$ and $B$ be $\gG$-Banach algebras and $\mA(\gG)$ an unconditional completion of $\Cont_c(\gG)$. Then $j_{\mA}$ induces a group homomorphism
\[
j_{\mA}\colon \KKbanW{\gG}\left(A,B\right) \ \to \ \KKban\left(\mA(\gG, A),\ \mA(\gG, B)\right).
\]
It is natural with respect to $\gG$-equivariant homomorphisms in both variables.
\end{theorem}

\subsection{The descent and equivalences of groupoids}

Let $\gG$ and $\gH$ be equivalent locally compact Hausdorff groupoids equipped with a Haar system and let $\Omega$ be an equivalence. Assume that $X:=\gG^{(0)}$ and $Y:=\gH^{(0)}$ are $\sigma$-compact (we need this condition only for the results concerning $\KK$-theory). Let $\gL$ be the linking groupoid for the equivalence of $\gG$ and $\gH$. As above, let $\mA(\gL)$ be an unconditional completion of $\Cont_c(\gL)$, inducing unconditional completions $\mA(\gG)$ and $\mA(\gH)$. Let $A$ and $B$ be non-degenerate $\gG$-Banach algebras.

We want to show that the following diagram is commutative:
\[
\xymatrix{
\KKbanW{\gG}\left(A,B\right) \ar[r]\ar[d]_{\cong}& \KKban\left(\mA(\gG,A),\ \mA(\gG,B)\right)\ar[d]^{\cong}\\
\KKbanW{\gH} \left(\Ind A, \Ind B\right) \ar[r]& \KKban\left(\mA(\gH, \Ind A), \mA(\gH, \Ind B)\right)
}
\]
Of course, one has to say what the vertical arrow on the right hand side should be. This is going to cause some problems, but we are going to find a commutative diagram which is at least pretty close to the one above.

Technically, we first show the result for the equivalent groupoids $\gG$ and $\gL$. By symmetry, the result also holds for $\gH$ and $\gL$, and by a transitivity argument we can thus show it to some extend for $\gG$ and $\gH$. We have said above that the categories of $\gL$-Banach algebras and of $\gG$-algebras are equivalent and that the equivalence is given by restriction onto $\gG^{(0)}=X$. For notational convenience, instead of starting with $\gG$-algebras, we prefer to consider $\gL$-Banach algebras. So let $A$ and $B$ be non-degenerate $\gL$-Banach algebras. Then $B\restr_X$ is the non-degenerate $\gG$-Banach algebra corresponding to the $\gL$-Banach algebra $B$. We have a canonical inclusion of $\mA(\gG, B):=\mA(\gG,\ B\restr_X)$ into $\mA(\gL,B)$ as a closed subalgebra via some inclusion map $\iota$ (actually, as a full corner as we have shown above).

We consider the diagram
\begin{equation}\label{Diagram:MoritaDescentSquare}
\xymatrix{
\KKbanW{\gL}\left(A,B\right) \ar[r]\ar[d]_{\cong}& \KKban\left(\mA(\gL,A),\ \mA(\gL,B)\right)\ar[d]\\
\KKbanW{\gG}\left(A\restr_X, B\restr_X\right) \ar[r]& \KKban\left(\mA(\gG, A), \mA(\gG, B)\right)
}
\end{equation}
The left vertical arrow is given by the induction isomorphism which happens to be restriction onto $X$, the horizontal arrows are the descent homomorphisms. The right vertical arrow is given by the pullback along the inclusion of $\mA(\gG,A)$ into $\mA(\gL,A)$ in the first component and the pushforward along the Morita equivalence of $\mA(\gL,B)$ and $\mA(\gG,B)$ in the second component (compare Theorem~5.29 of \cite{Paravicini:07:Morita:erschienen}).  We now show the following theorem:

\begin{theorem} Diagram~(\ref{Diagram:MoritaDescentSquare}) is commutative.
\end{theorem}

\begin{proof}
Let $(E,T)$ be in $\EbanW{\gL}(A,B)$. We have to trace $(E,T)$ through diagram (\ref{Diagram:MoritaDescentSquare}) and prove that the two cycles that we get in the lower right corner are homotopic. We do this by considering the following extended diagram:
\begin{equation}\label{Diagram:MoritaDescentSquare:Extended}
\xymatrix{
\KKbanW{\gL}\left(A,B\right) \ar[r]\ar[dd]_{\cong}& \KKban\left(\mA(\gL,A),\ \mA(\gL,B)\right)\ar[d]\\
&\KKban\left(\mA(\gG,A),\ \mA(\gL,B)\right)\\
\KKbanW{\gG}\left(A\restr_X, B\restr_X\right) \ar[r]& \KKban\left(\mA(\gG, A), \mA(\gG, B)\right)\ar[u]_{\cong}
}
\end{equation}
The upper right vertical arrow is the pullback along the inclusion of $\mA(\gG, A)$ into $\mA(\gL,A)$ in the first component, the lower vertical arrow is the pushforward along the inclusion of $\mA(\gG, B)$ into $\mA(\gL,B)$ in the second component. The latter map is an isomorphism because the two algebras are Morita equivalent, the composition of the upper vertical map and the inverse of the lower vertical map is the right vertical map in Diagram\ (\ref{Diagram:MoritaDescentSquare}).

If we go right and down in Diagram\ (\ref{Diagram:MoritaDescentSquare:Extended}), then we get the cycle $(\mA(\gL,E),\mA(\gL,T))$ where we regard $\mA(\gL,E)$ as a Banach $\mA(\gG,A)$-$\mA(\gL,B)$-pair. If we start with going down, then we get the cycle $(E\restr_X, T\restr_X) \in \EbanW{\gG}(A\restr_X,B\restr_X)$. If we go down and right, then we are left with the cycle $(\mA(\gG,E\restr_X),\ \mA(\gG,T\restr_X))$ regarded as a Banach $\mA(\gG,A)$-$\mA(\gG,B)$-pair. Finally, if we go down-right-up, then we get the cycle $(\mA(\gG,E\restr_X)\otimes_{\mA(\gG,B)}\mA(\gL,B),\ \mA(\gG,T\restr_X) \otimes 1)$ -- here we suppress the unitalisations that appear in the definition of the pushforward as they can be neglected up to homotopy.

We now define a homomorphism $\Phi$ from $\mA(\gG,E\restr_X)\otimes_{\mA(\gG,B)} \mA(\gL,B)$ to $\mA(\gL,E)$ with coefficient maps $\id_{\mA(\gG,A)}$ and $\id_{\mA(\gL,B)}$ by
\begin{eqnarray*}
\Phi^>\colon \mA(\gG,\,E^>\restr_X)\otimes_{\mA(\gG,B)}\mA(\gL,B)&\to&  \mA(\gL, E^>),\\
\xi^>\otimes \beta &\mapsto & \xi^>* \beta
\end{eqnarray*}
where we regard $\xi^> \in \mA(\gG,\, E^>\restr_X)$ as an element of $\mA(\gL,E^>)$; define $\Phi^<$ similarly. By the associativity of the convolution, the pair $\Phi:=(\Phi^<,\Phi^>)$ is a concurrent homomorphism. We now show that it induces a homotopy using Theorem~3.1 of \cite{Paravicini:07:Morita:erschienen}: 
Let $\alpha \in \ContSect_c(\gG,r^*A)$ and $\varepsilon>0$. Then $[\alpha,\, \mA(\gL,T)]$ is given by convolution with the compact continuous field of operators with compact support (compare the discussion after Theorem~\ref{Theorem:TheDescentIsAHomomorphism}):
\[
\gL\ni \gamma \mapsto \alpha(\gamma) \gamma T_{s(\gamma)} - T_{r(\gamma)} \alpha(\gamma) \in \Komp_{r^* B} \left(r^*E\right)_c.
\]
The support of this field is actually contained in $\gG$ because $\alpha$ is supported in $\gG$.

By Lemma~\ref{LemmaConvolutionWithCompactOperatorsOpenSubgroupoids} below, we can approximate $[\alpha,\, \mA(\gL,T)]$ by sums of operators of the form $\ketbra{\eta^>}{\xi^<}$ with $\xi^>\in \ContSect_c(\gL,\, r^*E^>)$ and $\xi^<\in \ContSect_c(\gL,\, r^*E^<)$, both having their support in $\gG$. Because $\mA(\gG, E^>)$ is a non-degenerate right Banach $\mA(\gG,B)$-module and $\mA(\gG, E^<)$ is a non-degenerate left Banach $\mA(\gG,B)$-module, we can actually approximate $[\alpha,\ \mA(\gL,T)]$ as follows: We can find an $n\in \N$ and $\xi^<_1,\ldots,\xi^<_n \in \ContSect_c(\gL,\, r^*E^<)$, $\xi^>_1,\ldots,\xi^>_n \in \ContSect_c(\gL,\, r^*E^>)$ and $\beta^<_1,\ldots,\beta^<_n,\beta^>_1,\ldots,\beta^>_n \in \ContSect_c(\gL,\, r^*B)$ which all are supported in $\gG$ such that
\[
\norm{[\alpha,\ \mA(\gL,T)] - \sum_{i=1}^n \ketbra{\xi^>_i * \beta^>_i}{\beta^<_i*\xi^<_i}}\leq \varepsilon.
\]

\noindent Note that we can regard the $\xi^>_i$ and the $\xi^<_i$ also as sections living on $\gG$. If we do so, we have $\xi^>_i * \beta^>_i = \Phi^>(\xi^>_i \otimes \beta_i^>)$ and $\beta^<_i * \xi^<_i = \Phi^<(\beta^<_i \otimes \xi_i^<)$ for all $i\in \{1,\ldots,n\}$.

The operator $[\alpha,\, \mA(\gL,T)] - \sum_{i=1}^n \ketbra{\xi^>_i * \beta^>_i}{\beta^<_i*\xi^<_i}$ leaves the subspace $\mA(\gG, E\restr_X)$ invariant. The norm of the restricted operator is of course less than or equal to the norm of the operator itself.

Note that $\ketbra{\xi^>_i \otimes \beta_i^>}{\beta_i^<\otimes \xi_i^<} = \ketbra{\xi^>_i *\beta_i^>} {\beta_i^<*\xi_i^<} \otimes 1$ and hence
\begin{eqnarray*}
&&\norm{\left[\alpha \otimes 1,\, \mA\left(\gG,\, T\restr_X\right) \otimes 1\right] - \sum_{i=1}^n \ketbra{\xi^>_i \otimes \beta_i^>}{\beta_i^<\otimes \xi_i^<}}\\
&=& \norm{\left(\left[\alpha,\, \mA\left(\gG,\, T\restr_X\right) \right] - \sum_{i=1}^n \ketbra{\xi^>_i *\beta_i^>}{\beta_i^<*\xi_i^<}\right)\otimes 1}\\
&\leq & \norm{\left[\alpha,\, \mA\left(\gG,\, T\restr_X\right) \right] - \sum_{i=1}^n \ketbra{\xi^>_i *\beta_i^>}{\beta_i^<*\xi_i^<}} \leq \varepsilon.
\end{eqnarray*}

\noindent A similar calculation can be done for $\alpha (\mA(\gL, T)^2-1)$. This shows that $\Phi$ induces a homotopy. Hence the above diagram is commutative.
\end{proof}

\noindent In the preceding proof, we have used the following lemma which is a variant of the result mentioned right after Theorem~\ref{Theorem:TheDescentIsAHomomorphism}.

\begin{lemma}\label{LemmaConvolutionWithCompactOperatorsOpenSubgroupoids}
Let $E$ and $F$ be $\gL$-Banach $B$-pairs. Let $S\in \Komp_{r^*B} (r^*E,\, r^*F)$ have compact support contained in $\gG$. Then the convolution by $S$ as an operator from $\mA(\gL,E)$ to $\mA(\gL,F)$, denoted above by $\hat{S}$, is not only in $\Komp_{\mA(\gL,B)} (\mA(\gL,E),\ \mA(\gL,F))$, but can be approximated by sums of operators of the form $\ketbra{\eta^>}{\xi^<}$ with $\eta^>\in \ContSect_c(\gL,\, r^*F^>)$ and $\xi^<\in \ContSect_c(\gL,\, r^*E^<)$, \emph{both having their support in $\gG$}.
\end{lemma}

\begin{remark} We do not know whether the right vertical arrow in Diagram\ (\ref{Diagram:MoritaDescentSquare}) is an isomorphism because we do not know whether $\KKban$ is also invariant under Morita equivalences in the \emph{first} component. However, the commutativity of the diagram will be sufficient for applications to the Bost conjecture.
\end{remark}

\noindent Let us now go back to the equivalent groupoids $\gG$ and $\gH$. By what we have just said we know that the following diagram commutes:
\begin{equation}\label{Diagram:MoritaDescentSquare:Double}
\xymatrix{
\KKbanW{\gG}\left(A\restr_X, B\restr_X\right) \ar[r]& \KKban\left(\mA(\gG, A\restr_X), \mA(\gG, B\restr_X)\right)\\
\KKbanW{\gL}\left(A,B\right) \ar[r]\ar[d]_{\cong}\ar[u]^{\cong}& \KKban\left(\mA(\gL,A),\ \mA(\gL,B)\right)\ar[d]\ar[u]\\
\KKbanW{\gH}\left(A\restr_Y, B\restr_Y\right) \ar[r]& \KKban\left(\mA(\gH, A\restr_Y), \mA(\gH, B\restr_Y)\right)
}
\end{equation}
The two vertical arrows on the left are given by induction isomorphisms, so their composition from the bottom to the top is also an induction isomorphism, namely the induction by the equivalence $\Omega$.

\section{The Bost conjecture with Banach algebra coefficients}

\subsection{The assembly map and the Bost conjecture}

Let $\gG$ be a locally compact Hausdorff groupoid equipped with a Haar system. Assume\footnote{In \cite{Tu:99:BCGroupoids} it is said that such a space always exists (at least if everything is assumed to be $\sigma$-compact), the given reference \cite{Tu:99} shows this in the case of \'{e}tale metrisable groupoids. We do not venture into the details but content ourselves with the assumption that $\uEgG$ exists.} that there is a locally compact classifying space $\uEgG$ for proper actions of $\gG$, which is then unique up to homotopy. Let $\mA(\gG)$ be an unconditional completion of $\Cont_c(\gG)$. If $B$ is a $\gG$-Banach algebra, then we define a $\gG$-Banach algebra $SB:= B]0,1[$ just as we have defined $B[0,1]$ in Paragraph~\ref{Subsubsection:ThegGBABNullEins}.

\begin{definition}[Topological $\KTh$-theory] For every $\gG$-Banach algebra $B$, define
\[
\KTh^{\top2,\ban}_0\left(\gG,\ B\right):= \lim_{\to}\, \KKbanW{\gG}\left(\Cont_0(X),\ B\right),
\]
where $X$ runs through the closed proper $\gG$-compact subspaces of $\uEgG$. Define $\KTh^{\top2,\ban}_n\left(\gG,\ B\right):= \KTh^{\top2,\ban}_0\left(\gG,\ S^nB\right)$ for $n\in \N$.
\end{definition}

\noindent Note that, if $X$ is a locally compact Hausdorff left $\gG$-space (with anchor map $\rho$), then we would like to think of $\Cont_0(X)$ as a $\gG$-Banach algebra. In particular, we should consider $\Cont_0(X)$ as a field over $\gG^{(0)}$; the fibre of $\Cont_0(X)$ as a $\gG$-Banach algebra over $g\in \gG^{(0)}$ is $\Cont_0(\rho^{-1}(\{g\}))$. A more systematic notation for $\Cont_0(X)$ as a $\gG$-Banach algebra would be $\rho_*(\C_X)$, but we stick to $\Cont_0(X)$ to obtain formulas which look similar to the formulas from the C$^*$-world.

If $B$ is a $\gG$-C$^*$-algebra, then there is a canonical homomorphism from the C$^*$-algebraic version of topological $\KTh$-theory to the Banach algebraic version:
\[
\KTh^{\top2}_*\left(\gG,\ B\right)\to \KTh^{\top2,\ban}_*\left(\gG,\ B\right).
\]

\begin{definition}[Bost assembly map]\label{Definition:TheBostMap} Let $B$ a $\gG$-Banach algebra. Define the homomorphism of abelian groups
\[
\mu_{\mA}^B\colon \KTh^{\top2,\ban}_0\left(\gG,\ B\right) \to \KTh_0\left(\mA\left(\gG,B\right)\right)
\]
to be the direct limit of the group homomorphisms $\mu_{\mA,X}^B$ given by
\[
\KKbanW{\gG}\left(\Cont_0(X),\ B\right) \stackrel{j_{\mA}}{\to}\KKban\left(\mA\left(\gG, \Cont_0(X)\right), \
\mA\left(\gG,B\right)\right) \stackrel{\Sigma(\cdot)\left(\lambda_{X,\gG,\mA}\right)}{\to} \KTh_0\left(\mA\left(\gG, B\right)\right)
\]
where $X$ runs through all closed, $\gG$-compact, proper subspaces of $\uEgG$.
\end{definition}

\noindent We discuss some details of this definition:

\smallskip

\noindent{\bf What is $\lambda_{X,\gG,\mA}$?} If $X$ is a $\gG$-compact proper $\gG$-space, then the element $\lambda_{X,\gG,\mA}$ of $\KTh_0\left(\mA(\gG, \Cont_0(X))\right)$ was defined in \cite{Lafforgue:06}, Paragraph~1.5.2, as follows: Consider the groupoid $\gG \ltimes X$. It is locally compact Hausdorff and proper and satisfies $(\gG \ltimes X)^{(0)} =X$ and $(\gG \ltimes X)\backslash X \cong \gG\backslash X$, this space being compact. We can hence find\footnote{See \cite{Tu:04} for a sufficiently strong existence result for cut-off functions.} a cut-off-function for $\gG \ltimes X$.

A cut-off function for $\gG \ltimes X$ is a function from $X$ to $\R_{\geq 0}$ with compact support such that $\int_{\gG^x} c(\gamma^{-1} x) \rmd \gamma =1$ for all $x\in X$.
\noindent Now consider the function
\[
\gamma \mapsto \left(X_{r_{\gG}(\gamma)} \ni x\mapsto c^{1/2} (x) \ c^{1/2} (\gamma^{-1}x) \right)
\]
with $\gamma \in \gG$. This is an idempotent element of $\ContSect_c\left(\gG,\ r_{\gG}^*\Cont_0(X) \right)$ (actually, we can think of it as an idempotent element of the algebra $\ContSect_c\left(\gG \ltimes X,\ r_{\gG \ltimes X}^* \C_X\right)= \Cont_c\left(\gG \ltimes X\right)$). It therefore gives an idempotent element of $\mA\left(\gG,\ \Cont_0(X)\right)$, and the element of $\KTh_0\left(\mA(\gG,\ \Cont_0(X))\right)$ that it determines is denoted by $\lambda_{X,\gG,\mA}$. One can show that this definition is independent of the choice of the cut-off function $c$.

\noindent {\bf What is $\Sigma(\cdot)\left(\lambda_{X,\gG,\mA}\right)$?} The action $\Sigma$ of $\KKban$ on the $\KTh$-theory was defined in\footnote{See Proposition 1.2.9 of \cite{Lafforgue:02} and the discussion thereafter.} \cite{Lafforgue:02}. In our case, $\Sigma$ is a homomorphism from $\KKban\left(\mA(\gG,\Cont_0(X)),\ \mA(\gG,B)\right)$ to the group of homomorphisms from $\KTh_0\left(\mA(\gG,\Cont_0(X))\right)$ to $\KTh_0\left(\mA(\gG, B)\right)$. Evaluating at $\lambda_{X,\gG,\mA}$ gives a homomorphism from $\KKban\left(\mA(\gG,\Cont_0(X)),\ \mA(\gG,B)\right)$ to $\KTh_0\left(\mA(\gG, B)\right)$.

\medskip

To define the Bost-assembly map also for higher $\KTh$-groups note that there is a the canonical homomorphism $\iota_B\colon \mA(\gG, SB) \to S\mA(\gG,B)$ for every $\gG$-Banach algebra $B$. We can define $\mu_{\mA}^B$ also for $\KTh^{\top2,\ban}_1\left(\gG,\ B\right)$ as the composition
\[
\xymatrix{
\KTh^{\top2,\ban}_1\left(\gG,\ B\right) = \KTh^{\top2,\ban}_0\left(\gG,\ SB\right) \ar[r]^-{\mu_{\mA}^B} & \KTh_0(\mA(\gG,SB)) \ar[r]^-{\iota_{B,*}} & \KTh_0(S\mA(\gG,B)) = \KTh_1(\mA(\gG,B)).
}
\]
Proceed inductively to define the assembly map for all $n\in \N_0$. Note that $\iota_B$ is an isomorphism in $\KTh$-theory as stated in \cite{Lafforgue:06}.

\noindent The \emph{(na\"{\i}ve) Banach algebraic version of the Bost conjecture} for $\gG$ and $\mA(\gG)$ with coefficients in a $\gG$-Banach algebra $B$ asserts that $\mu^B_{\mA}$ is an isomorphism. It is na\"{\i}ve in the sense that the construction of the involved assembly map is analogous to the construction of the Baum-Connes assembly map, but the conjecture itself might not be very helpful because, so far, it is not known how to calculate its left-hand side. There are probably versions of the Bost assembly map with Banach algebra coefficients of the type introduced in \cite{DaLueck:98} which have a more manageable left-hand side, see also \cite{BarEchLueck:07}. For now, we confine ourselves with our ``na\"{i}ve'' version for Banach algebra coefficients; for C$^*$-algebra coefficients there already is a well-established feasible variant: If $B$ is a $\gG$-C$^*$-algebra, then the assembly map introduced in \cite{Lafforgue:06} factors through the Banach algebraic topological $\KTh$-theory, i.e., it is given by the composition
\[
\xymatrix{
\KTh^{\top2}_*(\gG,B) \ar[r] &  \KTh^{\top2, \ban}_*(\gG,B) \ar[r]^{\mu^B_{\mA}} & \KTh_*(\mA(\gG, B)).
}
\]
The \emph{C$^*$-algebraic version of the Bost conjecture} for $\gG$ and $\mA(\gG)$ with coefficients in a $\gG$-C$^*$-algebra $B$ asserts that the composed assembly map is an isomorphism. The name ``Bost conjecture'' is a handy analogue to the name ``Baum-Connes conjecture'' and is justified by the remark on page 10 of \cite{Lafforgue:02} that Jean-Beno\^{i}t Bost has first conjectured a statement of this kind for $\Leb^1$-algebras of closed subgroups of semi-simple Lie groups.\footnote{As remarked by the referee, Bost also considered $\Leb^1$-algebras with coefficients; in \cite{Bost:90}, for example, there is a result which discusses the $\KTh$-theory of $\Leb^1(G,B)$ where $G$ is a group (of a special type) and $B$ is a $G$-Banach algebra.} The C$^*$-algebraic version of the Bost conjecture is an instance of an ``isomorphism conjecture'' as elaborated in \cite{BarEchLueck:07}.

\subsubsection*{The Bost map and varying unconditional completions}

Let $\mA(\gG)$ and $\mB(\gG)$ be unconditional completions of $\Cont_c(\gG)$ such that $\norm{\chi}_{\mB} \leq \norm{\chi}_{\mA}$ for all $\chi\in \Cont_c(\gG)$. Then the following result is a consequence of Proposition~1.4.8 in \cite{Lafforgue:06}, compare also Proposition~1.5.4 of the same article which is the analogous assertion for $\gG$-C$^*$-algebras.

\begin{defprop}\label{DefProp:BostMapAndVaryingCompletions}
Let $B$ be a $\gG$-Banach algebra and let $\iota_{\mA}$ and $\iota_{\mB}$ be the canonical maps from $\ContSect_c(\gG,\ r^* B)$ to $\mA(\gG,B)$ and $\mB(\gG,B)$, respectively. Let $\psi\colon \mA(\gG,B) \to \mB(\gG,B)$ be the homomorphism of Banach algebras such that $\psi \circ \iota_{\mA} = \iota_{\mB}$. Then
\[
\psi_*\colon \KTh_*\left(\mA(\gG,B)\right) \to \KTh_*\left(\mB(\gG,B)\right)
\]
is a homomorphism making the following diagram commutative
\[
\xymatrix{
\KTh^{\top2,\ban}_*\left(\gG,\ B\right) \ar[rr]^{\mu_{\mA}^B} \ar[rrd]_{\mu_{\mB}^B} && \KTh_*\left(\mA\left(\gG,B\right)\right)\ar[d]^{\psi_*}\\
&& \KTh_*\left(\mB\left(\gG,B\right)\right)
}
\]
\end{defprop}

\subsection{Induction and the assembly map}\label{Subsection:InductionAndAssembly}

Let $\gG$ and $\gH$ be locally compact Hausdorff groupoids with Haar systems and $\sigma$-compact unit spaces. Let $\Omega$ be an equivalence between $\gG$ and $\gH$.  Recall from Section~\ref{Subsection:EquivalencesOfGroupoids} that $\Ind_{\gH}^{\gG}$ is an equivalence between the categories of locally compact Hausdorff $\gH$-spaces and of locally compact Hausdorff $\gG$-spaces. It maps proper spaces to proper spaces and a universal locally compact proper $\gH$-space to a universal locally compact proper $\gG$-space. Assume that such universal spaces exist.

We have seen in Section~\ref{Subsection:PullbackAlongGeneralisedMorphisms} that the categories of non-degenerate $\gH$-Banach algebras and of non-degenerate $\gG$-Banach algebras are equivalent via the induction functor $B\mapsto \Ind_{\gH}^{\gG} B$. If $X$ is a locally compact Hausdorff proper $\gH$-space, then we can regard $\Cont_0(X)$ as a non-degenerate $\gH$-Banach algebra. The induced algebra $\Ind_{\gH}^{\gG} \Cont_0(X)$ is then naturally isomorphic to $\Cont_0(\Omega \times_{\gH} X)=\Cont_0(\Ind_{\gH}^{\gG} X)$ as a $\gG$-Banach algebra. If $A$ and $B$ are non-degenerate $\gH$-Banach algebras, then induction gives a canonical isomorphism
\[
\KKbanW{\gH}(A,B) \cong \KKbanW{\gG} (\Ind_{\gH}^{\gG} A,\Ind_{\gH}^{\gG} B).
\]
\noindent As a consequence, we have:

\begin{proposition} Let $B$ be a non-degenerate $\gH$-Banach algebra. Then there is a natural isomorphism given by induction
\[
\KTh^{\top2, \ban}_*(\gH,\ B) \cong \KTh^{\top2,\ban}_* (\gG, \Ind_{\gH}^{\gG} B),
\]
functorial in the groupoid.
\end{proposition}

\noindent Now let $\gL$ denote the linking groupoid (see Paragraph\ \ref{SubsectionLinkingGroupoid}) and let $\mA(\gL)$ be an unconditional completion of $\Cont_c(\gL)$ (this automatically gives unconditional completions $\mA(\gG)$ of $\Cont_c(\gG)$ and $\mA(\gH)$ of $\Cont_c(\gH)$). Let $B$ be a non-degenerate $\gH$-Banach algebra.

\begin{theorem}\label{Theorem:BostAndInduction:Groupoids}
The following square commutes
\[
\xymatrix{
\KTh^{\top2,\ban}_*\left(\gH,\  B\right) \ar[rrr]^{\mu_{\mA}^B} \ar[d]_{\cong} &&& \KTh_*(\mA(\gH,\ B))\ar[d]^{\cong} \\
\KTh^{\top2,\ban}_*\left(\gG,\ \Ind_{\gH}^{\gG} B\right) \ar[rrr]^-{\mu^{\Ind\! B}_{\mA}} &&& \KTh_*(\mA(\gG,\ \Ind_{\gH}^{\gG} B)) \\
}
\]
\end{theorem}

\noindent Before we prove this theorem, we discuss some of its consequences which are valid for a fixed choice of $\mA(\gL)$ and hence of $\mA(\gG)$ and $\mA(\gH)$.

\begin{corollary}
The Banach algebraic version of the Bost conjecture is true for $\gH$ and $B$ if and only if it is true for $\gG$ and $\Ind_{\gH}^{\gG} B$.
\end{corollary}

\begin{corollary}
The Banach algebraic version of the Bost conjecture with arbitrary coefficients is true for $\gH$ if and only if it is true for $\gG$.
\end{corollary}

\noindent It is not hard to check that the map which connects the C$^*$-algebraic $\KK$-theory with $\KKban$ is compatible with induction. More precisely: If $B$ is an $\gH$-C$^*$-algebra, then $\Ind_{\gH}^{\gG} B$ is a $\gG$-C$^*$-algebra and the following diagram commutes
\[
\xymatrix{
\KTh^{\top2}_*(\gH, B) \ar[d]_{\cong} \ar[r] & \KTh^{\top2,\ban}_*(\gH, B) \ar[rr]^{\mu_{\mA}^{B}} \ar[d]_{\cong} && \KTh_*(\mA(\gH,\ B))\ar[d]^{\cong} \\
\KTh^{\top2}_*\left(\gG, \Ind_{\gH}^{\gG} B\right) \ar[r]   & \KTh^{\top2,\ban}_*\left(\gG,  \Ind_{\gH}^{\gG} B\right) \ar[rr]^-{\mu^{\Ind\! B}_{\mA}} && \KTh_*(\mA(\gG,\ \Ind_{\gH}^{\gG} B)) \\
}
\]
Hence we have the following result for C$^*$-coefficients:

\begin{theorem}\label{Theorem:BostInduction:CStar}
The version for C$^*$-algebra coefficients of the Bost conjecture is true for $\gH$ if and only if it is true for $\gG$.
\end{theorem}

\begin{proof}[Proof of Theorem~\ref{Theorem:BostAndInduction:Groupoids}]
For every locally compact Hausdorff proper $\gH$-compact $\gH$-space $X$ we want to show that the following square commutes
\[
\xymatrix{
\KKbanW{\gH}\left(\Cont_0(X),\ B\right) \ar[rrr]^{\mu_{\mA,X}^B} \ar[d]_{\cong} &&& \KTh_0(\mA(\gH,\ B))\ar[d]^{\cong} \\
\KKbanW{\gG}\left(\Cont_0(\Omega \times_{\gH} X),\ \Ind_{\gH}^{\gG} B\right) \ar[rrr]^-{\mu^{\Ind\! B}_{\mA,\Ind\! X}} &&& \KTh_0(\mA(\gG,\ \Ind_{\gH}^{\gG} B)) \\
}
\]
As $B$ is a non-degenerate $\gH$-Banach algebra, we have $\Ind_{\gH}^{\gG} B \cong \Ind_{\gL}^{\gG} \Ind_{\gH}^{\gL} B$, and it suffices to show that the following diagram commutes:
\[
\xymatrix{
\KKbanW{\gH}\left(\Cont_0(X),\ B\right) \ar[rrr] \ar[d]_{\cong} &&& \KTh_0(\mA(\gH,\ B))\ar[d]^{\cong} \\
\KKbanW{\gL}\left(\Cont_0(\Ind_{\gH}^{\gL} X),\ \Ind_{\gH}^{\gL} B\right) \ar[rrr]  &&& \KTh_0(\mA(\gL,\ \Ind_{\gH}^{\gL} B)) \\
\KKbanW{\gG}\left(\Cont_0(\Omega \times_{\gH} X),\ \Ind_{\gH}^{\gG} B\right) \ar[rrr] \ar[u]^{\cong} &&& \KTh_0(\mA(\gG,\ \Ind_{\gH}^{\gG} B)) \ar[u]_{\cong} \\
}
\]
By symmetry, we can concentrate on one of the involved squares. It hence suffices to show for every non-degenerate $\gL$-Banach $B$ and every locally compact Hausdorff proper $\gL$-compact $\gL$-space $X$ (which we fix for the rest of this section) that the following diagram commutes:
\[
\xymatrix{
\KKbanW{\gH}\left(\Cont_0(X\restr_{\gH^{(0)}}),\ B\restr_{\gH^{(0)}}\right) \ar[r]\ar[d]_{\cong}  &\KKban\left(\mA(\gH,\Cont_0(X\restr_{\gH^{(0)}})),\ \mA(\gH, B\restr_{\gH^{(0)}})\right)\ar[r]\ar[d]& \KTh_0(\mA(\gH,\ B\restr_{\gH^{(0)}})) \ar[d]_{\cong}  \\
\KKbanW{\gL}\left(\Cont_0(X),\ B\right) \ar[r]  &\KKban\left(\mA(\gL,\Cont_0(X)),\mA(\gL,B)\right) \ar[r]& \KTh_0(\mA(\gL,\ B))
}
\]
Recall that the algebra $\mA(\gH, B\restr_{\gH^{(0)}})$ is contained as a full corner in $\mA(\gL, B)$, the inclusion will be denoted by $\varphi_B$; a similar statement is true for the inclusion $\varphi_{\Cont_0(X)} \colon \mA(\gH, \Cont_0(X\restr_{\gH^{(0)}})) \hookrightarrow \mA(\gL, \Cont_0(X))$. The vertical arrow in the middle of the diagram is given by the pullback along $\varphi_{\Cont_0(X)}$ in the first component and the pushforward along the Morita equivalence of $\mA(\gH,B)$ and $\mA(\gL,B)$ in the second component. Note that the pushforward along this Morita equivalence inverts the pushforward along the inclusion $\varphi_{B}\colon \mA(\gH,B) \to \mA(\gL,B)$ in the other direction. Note moreover that we do not know and do not need to know that the middle arrow is an isomorphism.

That the left half of the diagram commutes is a special case of Diagram~(\ref{Diagram:MoritaDescentSquare}). To see that the right half of the diagram commutes, one uses the following observation:

\begin{lemma} Let $c$ be a cut-off function on $X\restr_{\gH^{(0)}}$ with respect to $\gH$. Extend $c$ by zero to a (continuous) function $\tilde{c}$ on the whole of $X$. Then $\tilde{c}$ is a cut-off function on $X$ with respect to $\gL$. If $p$ is the projection in $\mA(\gH, \Cont_0(X\restr_{\gH^{(0)}}))$ given by $c$ and if $\tilde{p}$ is the projection in $\mA(\gL, \Cont_0(X))$ given by $\tilde{c}$, then $\varphi_{\Cont_0(X)}(p) = \tilde{p}$. In particular, we have
\[
\varphi_{\Cont_0(X),*} (\lambda_{X\restr_{\gH^{(0)}},\gH,\mA}) =\varphi_{\Cont_0(X),*}([p]) =[\tilde{p}]= \lambda_{X,\gL,\mA}.
\]
\end{lemma}

\noindent Now let $(E,T) \in \Eban(\mA(\gH, \Cont_0(X\restr_{\gH^{(0)}})), \mA(\gH,B\restr_{\gH^{(0)}}))$. Then we have
\begin{eqnarray*}
&&\varphi_{B,*}(\Sigma([E,T])(\lambda_{X\restr_{\gH^{(0)}},\gH,\mA})) = \varphi_{B,*}(\Sigma([E,T])(\varphi_{\Cont_0(X),*}(\lambda_{X,\gL,\mA})))\\ &=& \varphi_{B,*}(\Sigma(\varphi_{\Cont_0(X)}^*[E,T])(\lambda_{X,\gL,\mA})) =\Sigma(\varphi_{B,*}\varphi_{\Cont_0(X)}^*[E,T])(\lambda_{X,\gL,\mA});
\end{eqnarray*}
here we use the functoriality of the action of $\KKban$ on $\KTh$-theory in both variables (see \cite{Lafforgue:02}, Proposition~1.2.9). The formula means that the right half of the above diagram commutes.
\end{proof}

\end{document}